\RequirePackage{ifpdf}
\ifpdf 
\documentclass[pdftex]{sigma}
\else
\documentclass{sigma}
\fi

\usepackage[all]{xy}

\begin{document}


\renewcommand{\PaperNumber}{099}

\FirstPageHeading

\ShortArticleName{Geometric Structures on Spaces of Weighted Submanifolds}

\ArticleName{Geometric Structures\\ on Spaces of Weighted Submanifolds}

\Author{Brian LEE}

\AuthorNameForHeading{B. Lee}

\Address{Department of Mathematics, University of Toronto, Toronto, ON M5S 2E4 Canada}

\Email{\href{mailto:brianlee@math.toronto.edu}{brianlee@math.toronto.edu}}

\ArticleDates{Received May 31, 2009, in f\/inal form October 25, 2009;  Published online November 02, 2009}

\Abstract{In this paper we use a dif\/feo-geometric framework based on manifolds
that are locally modeled on ``convenient'' vector spaces to study
the geometry of some inf\/inite dimensional spaces. Given a f\/inite dimensional
symplectic manifold $\left(M,\omega\right)$, we construct a~weak
symplectic structure on each leaf ${\textbf{I}}_{w}$ of a foliation
of the space of compact oriented isotropic submanifolds in $M$ equipped
with top degree forms of total measure 1. These forms are called \emph{weightings}
and such manifolds are said to be \emph{weighted}. We show that this
symplectic structure on the particular leaves consisting of weighted
Lagrangian submanifolds is equivalent to a heuristic weak symplectic
structure of Weinstein [\emph{Adv. Math.}
\textbf{82} (1990), 133--159]. When the weightings are
positive, these symplectic spaces are symplectomorphic to reductions
of a~weak symplectic structure of Donaldson [\emph{Asian J. Math.} \textbf{3} (1999), 1--15] on
the space of embeddings of a f\/ixed compact oriented manifold into
$M$. When $M$ is compact, by generalizing a moment map of Weinstein
we construct a symplectomorphism of each leaf ${\textbf{I}}_{w}$
consisting of positive weighted isotropic submanifolds onto a coadjoint
orbit of the group of Hamiltonian symplectomorphisms of $M$ equipped
with the Kirillov--Kostant--Souriau symplectic structure. After def\/ining
notions of Poisson algebras and Poisson manifolds, we prove that each
space ${\textbf{I}}_{w}$ can also be identif\/ied with a~symplectic
leaf of a Poisson structure. Finally, we discuss a kinematic description
of spaces of weighted submanifolds.}

\Keywords{inf\/inite dimensional manifolds; weakly symplectic structures; convenient vector spaces; Lagrangian submanifolds; isodrastic foliation}

\Classification{58B99}

\section{Introduction}

In the same way that f\/inite dimensional manifolds are locally modeled
on $\mathbb{R}^{n}$, many collections of geometric objects can be
viewed as inf\/inite dimensional manifolds locally modeled on interesting
geometric spaces. For example, if $N$ and $M$ are smooth manifolds
then the following local models are known:

\begin{center}
\begin{tabular}{|l||l|}
\hline
\textbf{Collection $\mathcal{M}$}  & \textbf{Modeling Space at $x\in\mathcal{M}$} \tabularnewline
\hline
\hline
dif\/feomorphisms of $N$  & vector f\/ields on $N$\tabularnewline
\hline
\hline
Riemannian metrics on $N$  & symmetric 2-tensors on $N$\tabularnewline
\hline
\hline
smooth maps from $N$ to $M$  & sections of $x^{*}TM$\tabularnewline
\hline
\hline
Lagrangian submanifolds  & closed 1-forms on $x$\tabularnewline
\hline
\end{tabular}
\end{center}

These local models represent certain choices, as many geometric structures
coincident in f\/inite dimensions diverge in inf\/inite dimensions. For
example, there are typically more derivations than equivalence classes
of paths; there are many ways to def\/ine the dual of a tangent space;
there may fail to exist holomorphic charts even when the Nijenhuis
tensor vanishes, etc. Accordingly, there are many frameworks available
to study dif\/ferential geometric structures in inf\/inite dimensions.
Depending on the problem, one might choose to work with Fr\'echet manifolds
\cite{Hamilton}, Dif\/feology \cite{Souriau}, Dif\/ferential Spaces
\cite{Sikorski}, the Global Analysis framework of Fr\"{o}licher, Kriegl
and Michor~\cite{FrolicherKriegl,KrieglMichor}, etc., or just choose
to work heuristically.

Once a framework has been chosen, and a local model identif\/ied, the
geometry of a collection~$\mathcal{M}$ can be explored using the
following correspondence: structures inherent to objects in~$\mathcal{M}$
induce global structures on~$\mathcal{M}$. For example, if~$N$ and~$M$ are as above then
\begin{itemize}\itemsep=0pt
\item the set of Riemannian metrics on $N$ inherits weak Riemannian structures
(Ebin 1970 \cite{Ebin}, Smolentzev 1994 \cite{Smolentsev});
\item if $M$ is symplectic and $L\to M$ is a prequantization line bundle,
then the space of sections~$\Gamma\left(L\right)$ inherits a weak
symplectic structure (Donaldson 2001 \cite{Donaldson2});
\item the set of embeddings of $N$ into $M$ is the total space of a principal
f\/iber bundle with structure group $\text{Dif\/f}\left(N\right)$, the
dif\/feomorphisms of $N$, and base the set of submanifolds of~$M$
dif\/feomorphic to $N$ (Binz, Fischer 1981~\cite{Binz}).
\end{itemize}
In this paper we study a particularly interesting example of this
phenomenon involving Lag\-rangian submanifolds equipped with certain
measures. From the very beginning, we study these objects in the ``Convenient
Setup'' of Fr\"{o}licher, Kriegl, and Michor (see \cite{KrieglMichor}).

The starting point for this framework is the def\/inition of smooth
curves in locally convex spaces called \emph{convenient} vector spaces.
Once the smooth curves have been specif\/ied, smooth maps between convenient
vector spaces can be def\/ined as maps which send smooth curves to smooth
curves. Smooth manifolds then are def\/ined as sets that can be modeled
on convenient vector spaces via charts, whose transition functions
are smooth. Once the appropriate notions of smoothness are specif\/ied,
objects in dif\/ferential geometry are def\/ined by choosing how to generalize
f\/inite dimensional constructions to inf\/inite dimensions (e.g.\ Lie
groups, principal $G$ bundles, vector f\/ields, dif\/ferential forms,
etc.) An important feature of this approach is that the modeling space
$E_{U}$ for each chart $\left(\varphi,U\right)$ can be dif\/ferent
for dif\/ferent chart neighbourhoods~$U$. This dif\/fers from the usual
description of f\/inite dimensional manifolds which are always modeled
on the same vector space $\mathbb{R}^{n}$. This f\/lexibility is useful
in describing the local structure of many inf\/inite dimensional manifolds,
including the collection of Lagrangian submanifolds in a~symplectic
manifold.

In 1990 Alan Weinstein \cite{Weinstein3} introduced a foliation $\mathcal{F}$
of the space of Lagrangian submanifolds in a f\/ixed symplectic manifold
$\left(M,\omega\right)$. A leaf of $\mathcal{F}$ consists of Lagrangian
submanifolds that can be joined by f\/lowing along Hamiltonian vector
f\/ields. $\mathcal{F}$ lifts to a foliation $\mathcal{F}_{w}$ of
the space of pairs $\left(L,\rho\right)$, where $L$ is a Lagrangian
submanifold in $M$ equipped with a smooth density $\rho$ of total
measure 1. Weinstein called such pairs \emph{weighted Lagrangian submanifolds}
and leaves of~$\mathcal{F}$ and~$\mathcal{F}_{w}$ \emph{isodrasts.}
He showed that each leaf ${\textbf{I}}_{w}$ of $\mathcal{F}_{w}$
can be given a weakly nondegenerate symplectic structure $\Omega^{W}$.
He also showed that the leaves consisting of Lagrangian submanifolds
equipped with positive densities can be identif\/ied with coadjoint
orbits of the group of Hamiltonian symplectomorphisms. All of these
constructions were done on a heuristic level.

Instead of starting with the Lagrangian submanifolds directly, we
instead begin by showing that the set of Lagrangian embeddings of
a f\/ixed compact oriented manifold $L_{0}$ into $M$ is the total
space of a principal f\/iber bundle with structure group $\text{Dif\/f}_{+}\left(L_{0}\right)$,
the orientation preserving dif\/feomorphisms of $L_{0}$. The base $\text{Lag}\left(M\right)$
is naturally identif\/ied with the space of oriented Lagrangian submanifolds
in $M$ dif\/feomorphic to $L_{0}$. We def\/ine a foliation $\mathcal{E}$
of the total space which descends to the isodrastic foliation $\mathcal{F}$
of the space of Lagrangian submanifolds. Similarly, the product of
the space of Lagrangian embeddings with the space of top degree forms
on $L_{0}$ that integrate to 1 is the total space of a principal
$\text{Dif\/f}_{+}\left(L_{0}\right)$ bundle. The base of this bundle
can be identif\/ied with the set of pairs $\left(L,\rho\right)$, where
$L$ is an oriented Lagrangian submanifold in $M$ dif\/feomorphic to
$L_{0}$ equipped with a top degree form $\rho$ (not necessarily
non-vanishing) satisfying $\int_{L}\rho=1$. The foliation $\mathcal{E}$
gives a foliation $\mathcal{E}_{w}$ of the total space that descends
to the isodrastic foliation $\mathcal{F}_{w}$ of the base. We def\/ine
a basic 2-form $\Omega$ on the leaves of $\mathcal{E}_{w}$ which
descends to a weakly nondegenerate symplectic structure on the leaves
of $\mathcal{F}_{w}$. We then show that the tangent spaces to the
space of pairs $\left(L,\rho\right)$ in the ``Convenient Setup''
can be identif\/ied with the tangent spaces in Weinstein's heuristic
construction, and that $\Omega$ corresponds to $\Omega^{W}$. In
this way we make rigourous Weinstein's original construction.

The set of pairs $\left(L,\rho\right)$ consisting of Lagrangian submanifolds
equipped with volume forms of total measure 1 is an open subset of
the set of all weighted Lagrangian submanifolds. The leaves of $\mathcal{F}_{w}$
in this open subset of positive weighted Lagrangian submanifolds inherit
the symplectic structure $\Omega$ and provide a link between Weinstein's
symplectic structure and a symplectic structure def\/ined by Simon Donaldson
on the space of smooth mappings between manifolds described brief\/ly
as follows.

In 1999 Donaldson \cite{Donaldson1} heuristically wrote down a symplectic
structure $\Omega^{D}$ on the space of smooth mappings $C^{\infty}\left(S_{0},M\right)$
of a compact oriented manifold $S_{0}$, equipped with a f\/ixed vo\-lu\-me
form $\eta_{0}$, into a symplectic manifold $\left(M,\omega\right)$.
Under some topological restrictions on $\omega$ and $S_{0}$, Donaldson
described a moment map $\mu$ for the $\text{Dif\/f}\left(S_{0},\eta_{0}\right)$-action
of volume preserving dif\/feomorphisms on $C^{\infty}\left(S_{0},M\right)$.
This $\text{Dif\/f}\left(S_{0},\eta_{0}\right)$-action restricts to
a Hamiltonian action on the space of embeddings $\text{Emb}\left(S_{0},M\right)\subset C^{\infty}\left(S_{0},M\right)$,
with respect to the restrictions of~$\Omega^{D}$ and~$\mu$. By a
lemma of Moser, symplectic quotients of $\text{Emb}\left(S_{0},M\right)$
by $\text{Dif\/f}\left(S_{0},\eta_{0}\right)$ can be identif\/ied with
spaces of submanifolds in $M$ equipped with volume forms of f\/ixed
total measure. In fact when $S_{0}$ is half the dimension of $M$
the level surface $\mu^{-1}\left\{ 0\right\} $ consists of Lag\-rangian
embeddings. This suggests that when $\eta_{0}$ has total measure
1 the symplectic quotients of $\left(\text{Emb}\left(S_{0},M\right),\Omega^{D}\right)$
should be related to the leaves of $\mathcal{F}_{w}$ consisting of
positive weighted Lagrangian submanifolds.

The main result of this paper is that reductions of $\left(\text{Emb}\left(S_{0},M\right),\Omega^{D}\right)$
can be def\/ined, in the ``Convenient Setup'', without any topological
restrictions on $\omega$ or $S_{0}$ and that these reductions are
symplectomorphic to leaves of $\mathcal{F}_{w}$ consisting of positive
weighted Lagrangian submanifolds when the dimension of $S_{0}$ is
half the dimension of $M$. In this way we obtain not only a rigorous
formulation of Donaldson's heuristic constructions, but also a precise
relationship between Weinstein's symplectic structure and Donaldson's
symplectic structure. Namely, symplectic quotients of Donaldson's
symplectic space can be identif\/ied with Weinstein's symplectic spaces
in the particular case of leaves consisting of positive weighted Lagrangian
submanifolds.

For $S_{0}$ of dimension less than or equal to half the dimension
of $M$, symplectic reductions of $\left(\text{Emb}\left(S_{0},M\right),\Omega^{D}\right)$
are still well def\/ined in the ``Convenient Setup'' and yield symplectic
spaces consisting of positive weighted isotropic submanifolds in $M$.
This suggests that the symplectic structure $\Omega$ on weighted
Lagrangian submanifolds should have a generalization to weighted isotropic
submanifolds. We show that indeed such a generalization exists, and
that the corresponding symplectic spaces in the particular case of
leaves consisting of positive weighted isotropic submanifolds are
symplectomorphic to reductions of $\left(\text{Emb}\left(S_{0},M\right),\Omega^{D}\right)$.
In this way we obtain a generalization of our observed relationship
between Weinstein's symplectic structure and Donaldson's symplectic
structure to the case of weighted isotropic submanifolds.

Our next result takes its cue from this generalization to weighted
isotropic submanifolds. Namely, we show that the symplectic spaces
of positive weighted isotropic submanifolds are symplectomorphic to
coadjoint orbits of the group $\text{Ham}\left(M\right)$ of Hamiltonian
symplectomorphisms of~$M$ equipped with the Kirillov--Kostant--Souriau
symplectic structure. This symplectomorphism is given by a generalization
of the moment map written down by Weinstein in his identif\/ication
of positive weighted Lagrangian submanifolds with coadjoint orbits
of $\text{Ham}\left(M\right)$. The heuristic idea is that any submanifold
$I$ equipped with a volume form $\rho$ can be viewed as an element
of the dual of the Lie algebra of Hamiltonian vector f\/ields via the
mapping $\left(I,\rho\right)\mapsto\left(f\mapsto\int_{I}f\rho\right)$.
This mapping is equivariant, injective, and hence induces a coadjoint
orbit symplectic structure on spaces of positive weighted submanifolds
that can be joined by Hamiltonian deformations.

These positive weighted isotropic submanifolds have yet another interpretation
akin to leaves of Poisson manifolds in f\/inite dimensions. Given a
f\/inite dimensional Poisson manifold $\left(P,\left\{ \cdot,\cdot\right\} \right)$,
for each smooth function $f\in C^{\infty}\left(P,\mathbb{R}\right)$
on $P$ there exists a unique vector f\/ield $X_{f}$ on $P$ satisfying
$\text{d}g\left(X_{f}\right)=\left\{ f,g\right\} $ for all $g\in C^{\infty}\left(P,\mathbb{R}\right)$.
The leaves swept out by integral curves to such vector f\/ields $X_{f}$
are symplectic manifolds. This picture can be adapted to inf\/inite
dimensions in the following sense. Given an inf\/inite dimensional manifold
$P$, for a subalgebra $\mathcal{A}\subset C^{\infty}\left(P,\mathbb{R}\right)$
we def\/ine a Poisson bracket $\left\{ \cdot,\cdot\right\} $ on $\mathcal{A}$
and a Poisson algebra $\left(\mathcal{A},\left\{ \cdot,\cdot\right\} \right)$
in the usual way. If for every $f\in\mathcal{A}$ there exists a unique
vector f\/ield $X_{f}$ on $P$ satisfying $\text{d}g\left(X_{f}\right)=\left\{ f,g\right\} $
for all $g\in\mathcal{A}$, then the directions swept out by such
vector f\/ields on each point in $P$ def\/ine a distribution on $P$.
We call maximal integral manifolds of this distribution \emph{leaves}.
By def\/ining a Poisson algebra on $\text{Emb}\left(S_{0},M\right)$,
which restricts to a Poisson algebra on the space of isotropic embeddings,
which descends to a Poisson algebra on the space of positive weighted
isotropic submanifolds, we show that the reductions of $\left(\text{Emb}\left(S_{0},M\right),\Omega^{D}\right)$
are symplectic leaves of a~Poisson structure.

As a result we arrive at three dif\/ferent interpretations of the symplectic
spaces consisting of positive weighted isotropic submanifolds. Namely,
they can be identif\/ied with reductions of the space of embeddings
$\left(\text{Emb}\left(S_{0},M\right),\Omega^{D}\right)$, with coadjoint
orbits of the group $\text{Ham}\left(M\right)$ of Hamiltonian symplectomorphisms,
and with symplectic leaves of Poisson structures.

We then take a kinematic approach to the leaves of the foliation $\mathcal{F}$
of the space of Lagrangian submanifolds to obtain a phase space symplectic
structure. That is, by viewing the leaves of $\mathcal{F}$ as possible
conf\/igurations for a submanifold moving in $M$, on each Lagrangian
we can associate ``conjugate momenta'' with top degree forms that
integrate to $0$. We call such pairs $\left(L,\chi\right)$ with
$L$ in a leaf of $\mathcal{F}$ satisfying $\int_{L}\chi=0$ \emph{momentum
weighted Lagrangian submanifolds}. By writing down what should be
the canonical 1-form on this set of momentum weighted Lagrangian submanifolds
and calculating its exterior derivative, we obtain a weakly symplectic
structure.

Finally, we apply this kinematic approach to the set of pseudo Riemannian
metrics of a~f\/ixed signature on a~f\/inite dimensional manifold $N$.
This collection can be viewed as a set of submanifolds by identifying
each metric with its graph as a section. Weightings then can be assigned
to each metric by pulling up a structure assigned to $N$. By equipping
each metric in this way with a compactly supported symmetric 2-tensor
on $N$, we show that the set of all such weighted metrics has a natural
exact symplectic structure.

\subsection*{Conventions}
Unless stated otherwise, all f\/inite dimensional manifolds are smooth,
connected, and paracompact. For manifolds $M$ and $N$ and vector
bundle $E\to M$, we will use the following notation:{\samepage

\smallskip

\begin{tabular}{ll}
$C^{\infty}\left(M,\mathbb{R}\right)$  & smooth functions on $M$;\tabularnewline
$C_{c}^{\infty}\left(M,\mathbb{R}\right)$  & smooth functions of compact support on $M$;\tabularnewline
$\text{Dif\/f}\left(M\right)$  & dif\/feomorphisms of $M$;\tabularnewline
$\text{Dif\/f}_{+}\left(M\right)$  & orientation preserving dif\/feomorphisms of $M$;\tabularnewline
$C^{\infty}\left(M,N\right)$  & smooth mappings from $M$ to $N$; \tabularnewline
$\text{Emb}\left(M,N\right)$  & smooth embeddings from $M$ to $N$; \tabularnewline
$\mathfrak{X}\left(M\right)$  & vector f\/ields on $M$;\tabularnewline
$\Omega^{k}\left(M\right)$  & \emph{k}-forms on $M$;\tabularnewline
$\text{T}_{l}^{k}\left(M\right)$  & $\left(k,l\right)$-tensor f\/ields on $M$;\tabularnewline
$\Gamma\left(E\right)$  & sections of $E\to M$;\tabularnewline
$\imath\left(X\right)$  & interior derivative with respect to $X$;\tabularnewline
$\mathcal{L}_{X}$  & Lie derivative with respect to $X$.\tabularnewline
\end{tabular}
\smallskip

In the absence of summation signs repeated indices are summed over.}


\section{Basic def\/initions}

We begin by describing the ``Convenient Setup'' of Fr\"{o}licher,
Kriegl, and Michor in order to establish what we will mean by smoothness,
tangent vectors, etc. on some inf\/inite dimensional manifolds. Many
def\/initions will be taken verbatim from \cite{KrieglMichor}. All
references like \cite[X.X]{KrieglMichor} in this section refer to
sections in \cite{KrieglMichor}.

\subsection{Locally convex spaces}

A \emph{real topological vector space} $E$ is a vector space equipped
with a topology under which addition $+:E\times E\to E$ and scalar
multiplication $\mathbb{R}\times E\to E$ are continuous.

A subset $C$ in $E$ is said to be
\begin{enumerate}\itemsep=0pt
\item[1)] \emph{circled} if $\lambda C\subset C$ whenever $\left|\lambda\right|\leq1$;
\item[2)] \emph{convex} if $\lambda_{1}C+\lambda_{2}C\subset C$ for all $\lambda_{1},\lambda_{2}\geq0$
satisfying $\lambda_{1}+\lambda_{2}=1$;
\item[3)] \emph{absolutely convex} if $C$ is circled and convex.
\end{enumerate}
A \emph{locally convex space} is a Hausdorf\/f topological vector space
$E$, for which every neighbourhood of $0$ contains an absolutely
convex neighbourhood of $0$.

\subsection{Bounded sets}

A subset $C$ of a locally convex space $E$ is \emph{bounded} if
for each neighbourhood $U$ of $0$ there exists a $\rho>0$ such
that $C\subset\rho U$. The family of all bounded sets in $E$ is
called the \emph{bornology} of $E$. A linear map $T:E\to F$ between
locally convex spaces is \emph{bounded} if it maps bounded sets to
bounded sets \cite[1.1]{KrieglMichor}.

\subsection{Smooth curves}

Let $E$ be a locally convex space. A curve $c:\mathbb{R}\to E$ is
called \emph{differentiable} if the \emph{derivative} $c'\left(t\right):=\underset{h\to0}{\text{lim}}\frac{1}{h}\left[c\left(t+h\right)-c\left(t\right)\right]$
at $t$ exists for all $t$. A curve $c:\mathbb{R}\to E$ is called
\emph{smooth} if all iterated derivatives exist. The set of all smooth
curves in $E$ will be denoted by $C^{\infty}\left(\mathbb{R},E\right)$
\cite[1.2]{KrieglMichor}.

One would hope that reasonable def\/initions of smoothness would imply
that ``dif\/feo\-mor\-phisms'' are homeomorphisms. For this purpose
we will make use of another topology on locally convex spaces.

\subsection[The $c^{\infty}$-topology]{The $\boldsymbol{c^{\infty}}$-topology}\label{sub:cinftyTopology}

The \emph{$c^{\infty}$-topology} on a locally convex space $E$ is
the f\/inest topology for which all smooth curves $c:\mathbb{R}\to E$
are continuous \cite[2.12]{KrieglMichor}. The $c^{\infty}$-topology
is f\/iner than the locally convex topology on $E$ \cite[4.7]{KrieglMichor}.
If $E$ is a Fr\'echet space, (i.e.\ a complete and metrizable locally
convex space), then the two topologies coincide \cite[4.1, 4.11]{KrieglMichor}.

\subsection{Convenient vector spaces}

A \emph{convenient vector space} is a locally convex space $E$ with
the following property: For any $c_{1}\in C^{\infty}\left(\mathbb{R},E\right)$
there is a $c_{2}\in C^{\infty}\left(\mathbb{R},E\right)$ with $c_{2}'=c_{1}$.
Any $c^{\infty}$-closed subspace of a convenient vector space is
convenient \cite[2.12, 2.13, 2.14]{KrieglMichor}.

\subsection{Space of curves}

The set of smooth curves $C^{\infty}\left(\mathbb{R},E\right)$ in
a convenient vector space $E$ has a natural convenient structure.
Moreover, a locally convex space $E$ is convenient if and only if
$C^{\infty}\left(\mathbb{R},E\right)$ is convenient \cite[3.7]{KrieglMichor}.

We would like to study sets that can be locally modeled on convenient
vector spaces. To def\/ine ``smooth transition functions'' we need
to def\/ine smooth mappings between convenient vector spaces.

\subsection*{Convention}

For the rest of this section $E$ and $F$ will denote convenient
vector spaces.

\subsection{Mappings between convenient vector spaces}

Let $U\subset E$ be a $c^{\infty}$-open subset. A mapping $f:U\to F$
is called \emph{smooth} if it maps smooth curves in $U$ to smooth
curves in $F$. Let $C^{\infty}\left(U,F\right)$ denote the set of
all smooth mappings $f:U\to F$, equipped with the f\/inest topology
on $C^{\infty}\left(U,F\right)$ for which all maps $c^{*}:C^{\infty}\left(U,F\right)\to C^{\infty}\left(\mathbb{R},F\right)$,
given by pullback along smooth curves (i.e. $c^{*}\left(f\right)=f\circ c$),
are continuous. Then $C^{\infty}\left(U,F\right)$ is a convenient
vector space \cite[3.11]{KrieglMichor}.

\subsection{Spaces of linear mappings}

Let $L\left(E,F\right)$ denote the set of all bounded linear mappings
from $E$ to $F$. Then $L\left(E,F\right)$ is contained in $C^{\infty}\left(E,F\right)$
\cite[2.11]{KrieglMichor} and inherits a convenient structure \cite[3.17]{KrieglMichor}.
The set of invertible maps in $L(E,F)$ with bounded inverse will
be denoted by $GL(E,F)$.

\subsection{The dif\/ferentiation operator and chain rule}

Let $U\subset E$ be a $c^{\infty}$-open subset. The dif\/ferentiation
operator
\begin{gather*}
\text{d}: \ C^{\infty}(U,F)   \to   C^{\infty}(U,L(E,F)),\qquad
\text{d}f(x)v   :=   \lim_{t\to0}\frac{f(x+tv)-f(x)}{t}
\end{gather*}
exists, is linear and bounded (smooth). Note that the above limit
is taken in the locally convex topology of $F$. Also the chain rule
\[
\text{d}(f\circ g)(x)v=\text{d}f(g(x))\text{d}g(x)v\]
 holds \cite[3.18]{KrieglMichor}.

\subsection{Examples of convenient vector spaces}\label{sub:Examples-of-Convenient}

The following spaces have natural convenient structures:

\begin{tabular}{ll}
$L_{\text{alt}}^{k}\left(E,F\right)$  & bounded alternating multilinear maps $E\times\cdots\times E\to F$ \cite[5.9, 5.13]{KrieglMichor};\tabularnewline
$C^{\infty}\left(M,\mathbb{R}\right)$  & smooth functions on a f\/inite dimensional manifold $M$ \cite[6.1]{KrieglMichor};\tabularnewline
$C_{c}^{\infty}\left(M,\mathbb{R}\right)$  & smooth functions of compact support on a f\/inite\tabularnewline
 & dimensional manifold $M$ \cite[6.2]{KrieglMichor};\tabularnewline
$\Gamma_{c}\left(Q\right)$  & compactly supported smooth sections of a vector bundle $Q\to M$; \tabularnewline
 & with f\/inite dimensional total space, base, and f\/ibers \cite[30.4]{KrieglMichor}.\tabularnewline
\end{tabular}

\subsection{Manifolds}

A \emph{chart} $\left(U,\varphi\right)$ on a set $M$ is a bijection
$\varphi:U\to E_{U}$ from a subset $U\subset M$ onto a $c^{\infty}$-open
set in a convenient vector space $E_{U}$. A family of charts $\left(U_{\alpha},\varphi_{\alpha}\right)_{\alpha\in A}$
is called an \emph{atlas} for $M$, if the~$U_{\alpha}$ cover $M$
and all \emph{transition functions} $\varphi_{\alpha\beta}:=\varphi_{\alpha}\circ\varphi_{\beta}^{-1}:\varphi_{\beta}\left(U_{\alpha}\cap U_{\beta}\right)\to\varphi_{\alpha}\left(U_{\alpha}\cap U_{\beta}\right)$
are smooth. Two atlases are \emph{equivalent} if their union is again
an atlas on $M$. A \emph{smooth manifold}~$M$ is a set together
with an equivalence class of atlases on it \cite[27.1]{KrieglMichor}.

\subsection{Smooth mappings between manifolds}

A mapping $f:M\to N$ between manifolds is \emph{smooth} if for every
$x\in M$ and chart $\left(V,\psi\right)$ on $N$ with $f\left(x\right)\in V$
there exists a chart $\left(U,\varphi\right)$ on $M$ with $x\in U$
and $f\left(U\right)\subset V$ such that $\psi\circ f\circ\varphi^{-1}$
is smooth. So a mapping $f:M\to N$ is smooth if and only if it maps
smooth curves to smooth curves. A smooth mapping $f:M\to N$ is a
\emph{diffeomorphism} if it is a bijection and if its inverse is smooth
\cite[27.2]{KrieglMichor}. The set of smooth maps from $M$ to $N$
will be denoted by $C^{\infty}(M,N)$.

\subsection{Submanifolds}

A subset $N$ of a smooth manifold $M$ is called a \emph{submanifold},
if for each $x\in N$ there is a chart $\left(U,\varphi\right)$ of
$M$ such that $\varphi\left(U\cap N\right)=\varphi\left(U\right)\cap F_{U}$,
where $F_{U}$ is a $c^{\infty}$-closed linear subspace of the convenient
model space $E_{U}$ \cite[27.11]{KrieglMichor}. A curve in a submanifold
$N$ of $M$ is smooth if and only if it is smooth as a curve in $M$.

\subsection{Tangent spaces of a convenient vector space}

Let $a\in E$. A \emph{tangent vector} with base point $a$ is a pair
$\left(a,X\right)$ with $X\in E$. For each neighbourhood $U$ of
$a$ in $E$, a tangent vector $\left(a,X\right)$ def\/ines a derivation
$C^{\infty}\left(U,\mathbb{R}\right)\to\mathbb{R}$ by $X_{a}f:=\text{d}f\left(a\right)\left(X\right)$
\cite[28.1]{KrieglMichor}.

\begin{remark}
In \cite{KrieglMichor} these tangent vectors are called \emph{kinematic
tangent vectors} since they can be realized as derivatives $c'\left(0\right)$
at $0$ of smooth curves $c:\mathbb{R}\to E$. This is to distinguish
them from more general derivations which are called \emph{operational
tangent vectors}.
\end{remark}

\subsection{The tangent bundle}

Let $M$ be a smooth manifold with an atlas $\left(U_{\alpha},\varphi_{\alpha}\right)$
indexed by $\alpha\in A$. On the disjoint union\[
\bigsqcup_{\alpha\in A}\mathcal{U}_{\alpha}\times E_{\alpha}\times\left\{ \alpha\right\} \]
 def\/ine the following equivalence relation:\[
\left(x,v,\alpha\right)\sim\left(y,w,\beta\right)\Longleftrightarrow x=y\text{ and d}\varphi_{\alpha\beta}\left(\varphi_{\beta}\left(y\right)\right)w=v.\]
 A \emph{tangent vector at} $x\in M$ is an equivalence class $\left[\left(x,v,\alpha\right)\right]$.
The quotient $\bigsqcup_{\alpha\in A}\mathcal{U}_{\alpha}\times E_{\alpha}\times\left\{ \alpha\right\} /\sim$
will be called the \emph{tangent bundle of} $M$ and will be denoted
by $TM$.

Let $\pi:TM\to M$ denote the projection $\left[\left(x,v,\alpha\right)\right]\mapsto x$.
$TM$ inherits a smooth manifold structure from $M$. For $x\in M$
the set $T_{x}M:=\pi^{-1}\left(x\right)$ is called the \emph{tangent
space at $x$.} Since each transition function $\varphi_{\alpha\beta}$
is smooth, each dif\/ferential $\text{d}\varphi_{\alpha\beta}(x)$ is
bounded linear, which means each tangent space $T_{x}M$ has a well
def\/ined bornology independent of the choice of chart (cf. \cite[1.1, 2.11]{KrieglMichor}).

Alternatively, we can describe tangent vectors to a smooth manifold
by means of equivalence classes of smooth curves. We will say that
two smooth curves $c_{1}$ and $c_{2}$ in $M$ are equivalent at
$x\in M$, (and write $c_{1}\sim_{x}c_{2}$), if $c_{1}\left(0\right)=x=c_{2}\left(0\right)$
and $\left.\frac{d}{dt}\right|_{t=0}\varphi_{\alpha}\circ c_{1}\left(t\right)=\left.\frac{d}{dt}\right|_{t=0}\varphi_{\alpha}\circ c_{2}\left(t\right)$
for a~chart $\varphi_{\alpha}$ in an atlas $\left(U_{\alpha},\varphi_{\alpha}\right)_{\alpha\in A}$
on $M$. The tangent space at~$x$ then is equal to $C^{\infty}\left(\mathbb{R},M\right)/\sim_{x}$
(compare with \cite[28.12]{KrieglMichor}).

\subsection{Tangent mappings}

Let $f:M\to N$ be a smooth mapping between manifolds. Then $f$ induces
a linear map $\text{d}f\left(x\right):T_{x}M\to T_{f\left(x\right)}N$
for each $x\in M$ by the following formula. If $X=\left[c\right]\in C^{\infty}\left(\mathbb{R},M\right)/\sim_{x}$
then $\text{d}f\left(x\right)X:=\left[f\circ c\right]$. This def\/ines
a f\/iberwise linear map $\text{d}f:TM\to TN$ called the \emph{differential
of} $f$ (compare with \cite[28.15]{KrieglMichor}).

\subsection{Distributions}

A \emph{distribution} on a smooth manifold $M$ is an assignment to
each point $x\in M$ a $c^{\infty}$-closed subspace $\mathcal{D}_{x}$
of $T_{x}M$. If $\mathcal{D}=\left\{ \mathcal{D}_{x}\right\} $ is
a distribution on a manifold $M$ and $i:N\hookrightarrow M$ is the
inclusion map of a connected submanifold $N$ of $M$, then $N$ is
called an \emph{integral manifold} of $\mathcal{D}$ if $\text{d}i\left(T_{x}N\right)=\mathcal{D}_{i\left(x\right)}$
for all $x\in N$. An integral manifold of $\mathcal{D}$ is called
\emph{maximal} if it is not properly contained in any other integral
manifold.

Let $\mathcal{D}$ be a distribution on a manifold $M$. The set of
locally def\/ined vector f\/ields $X$ on $M$ satisfying $X\left(x\right)\in\mathcal{D}_{x}$
will be denoted by $\mathfrak{X}_{\mathcal{D}}\left(M\right)$.

\begin{remark}
In f\/inite dimensions such distributions def\/ined without any assumptions
regarding continuity or smoothness are sometimes called ``generalized
distributions''. If a generalized distribution $\mathcal{D}$ is
``smooth'' in the sense that every $v\in\mathcal{D}_{x}\subset T_{x}M$
can be realized as $X\left(x\right)$ for a locally def\/ined vector
f\/ield $X\in\mathfrak{X}_{\mathcal{D}}\left(M\right)$, then there
exist results on integrability of such distributions (see e.g.~\cite{Sussmann1,
Sussmann2, Stefan1, Stefan2}).
\end{remark}

\subsection{Foliations\label{sub:Foliations}}

Let $M$ be a smooth manifold. A \emph{foliation} of $M$ is a distribution
$\mathcal{F}=\left\{ \mathcal{F}_{x}\right\} $ on $M$, for which
there exists an atlas $\left(U_{\alpha},\varphi_{\alpha}\right)$
of charts $\varphi_{\alpha}:U_{\alpha}\to E_{\alpha}$ on $M$ and
a family of $c^{\infty}$-closed subspaces $\left\{ F_{\alpha}\subset E_{\alpha}\right\} $,
such that the inverse image under $\varphi_{\alpha}$ of translations
of $F_{\alpha}$ are integral manifolds of~$\mathcal{F}$, and such
that if $N\subset U_{\alpha}$ is an integral manifold of~$\mathcal{F}$
then $\varphi_{\alpha}\left(N\right)$ is contained in a~translation
of~$F_{\alpha}$. The charts $\varphi_{\alpha}$ will be called \emph{distinguished
charts}.

Let $\varphi_{\alpha}:U_{\alpha}\to E_{\alpha}$ be a distinguished
chart of a foliation $\mathcal{F}$ of $M$ and $y+F_{\alpha}$ a
translation of $F_{\alpha}\subset E_{\alpha}$. Then $\psi_{\alpha,y}:=\left.\varphi_{\alpha}\right|_{\varphi^{-1}\left(y+F_{\alpha}\right)}-y$
def\/ines a chart into $F_{\alpha}$, and the set of all such charts
$\psi_{\alpha,y}$ def\/ines an alternative smooth structure on the
set $M$ modeled on the spaces $F_{\alpha}$. The set $M$ equipped
with this alternative manifold structure will be denoted by $M^{\mathcal{F}}$.
A \emph{leaf} of the foliation $\mathcal{F}$ is a connected component
of $M^{\mathcal{F}}$. Since $T_{x}M^{\mathcal{F}}=\mathcal{F}_{x}$
for all $x\in M$, every leaf is a maximal integral manifold of $\mathcal{F}$
(compare with \cite[27.16]{KrieglMichor}).

\begin{remark}
Our def\/inition of foliation dif\/fers from the def\/inition in \cite[27.16]{KrieglMichor}
because we wish to describe foliations on manifolds modeled on dif\/ferent
convenient spaces in dif\/ferent charts.
\end{remark}

\subsection{Fiber bundles}

A \emph{fiber bundle $\left(Q,p,M\right)$} consists of manifolds
$Q$ (the \emph{total space}), $M$ (the \emph{base}), and a smooth
mapping $p:Q\to M$ (the \emph{projection}) such that for every $x\in M$
there exists an open neighbourhood $U$ of $x$, a smooth manifold
$S_{U}$, and a dif\/feomorphism $\psi$ such that the following diagram
commutes:\[
\xymatrix{p^{-1}\left(U\right)\ar[r]^{\psi}\ar[dr]^{p} & U\times S_{U}\ar[d]^{\text{pr}_{1}}\\
 & U.}
\]
 Such a pair $\left(U,\psi\right)$ as above is called a \emph{fiber
bundle chart}. A \emph{fiber bundle atlas} $\left(U_{\alpha},\psi_{\alpha}\right)_{\alpha\in A}$
is a~set of f\/iber bundle charts such that $\left\{ U_{\alpha}\right\} _{\alpha\in A}$
is an open cover of $M$. If we f\/ix a f\/iber bundle atlas, then $\psi_{\alpha}\circ\psi_{\beta}^{-1}\left(x,s\right)=\left(x,\psi_{\alpha\beta}\left(x,s\right)\right)$
where $\psi_{\alpha\beta}:\left(U_{\alpha}\cap U_{\beta}\right)\times S_{\beta}\to S_{\alpha}$
is smooth, and where $\psi_{\alpha\beta}\left(x,\cdot\right)$ is
a dif\/feomorphism of $S_{\beta}$ onto $S_{\alpha}$ for each $x\in U_{\alpha\beta}:=U_{\alpha}\cap U_{\beta}$.
The mappings~$\psi_{\alpha\beta}$ are called the \emph{transition
functions} of the bundle. When $S_{U_{\alpha}}=S$ for all charts
$\left(U_{\alpha},\psi_{\alpha}\right)$ for some smooth manifold
$S$, then $S$ is called the \emph{standard fiber} (compare with
\cite[37.1]{KrieglMichor}).

\begin{remark}
Our def\/inition of a f\/iber bundle dif\/fers from the def\/inition in \cite[37.1]{KrieglMichor}
in the sense that it allows for dif\/ferent $S_{U}$ for dif\/ferent neighbourhoods
$U$.
\end{remark}

\subsection{Vector bundles}

Let $\left(Q,p,M\right)$ be a f\/iber bundle. A f\/iber bundle chart
$\left(U,\psi\right)$ is called a \emph{vector bundle chart} if~$S_{U}$
is a convenient vector space. Two vector bundle charts $\left(U_{\alpha},\psi_{\alpha}\right)$
and $\left(U_{\beta},\psi_{\beta}\right)$ are \emph{compatible} if
the transition function $\psi_{\alpha\beta}$ is bounded and linear
in the f\/ibers, i.e.\ $\psi_{\alpha\beta}\left(x,s\right)=\left(x,\phi_{\alpha\beta}\left(x\right)s\right)$
for some mapping $\phi_{\alpha\beta}:U_{\alpha\beta}\to GL\left(S_{\beta},S_{\alpha}\right)\subset L(S_{\beta},S_{\alpha})$.
A \emph{vector bundle atlas} is a f\/iber bundle atlas $\left(U_{\alpha},\psi_{\alpha}\right)_{\alpha\in A}$
consisting of pairwise compatible vector bundle charts. Two vector
bundle atlases are \emph{equivalent} if their union is again a vector
bundle atlas. A \emph{vector bundle $\left(Q,p,M\right)$} is a~f\/iber
bundle together with an equivalence class of vector bundle atlases
(compare with~\mbox{\cite[29.1]{KrieglMichor}}).

\begin{remark}
Here again our def\/inition dif\/fers from that in \cite[29.1]{KrieglMichor}
in that we allow for dif\/ferent~$S_{U}$ in dif\/ferent neighbourhoods
$U$. However, this more general version of vector bundles is subsequently
used implicitly throughout the text (see e.g.~\cite[29.9]{KrieglMichor}
and \cite[29.10]{KrieglMichor} where the tangent bundle $TM$ of
any smooth manifold $M$ is taken to be a vector bundle).
\end{remark}

\subsection{Constructions with vector bundles}

If $Q\to M$ and $R\to M$ are vector bundles then we have vector
bundles $Q^{*}$, $L\left(Q,R\right)$, and $L_{\text{alt}}^{k}\left(Q,R\right)$
whose f\/ibers over $x\in M$ are $\left(Q_{x}\right)^{*}$ (the space
of bounded linear functionals on $Q_{x}$), $L\left(Q_{x},R_{x}\right)$
and $L_{\text{alt}}^{k}\left(Q_{x},R_{x}\right)$ respectively \cite[29.5]{KrieglMichor}.

\begin{remark}
We will always use $E^{*}$ to denote the space of bounded linear
functionals on a locally convex space $E$. In \cite{KrieglMichor}
$E^{*}$ is reserved for the space of continuous (in the locally convex
topology) linear functionals while $E'$ is used to denote the space
of bounded linear functionals.
\end{remark}

\subsection{Cotangent bundles}

Since $TM$ is a vector bundle for any manifold $M$, the bundle $\left(TM\right)^{*}$
with f\/iber over $x\in M$ equal to $\left(T_{x}M\right)^{*}$ is also
a vector bundle. This vector bundle is called the \emph{cotangent
bundle of}~$M$ and will be denoted by~$T^{*}M$ \cite[33.1]{KrieglMichor}.

\subsection{Spaces of sections of vector bundles}\label{sub:Spaces-of-Sections}

A \emph{section} of a vector bundle $p:Q\to M$ is a smooth map $s:M\to Q$
such that $p\circ s=\text{Id}_{M}$. The set of sections of $p:Q\to M$
will be denoted by $\Gamma\left(Q\right)$, and the set of sections
with compact support by $\Gamma_{c}\left(Q\right)$. The space $\Gamma\left(Q\right)$
has a natural convenient structure \cite[30.1]{KrieglMichor}. If
$M$ is f\/inite dimensional and paracompact then $\Gamma_{c}\left(Q\right)$
has a natural convenient structure \cite[30.4]{KrieglMichor}, and
if $W\subset Q$ is an open subset, then $\left\{ s\in\Gamma_{c}\left(Q\right)\mid s\left(M\right)\subset W\right\} $
is $c^{\infty}$-open in $\Gamma_{c}\left(Q\right)$ \cite[30.10]{KrieglMichor}.
If $p:Q\to M$ is a f\/inite dimensional vector bundle over a f\/inite
dimensional paracompact base, then the $c^{\infty}$-topology on $\Gamma_{c}\left(Q\right)$
is induced from the Whitney $C^{\infty}$-topology on $C^{\infty}\left(M,Q\right)$
(which coincides with the Fr\'echet topology; see Section~\ref{sub:cinftyTopology})
\cite[41.13]{KrieglMichor}.

We will be interested primarily in sets that can be locally modeled
on spaces of sections of vector bundles. To understand notions of
smoothness on such sets, it is enough to identify the smooth curves.

\subsection{Curves in spaces of sections}\label{sub:Curves-in-Spaces}

Let $Q\to M$ be a vector bundle. A curve $c:\mathbb{R}\to\Gamma\left(Q\right)$
is smooth if and only if the associated map $c^{\wedge}:\mathbb{R}\times M\to Q$
def\/ined by $c^{\wedge}\left(t,x\right):=c\left(t\right)\left(x\right)$
is smooth \cite[30.8]{KrieglMichor}.

\subsection{Example: manifold of mappings}\label{sub:Example:MMappings}

Let $M$ and $N$ be f\/inite dimensional manifolds. The space $C^{\infty}\left(M,N\right)$
is a smooth manifold modeled on spaces $\Gamma_{c}\left(f^{*}TN\right)$
of compactly supported sections of the pullback bundle along $f\in C^{\infty}\left(M,N\right)$
\cite[42.1]{KrieglMichor}. The charts can be described as follows.
Choose a Riemannian metric on $N$ and let $\exp:TN\supset U\to N$
be the smooth exponential map of this metric. If $\pi_{N}:TN\to N$
denotes the projection of the tangent bundle, then we can assume that
$\left(\pi_{N},\exp\right):U\to N\times N$ is a dif\/feomorphism onto
an open neighbourhood $W$ of the diagonal. For $f,g\in C^{\infty}\left(M,N\right)$,
we will write $f\sim g$ if $f$ and $g$ dif\/fer only on a compact
set in $M$. The charts $\left(V_{f},\psi_{f}\right)$ indexed by
$f\in C^{\infty}\left(M,N\right)$ are given by
\begin{gather*}
V_{f}  :=  \big\{ g\in C^{\infty}\left(M,N\right)\mid g\sim f,\,\left(f\left(x\right),g\left(x\right)\right)\in W\text{ for all }x\in M\big\}, \\
\psi_{f}   : \  V_{f}\to\Gamma_{c}\left(f^{*}TN\right),\\
\psi_{f}\left(g\right)\left(x\right)   :=  \big(x,\exp_{f\left(x\right)}^{-1}\left(g\left(x\right)\right)\big)
=\big(x,\big(\left(\pi_{N},\exp\right)^{-1}\left(f\left(x\right),g\left(x\right)\right)\big)\big).
\end{gather*}

\subsection{Vector f\/ields}

Let $M$ be a smooth manifold. A \emph{vector field} $X$ on $M$
is a smooth section of the tangent bundle $TM$ \cite[32.1]{KrieglMichor}.
The set of all vector f\/ields on $M$ will be denoted by $\mathfrak{X}\left(M\right)$.
Each vector f\/ield $X$ specif\/ies a map\begin{gather*}
C^{\infty}(M,\mathbb{R})   \to   C^{\infty}(M,\mathbb{R}),\\
f   \mapsto   Xf,\\
Xf(x)   :=   \text{d}f(x)X(x).
\end{gather*}

\subsection{The Lie bracket}

Let $X$ and $Y$ be smooth vector f\/ields on a manifold $M$. Each such vector f\/ield is a smooth mapping $M\to TM$ between manifolds, and so it makes sense to compute the dif\/ferentials $\text{d}X$ and $\text{d}Y$. The
\emph{Lie bracket} $\left[X,Y\right]$ of $X$ and $Y$ is the vector
f\/ield on $M$ given by the expression\[
\left[X,Y\right]=\text{d}Y\left(X\right)-\text{d}X\left(Y\right).
\]
 The bracket $\left[\cdot,\cdot\right]:\mathfrak{X}\left(M\right)\times\mathfrak{X}\left(M\right)\to\mathfrak{X}\left(M\right)$
def\/ines a Lie algebra structure on $\mathfrak{X}\left(M\right)$ \cite[32.5, 32.8]{KrieglMichor}.

\subsection{Dif\/ferential forms}

A \emph{differential $k$-form} on a manifold $M$ is a section $\omega\in\Gamma\left(L_{\text{alt}}^{k}\left(TM,M\times\mathbb{R}\right)\right)$.
The set of all dif\/ferential $k$-forms will be denoted by $\Omega^{k}\left(M\right)$
\cite[33.22]{KrieglMichor}.

\subsection{The pullback of a dif\/ferential form}

Let $f:N\to M$ be a smooth mapping and $\omega\in\Omega^{k}\left(M\right)$
be a dif\/ferential $k$-form on $M$. The \emph{pullback} $f^{*}\omega\in\Omega^{k}\left(N\right)$
of $\omega$ is def\/ined by\[
f^{*}\omega_{x}\left(X_{1},\ldots,X_{k}\right):=\omega_{f\left(x\right)}(\text{d}f\left(x\right)X_{1},\ldots,\text{d}f\left(x\right)X_{k}),
\]
see \cite[33.9]{KrieglMichor}.

\subsection{The insertion operator}

For a vector f\/ield $X\in\mathfrak{X}\left(M\right)$ on a manifold
$M$, the \emph{insertion operator} $\imath\left(X\right)$ is def\/ined
by\begin{gather*}
\imath\left(X\right): \ \Gamma\big(L_{\text{alt}}^{k}(TM,M\times\mathbb{R})\big)   \to   \Gamma\big(L_{\text{alt}}^{k-1}(TM,M\times\mathbb{R})\big)\\
\phantom{\imath\left(X\right):} \   \left(\imath\left(X\right)\omega\right)\left(Y_{1},\ldots,Y_{k-1}\right)   :=   \omega(X,Y_{1},\ldots,Y_{k-1}),
\end{gather*}
see \cite[33.10]{KrieglMichor}.

\subsection{The exterior derivative}

Let $U\subset E$ be a $c^{\infty}$-open subset and let $\omega\in C^{\infty}\left(U,L_{\text{alt}}^{k}\left(E,\mathbb{R}\right)\right)$
be a dif\/ferential $k$-form on~$U$. The \emph{exterior derivative}
$d\omega\in C^{\infty}\big(U,L_{\text{alt}}^{k+1}\left(E,\mathbb{R}\right)\big)$
of $\omega$ is the skew symmetrization of the dif\/fe\-ren\-tial $\text{d}\omega$:
\[
(d\omega)(x)(X_{0},\ldots,X_{k})=\sum_{i=0}^{k}(-1)^{i}\text{d}\omega(x)(X_{i})\big(X_{0},\ldots,\hat{X}_{i},\ldots,X_{k}\big).\]
 (Note that the \emph{differential} $\text{d}\omega$ with plain text $\text{d}$ is used to def\/ine the \emph{exterior derivative} $d\omega$ with italicized $d$.) If $\omega$ is a dif\/ferential $k$-form on a manifold $M$, then
this local formula def\/ines a~dif\/ferential $k+1$-form $d\omega$ on
$M$. The above local expression for the exterior derivative induces
the global formula
\begin{gather*}
(d\omega)(x)(X_{0},\ldots,X_{k})  =  \sum_{i=0}^{k}(-1)^{i}X_{i}\big(\omega\circ\big(X_{0},\ldots,\hat{X}_{i},\ldots,X_{k}\big)\big)\\
\phantom{(d\omega)(x)(X_{0},\ldots,X_{k})  =}{}+ \sum_{i<j}(-1)^{i+j}\omega\circ\big([X_{i},X_{j}],X_{0},\ldots,\widehat{X}_{i},\ldots,\widehat{X}_{j},\ldots,X_{k}\big),
\end{gather*}
 where $X_{0},\ldots, X_{k}\in\mathfrak{X}(M)$ \cite[33.12]{KrieglMichor}.

\subsection{Lie groups}

A \emph{Lie group} $G$ is a smooth manifold and a group such that
multiplication $\mu:G\times G\to G$ and inversion $\nu:G\to G$ are
smooth. The \emph{Lie algebra} of a Lie group $G$ is the tangent
space to $G$ at the identity $e$, which inherits a Lie bracket from
the identif\/ication with left invariant vector f\/ields. The Lie algebra
will be denoted either by $\mathfrak{g}$ or $\text{Lie}(G)$
\cite[36.1, 36.3]{KrieglMichor}.

\subsection{Basic dif\/ferential forms}

Let $l:G\times M\to M$ be a smooth action of a Lie group $G$ on
a smooth manifold $M$. Let $l_{g}:M\to M$ denote the left multiplication
mapping $x\mapsto g\cdot x$. For $\xi\in\mathfrak{g}$ the \emph{generating
vector field} $\xi_{M}$ is def\/ined by $\xi_{M}\left(x\right):=\text{d}l_{\left(e,x\right)}\left(\xi,0\right)$.
A dif\/ferential $k$-form $\omega\in\Omega^{k}\left(M\right)$ on $M$
is \emph{$G$-invariant} if $l_{g}^{*}\omega=\omega$ for all $g\in G$
and \emph{horizontal} if $\omega\left(\xi_{M},\cdot\right)=0\in\Omega^{k-1}\left(M\right)$
for all $\xi\in\mathfrak{g}$. A~dif\/ferential $k$-form $\omega\in\Omega^{k}\left(M\right)$
is \emph{basic} if it is both $G$-invariant and horizontal. The set
of all basic $k$-forms on $M$ will be denoted by $\Omega_{\text{hor}}^{k}\left(M\right)^{G}$
\cite[37.23]{KrieglMichor}.

\subsection[Principal $G$ bundles]{Principal $\boldsymbol{G}$ bundles}\label{sub:PrincipalBundles}

Let $G$ be a Lie group. A \emph{principal $G$ bundle} is a f\/iber
bundle $\left(P,p,M,G\right)$ with standard f\/iber $G$ whose transition
functions act on $G$ via left translation: There is a family of smooth
mappings $\left\{ \phi_{\alpha\beta}:U_{\alpha\beta}\to G\right\} $
that satisfy the cocycle condition $\phi_{\alpha\beta}\left(x\right)\phi_{\beta\gamma}\left(x\right)=\phi_{\alpha\gamma}\left(x\right)$
for $x\in U_{\alpha}\cap U_{\beta}\cap U_{\gamma}$, $\phi_{\alpha\alpha}\left(x\right)=e$
(the identity in $G$), and $\psi_{\alpha\beta}\left(x,g\right)=\phi_{\alpha\beta}\left(x\right)\cdot g$
\cite[37.7, 37.8]{KrieglMichor}. The pull back through the projection
$p^{*}:\Omega^{k}\left(M\right)\to\Omega_{\text{hor}}^{k}\left(P\right)^{G}$
is an isomorphism \cite[37.30]{KrieglMichor}.

\subsection{Dif\/feomorphism groups}

The following dif\/feomorphism groups are examples of inf\/inite dimensional
Lie groups:
\begin{itemize}\itemsep=0pt
\item the group $\text{Dif\/f}\left(M\right)$ of dif\/feomorphisms of a f\/inite
dimensional paracompact manifold $M$; the Lie algebra is the space
$\mathfrak{X}_{c}\left(M\right)$ of compactly supported vector f\/ields
on $M$. In fact, $\text{Dif\/f}\left(M\right)$ is open in $C^{\infty}\left(M,M\right)$
so the tangent space at each $f\in\text{Dif\/f}\left(M\right)$ is equal
to $\Gamma\left(f^{*}TM\right)$ \cite[43.1]{KrieglMichor}.
\item the group $\text{Symp}\left(M\right)$ of symplectomorphisms of a
(f\/inite dimensional) symplectic manifold $\left(M,\sigma\right)$;
the Lie algebra is the space $\mathfrak{X}_{c}^{\text{symp}}\left(M\right)$
of compactly supported symplectic vector f\/ields \cite[43.12]{KrieglMichor}.
($\phi\in\text{Dif\/f}\left(M\right)$ is a symplectomorphism if $\phi^{*}\sigma=\sigma$;
$X\in\mathfrak{X}\left(M\right)$ is a~symplectic vector f\/ield if
$\mathcal{L}_{X}\sigma=0$.)
\item the group $\text{Ham}\left(M\right)$ of Hamiltonian symplectomorphisms
of a (f\/inite dimensional) symplectic manifold $\left(M,\sigma\right)$;
the Lie algebra is the space $\mathfrak{X}_{c}^{\text{ham}}\left(M\right)$
of compactly supported Hamiltonian vector f\/ields \cite[43.12, 43.13]{KrieglMichor}.
($X\in\mathfrak{X}\left(M\right)$ is a Hamiltonian vector f\/ield if~$\imath\left(X\right)\sigma$ is exact; $\phi\in\text{Symp}\left(M\right)$
is a Hamiltonian symplectomorphism if it is the time 1 f\/low of a time
dependent Hamiltonian vector f\/ield.)
\end{itemize}

\begin{remark}
In heuristic approaches to inf\/inite dimensional Lie groups, the Lie
algebra to the group $\text{Dif\/f}\left(M\right)$ of dif\/feomorphisms
of a f\/inite dimensional manifold $M$ is often taken to be the space
of smooth vector f\/ields on $M$. In the convenient setup, the Lie
algebra of $\text{Dif\/f}\left(M\right)$ is given by the space $\mathfrak{X}_{c}\left(M\right)$
of compactly supported vector f\/ields on $M$ because of the choice
of charts.
\end{remark}

\subsection{The adjoint representation}

Let $GL\left(E\right)$ denote the set of bounded invertible linear
transformations of $E$. Let $G$ be a Lie group with Lie algebra
$\mathfrak{g}$. Every element $g\in G$ def\/ines an automorphism $\psi_{g}:G\to G$
by conjugation: $\psi_{g}\left(a\right):=gag^{-1}$. The \emph{adjoint
representation of $G$} denoted by $\text{Ad}:G\to GL\left(\mathfrak{g}\right)\subset L\left(\mathfrak{g},\mathfrak{g}\right)$
is given by $\text{Ad}\left(g\right):=\text{d}_{e}\psi_{g}:\mathfrak{g}\to\mathfrak{g}$
for $g\in G$. The \emph{adjoint representation of} $\mathfrak{g}$
denoted by $\text{ad}:\mathfrak{g}\to\mathfrak{gl}\left(\mathfrak{g}\right):=L\left(\mathfrak{g},\mathfrak{g}\right)$
is given by $\text{ad}:=\text{d}_{e}\text{Ad}$ \cite[36.10]{KrieglMichor}.

\subsection{Weak symplectic manifolds}

A 2-form $\sigma\in\Omega^{2}\left(M\right)$ on a manifold $M$ is
called a \emph{weak symplectic structure on} $M$ if it is closed
($d\sigma=0)$ and if its associated vector bundle homomorphism $\sigma^{\flat}:TM\to T^{*}M$
is injective. This last condition is equivalent to weak nondegeneracy:
for every $x\in M$ and $v\in T_{x}M$ there exists a $w\in T_{x}M$
such that $\sigma_{x}\left(v,w\right)\neq0$. If $\sigma^{\flat}:TM\to T^{*}M$
is invertible with a smooth inverse then $\sigma$ is called a \emph{strong
symplectic structure on} $M$ \cite[48.2]{KrieglMichor}. A vector
f\/ield $X\in\mathfrak{X}\left(M\right)$ will be called \emph{Hamiltonian}
if $\imath\left(X\right)\sigma=dH$ for some $H\in C^{\infty}\left(M,\mathbb{R}\right)$,
and the function $H$ will called a~\emph{Hamiltonian} of $X$.

\section{Isodrastic foliations}

In this section we will describe our approach towards describing Lagrangian
submanifolds as Lagrangian embeddings modulo reparametrizations. We
will show that the space of Lagrangian embeddings into a f\/ixed symplectic
manifold $\left(M,\omega\right)$ is a smooth manifold which has a
natural foliation $\mathcal{E}$. Moreover, the space of Lagrangian
embeddings of the form $L_{0}\hookrightarrow M$ is the total space
of a principal $\text{Dif\/f}_{+}\left(L_{0}\right)$ bundle over the
space of Lagrangian submanifolds in $M$. The leaves of $\mathcal{E}$
will turn out to be orbits of the group of Hamiltonian symplectomorphisms
under the natural left composition action. Meanwhile the foliation
$\mathcal{E}$ descends to a foliation $\mathcal{F}$ of the space
of Lagrangian submanifolds in $M$. In all of these constructions,
the key will be to use Weinstein's Lagrangian Neighbourhood Theorem
which says that any symplectic manifold near a Lagrangian $L$ looks
like a neighbourhood of the zero section in the cotangent bundle $T^{*}L$.

Let $\left(M,\omega\right)$ be a f\/inite dimensional symplectic manifold.
Let $L_{0}$ be an oriented, compact manifold of half the dimension
of $M$.

\subsection*{Notation}

By $\text{Lag}\left(L_{0},M\right)$ we will denote the set of Lagrangian
embeddings of $L_{0}$ into $\left(M,\omega\right)$. That is, \[
\text{Lag}\left(L_{0},M\right):=\left\{ i\in\text{Emb}\left(L_{0},M\right)\mid i^{*}\omega=0\right\} .\]
 Let $Z^{k}\left(N\right)$ and $B^{k}\left(N\right)$ denote the
set of closed and exact $k$-forms respectively on a manifold~$N$. That is,
\begin{gather*}
Z^{k}(N)  :=  \text{Ker}(d)\cap\Omega^{k}(N),\qquad
B^{k}(N)  :=  \text{Im}(d)\cap\Omega^{k}(N).
\end{gather*}

We will show that $\text{Lag}\left(L_{0},M\right)$ is a smooth manifold
by def\/ining an atlas of charts using the following Lagrangian neighbourhood
theorem of Weinstein:

\begin{theorem}[see Theorem 6.1 and Corollary 6.2 in~\cite{Weinstein1}]\label{sub:TheoremWLN} Let $L$ be
a Lagrangian submanifold of a symplectic manifold $\left(M,\omega\right)$.
Then there exists an open neighbourhood $U$ of $L$ and a~symplectic
embedding $\psi:U\to T^{*}L$ such that $\left.\psi\right|_{L}=1_{L}$
and $\psi^{*}\omega_{T^{*}L}=\omega$.
\end{theorem}

\begin{proposition}\label{sub:PropLag(L_0,M)}
$\emph{Lag}\left(L_{0},M\right)$ is a smooth manifold modeled
on the space $Z^{1}\left(L_{0}\right)\oplus\mathfrak{X}\left(L_{0}\right)$.
\end{proposition}

\begin{proof}
The idea of the proof is as follows. By Theorem~\ref{sub:TheoremWLN},
Lagrangian submanifolds near a given Lagrangian submanifold can be
identif\/ied with the graphs of closed 1-forms in $T^{*}L_{0}$. It
follows that Lagrangian embeddings near a given one can be identif\/ied
with closed 1-forms viewed as maps from $L_{0}$ to $T^{*}L_{0}$
precomposed with dif\/feomorphisms of $L_{0}$.

Given $i\in\text{Lag}$$\left(L_{0},M\right)$, by Theorem~\ref{sub:TheoremWLN}
the embedding $i$ can be extended on a neighbourhood~$W_{i}$ of
the zero section in $T^{*}L_{0}$ to a~symplectic embedding $\lambda_{i}:W_{i}\to M$.
Let $V_{e}$ be a~chart neighbourhood of the identity map $e\in\text{Dif\/f}\left(L_{0}\right)$
and denote by $\psi_{e}:V_{e}\to\mathfrak{X}\left(L_{0}\right)$ the
corresponding chart as part of an atlas on $\text{Dif\/f}\left(L_{0}\right)$.
Def\/ine
\begin{gather*}
U_{i}  :=  \big\{ l\in\text{Lag}(L_{0},M)\mid l=\lambda_{i}\circ\alpha\circ a,\,\alpha\in Z^{1}(L_{0}),\alpha(L_{0})\subset W_{i},\, a\in V_{e}\big\} ,\\
\varphi_{i}  :  \  U_{i}\to Z^{1}(L_{0})\oplus\mathfrak{X}(L_{0}),\qquad
\varphi_{i}(l)  :=  (\alpha,\psi_{e}(a)).
\end{gather*}
 The space $\mathfrak{X}\left(L_{0}\right)=\Gamma\left(TL_{0}\right)$
is convenient by Section~\ref{sub:Examples-of-Convenient}. The space
$Z^{1}\left(L_{0}\right)$ is a $c^{\infty}$-closed subspace of $\Gamma\left(T^{*}L_{0}\right)$
since it is the kernel of the continuous map $d:\Gamma(T^{*}L_{0})\to\Gamma(\bigwedge^{2}T^{*}L_{0})$,
and therefore it is convenient. The set $\{\alpha\in Z^{1}\left(L_{0}\right)\mid\alpha\left(L_{0}\right)\subset W_{i}\}$
is $c^{\infty}$-open in $Z^{1}\left(L_{0}\right)$ (see Section~\ref{sub:Spaces-of-Sections}).
Thus $\varphi_{i}$ is a bijection of $U_{i}$ onto a $c^{\infty}$-open
subset of $Z^{1}\left(L_{0}\right)\oplus\mathfrak{X}\left(L_{0}\right)$.

The collection $\left(U_{i},\varphi_{i}\right)_{i\in\text{Lag}\left(L_{0},M\right)}$
def\/ines a smooth atlas on $\text{Lag}\left(L_{0},M\right)$, since
the chart chan\-gings $\varphi_{ik}$ are smooth by smoothness of the
exponential map, by smoothness of each symplectic embedding $\lambda_{i}$,
and by Section~\ref{sub:Curves-in-Spaces}.
\end{proof}

To explicitly describe the tangent space to $\text{Lag}\left(L_{0},M\right)$ at
a point $i$, we will make use of the following notation.

\subsection*{Notation}\label{sub:Notalpha_X}
If $S_{0}$ is a manifold (not necessarily of half the dimension of
$M$), then for every $i\in\text{Emb}\left(S_{0},M\right)$ we can
view the tangent bundle $TS_{0}$ as a subbundle of the pullback bundle~$i^{*}TM$. The symplectic form $\omega$ def\/ines a vector bundle
isomorphism $\omega^{\flat}:TM\to T^{*}M$, which induces a vector
space isomorphism $\mu:\Gamma\left(i^{*}TM\right)\to\Gamma\left(i^{*}T^{*}M\right)$.
There is a natural surjection from the pullback bundle~$i^{*}T^{*}M$
onto the cotangent bundle~$T^{*}S_{0}$. This induces a linear map
$\nu:\Gamma\left(i^{*}T^{*}M\right)\to\Gamma\left(T^{*}S_{0}\right)=\Omega^{1}\left(S_{0}\right)$.
For $X\in\Gamma\left(i^{*}TM\right)$, let $\alpha_{X}\in\Omega^{1}\left(S_{0}\right)$
denote the image of $X$ under the composition $\nu\circ\mu$. That
is,
\begin{gather*}
X  \mapsto  \alpha_{X}\in\Omega^{1}(S_{0}),\\
\alpha_{X}(\xi)  :=  \omega_{i(x)}(X(x),\text{d}i(x)\cdot\xi)\qquad\forall\,\xi\in T_{x}S_{0}.\end{gather*}
 Also, set
 \begin{gather*}
\Gamma_{\text{closed}}i^{*}TM)  :=  \big\{ X\in\Gamma(i^{*}TM)\mid\alpha_{X}\in Z^{1}(S_{0})\big\} ,\\
\Gamma_{\text{exact}} (i^{*}TM )   :=   \big\{ X\in\Gamma (i^{*}TM )\mid\alpha_{X}\in B^{1} (S_{0} )\big\} .
\end{gather*}

\begin{remark}
If $Y\in\mathfrak{X}\left(L_{0}\right)$ then $\alpha_{Y}=0$ for
all $i\in\text{Lag}\left(L_{0},M\right)$ since such embeddings are
Lagrangian.
\end{remark}

\begin{proposition}
For each $i\in\emph{Lag}\left(L_{0},M\right)$, the sequence\begin{equation}
0\xrightarrow{}\mathfrak{X}\left(L_{0}\right)\xrightarrow{f_{1}}\Gamma_{\emph{closed}}\left(i^{*}TM\right)
\xrightarrow{f_{2}}Z^{1}\left(L_{0}\right)\xrightarrow{}0,
\label{eq:LagSeq}\end{equation}
 where $f_{1}\left(Y\right)=\text{\rm d}i\left(Y\right)$ and
$f_{2}\left(X\right)=\alpha_{X}$, is a $\emph{Dif\/f}\left(L_{0}\right)$-equivariant
exact sequence.
\end{proposition}

\begin{proof}
If $i\in\text{Lag}\left(L_{0},M\right)$
and $X\in\Gamma\left(i^{*}TM\right)$, then $\alpha_{X}=0$ if and
only if $X$ is tangent to~$i\left(L_{0}\right)$ since $i$ is Lagrangian.
Thus $\text{Ker}\left(f_{2}\right)=\text{Im}\left(f_{1}\right)$.
To check that $f_{2}$ is onto, let $\alpha\in Z^{1}\left(L_{0}\right)$.
By Theorem~\ref{sub:TheoremWLN}, it is enough to prove the assertion
when $M=T^{*}L_{0}$ and $i$ is the zero section inclusion. Let $\pi:T^{*}L_{0}\to L_{0}$
denote the canonical projection. Def\/ine $Z_{\alpha}\in\mathfrak{X}\left(T^{*}L_{0}\right)$
by\[
\imath\left(Z_{\alpha}\right)\omega:=\pi^{*}\alpha.\]
 Then $\alpha_{Z_{\alpha}\circ i}=\alpha$, which means $f_{2}$ is
surjective and so the sequence is exact.
\end{proof}

\begin{remark}\label{sub:SplitMap}
Each symplectic embedding $\lambda_{i}:T^{*}L_{0}\supset W_{i}\to M$
def\/ines a splitting map $s_{i}:Z^{1}\left(L_{0}\right)\to\Gamma_{\text{closed}}\left(i^{*}TM\right)$
of the exact sequence~\eqref{eq:LagSeq} given by
\[
s_{i}\left(\alpha\right):=\text{d}\lambda_{i}\circ\left.Z_{\alpha}\right|_{L_{0}},
\]
 where $\imath\left(Z_{\alpha}\right)\omega_{T^{*}L_{0}}=\pi^{*}\alpha$.
Under this splitting $\Gamma_{\text{closed}}(i^{*}TM)=Z^{1}(L_{0})\oplus\mathfrak{X}(L_{0})$.
\end{remark}

\begin{proposition}
For each $i\in\emph{Lag}\left(L_{0},M\right)$,\[
T_{i}\emph{Lag}\left(L_{0},M\right)=\Gamma_{\emph{closed}}\left(i^{*}TM\right).\]
\end{proposition}

\begin{proof}
Let $j_{t}$ be a smooth curve in $\text{Lag}\left(L_{0},M\right)$
such that $j_{0}=i$. For each $x\in L_{0}$ we have that $\left.\frac{d}{dt}\right|_{t=0}j_{t}\left(x\right)\in T_{i\left(x\right)}M$,
which means $\left.\frac{d}{dt}\right|_{t=0}j_{t}\in\Gamma\left(i^{*}TM\right)$.
The fact that $j_{t}^{*}\omega=0$ implies that $\left.\frac{d}{dt}\right|_{t=0}j_{t}\in\Gamma_{\text{closed}}\left(i^{*}TM\right)$.
If $\left(U_{l},\varphi_{l}\right)$ is a chart corresponding to a
symplectic embedding $\lambda_{l}:T^{*}L_{0}\supset W_{l}\to M$ with
$i\in U_{l}$, then the derivative $\left.\frac{d}{dt}\right|_{t=0}\varphi_{l}\left(j_{t}\right)$
in $Z^{1}\left(L_{0}\right)\oplus\mathfrak{X}\left(L_{0}\right)\subset\left.T\left(T^{*}L_{0}\right)\right|_{L_{0}}$
is canonically identifed with $\left.\frac{d}{dt}\right|_{t=0}j_{t}$
via \[
\left.\frac{d}{dt}\right|_{t=0}\varphi_{l}\left(j_{t}\right)=\alpha_{\left.\frac{d}{dt}\right|_{t=0}j_{t}}+\text{d}\pi\circ\text{d}\lambda_{l}^{-1}\circ\left.\frac{d}{dt}\right|_{t=0}j_{t}.\]
 So, $T_{i}\text{Lag}\left(L_{0},M\right)\subset\Gamma_{\text{closed}}\left(i^{*}TM\right)$.

Conversely, suppose that $X\in\Gamma_{\text{closed}}\left(i^{*}TM\right)$
and denote by $\lambda_{i}:T^{*}L_{0}\supset W_{i}\to M$ the symplectic
embedding associated to the chart $\left(U_{i},\varphi_{i}\right)$.
Def\/ine a smooth curve in $\text{Symp}\left(T^{*}L_{0}\right)$ by
$\psi_{t}\left(x,p\right):=\left(x,p+t\alpha_{X}\left(x\right)\right)$
and a vector f\/ield on $L_{0}$ by $Y:=\text{d}\pi\circ\text{d}\lambda_{i}^{-1}\circ X$.
If $a_{t}$ denotes the f\/low generated by $Y$ on $L_{0}$ and $\text{\textbf{O}}$
denotes the zero section, then $\text{d}\lambda_{i}\circ\left.\frac{d}{dt}\right|_{t=0}\left(\psi_{t}\circ\text{\textbf{O}}\circ a_{t}\right)=X$.
Thus $X\in\Gamma_{\text{closed}}\left(i^{*}TM\right)$ corresponds
to the class $\left[\lambda_{i}\circ\psi_{t}\circ\text{\textbf{O}}\circ a_{t}\right]$
in $T_{i}\text{Lag}\left(L_{0},M\right)$.
\end{proof}

The set $\text{Ham}\left(M\right)$
of Hamiltonian symplectomorphisms is a subgroup of $\text{Symp}\left(M\right)$
(see e.g.~\cite{Polterovich}). So left composition def\/ines an action
of $\text{Ham}\left(M\right)$ on $\text{Lag}\left(L_{0},M\right)$
via\[
\phi\cdot i:=\phi\circ i.
\]

\begin{proposition}\label{sub:PropLag(L_0,M)Foliation}
The spaces $\Gamma_{\emph{exact}}(i^{*}TM){\subset}\Gamma_{\emph{closed}}(i^{*}TM)$
and charts $(U_{i},\varphi_{i})$ for $i\in\emph{Lag}(L_{0},M)$
define a foliation $\mathcal{E}$ of $\emph{Lag}(L_{0},M)$,
whose leaves consist of $\emph{Ham}(M)$ orbits.
\end{proposition}

\begin{proof}
Set $\mathcal{E}_{i}:=\Gamma_{\text{exact}}\left(i^{*}TM\right)$
for $i\in\text{Lag}\left(L_{0},M\right)$. We will f\/irst show that
$\Gamma_{\text{exact}}\left(i^{*}TM\right)$ is a~$c^{\infty}$-closed
(i.e.~convenient) subspace of $\Gamma_{\text{closed}}\left(i^{*}TM\right)$
for all $i\in\text{Lag}\left(L_{0},M\right)$. If $\Gamma_{\text{closed}}\left(i^{*}TM\right)$
is identif\/ied with $Z^{1}\left(L_{0}\right)\oplus\mathfrak{X}\left(L_{0}\right)$
via the splitting map $s_{i}:Z^{1}\left(L_{0}\right)\to\Gamma_{\text{closed}}\left(i^{*}TM\right)$,
then to show that $\Gamma_{\text{exact}}\left(i^{*}TM\right)$ is
$c^{\infty}$-closed in $\Gamma_{\text{closed}}\left(i^{*}TM\right)$
it is enough to show that $F:=B^{1}\left(L_{0}\right)\oplus\mathfrak{X}\left(L_{0}\right)$
is $c^{\infty}$-closed in $Z^{1}\left(L_{0}\right)\oplus\mathfrak{X}\left(L_{0}\right)$.
Let $c_{1}:\mathbb{R}\to F$ be a smooth curve in $F$. If $c_{1}\left(t\right)=\left(df_{t},Y_{t}\right)$
then $p_{1}:t\mapsto f_{t}$ can be chosen to be a smooth curve in
$C^{\infty}\left(L_{0},\mathbb{R}\right)$. Since $C^{\infty}\left(L_{0},\mathbb{R}\right)$
is convenient, there exists a smooth curve $p_{2}:t\mapsto g_{t}$
in $C^{\infty}\left(L_{0},\mathbb{R}\right)$ such that $p_{2}'=p_{1}$.
Similarly for the curve $q_{1}:t\mapsto Y_{t}$ there exists a smooth
curve $q_{2}:t\mapsto Z_{t}$ in $\mathfrak{X}\left(L_{0}\right)$
such that $q_{2}'=q_{1}$. Then $c_{2}\left(t\right):=\left(dg_{t},Z_{t}\right)$
is an antiderivative of $c_{1}$, i.e.\ $c_{2}'=c_{1}$. This means
$F$ is a convenient subspace.

We will next show that the $\text{Ham}\left(M\right)$ orbits in $\text{Lag}\left(L_{0},M\right)$
are maximal integral manifolds of $\mathcal{E}=\left\{ \mathcal{E}_{i}\right\} $.
The tangent vectors to a $\text{Ham}\left(M\right)$ orbit at a point
$i\in\text{Lag}\left(L_{0},M\right)$ are of the form $X_{H}\circ i$
where $X_{H}$ is a Hamiltonian vector f\/ield on $M$. Since $\alpha_{X_{H}\circ i}=i^{*}dH$
it follows that $T_{i}(\text{Ham}(M)\cdot i)\subset\mathcal{E}_{i}$.
Conversely, if $X\in\mathcal{E}_{i}$ then there exists a Hamiltonian
vector f\/ield~$X_{H}$ def\/ined on a neighbourhood of $i(L_{0})$ satisfying
$X=X_{H}\circ i$. By multiplying $H$ by a cutof\/f function which
is equal to 1 near $i(L_{0})$ we may assume that $X_{H}$ is def\/ined
on all of $M$. It follows that $\mathcal{E}_{i}\subset T_{i}(\text{Ham}(M)\cdot i)$.
So $\text{Ham}\left(M\right)$ orbits are integral manifolds. To show
they are maximal, we f\/irst consider the case when $M=T^{*}L_{0}$.
Let $i:L_{0}\hookrightarrow T^{*}L_{0}$ denote the zero section inclusion
and $\left(U_{i},\varphi_{i}\right)$ the corresponding chart on $\text{Lag}\left(L_{0},T^{*}L_{0}\right)$.
Let $j_{t}$ be a smooth curve in an integral manifold $N$ contained
in $U_{i}$. For every $t$, $\frac{d}{dt}j_{t}\in\Gamma_{\text{exact}}\left(j_{t}^{*}T\left(T^{*}L_{0}\right)\right)$
which means $\alpha_{\frac{d}{dt}j_{t}}=dh_{t}$ for a family of functions
$h_{t}\in C^{\infty}\left(L_{0},\mathbb{R}\right)$. This family $h_{t}$
can be chosen to be a~smooth curve in $C^{\infty}\left(L_{0},\mathbb{R}\right)$.
Since $j_{t}$ is contained in $U_{i}$, there exist smooth curves
$\beta_{t}$ in $Z^{1}\left(L_{0}\right)$ and $a_{t}$ in $\text{Dif\/f}\left(L_{0}\right)$
such that $j_{t}=\beta_{t}\circ a_{t}$. Since $\frac{d}{dt}(\beta_{t}\circ a_{t})=\left(\frac{d}{dt}\beta_{t}\right)\circ a_{t}+\text{d}\beta_{t}\left(\frac{d}{dt}a_{t}\right)$,
it follows that $-a_{t}^{*}\frac{d}{dt}\beta_{t}=dh_{t}$ for all
$t$. Thus, $\beta_{t}=\psi_{t}\circ\beta_{0}$ where $\psi_{t}\left(x,p\right)=\left(x,p-\left(d\int_{0}^{t}h_{s}\circ a_{s}^{-1}ds\right)_{x}\right)$,
which means $\beta_{t}$ is contained in the $\text{Ham}\left(T^{*}L_{0}\right)$
orbit through $\beta_{0}$. For each $t$ we can write $\beta_{0}\circ a_{t}=b_{t}\circ\beta_{0}$
where $b_{t}\in\text{Dif\/f}\left(\text{Graph}\left(\beta_{0}\right)\right)$
in such a way that $b_{t}$ is a smooth curve in $\text{Dif\/f}\left(\text{Graph}\left(\beta_{0}\right)\right)$.
Def\/ine $c_{t}:=b_{t}\circ b_{0}^{-1}$. By means of Theorem~\ref{sub:TheoremWLN}
we can identify an open neighbourhood of~$\text{Graph}\left(\beta_{0}\right)$
with an open neighbourhood of the zero section in $T^{*}\text{Graph}\left(\beta_{0}\right)$.
Under this identif\/ication, if~$C_{t}$ denotes the cotangent lift
of $c_{t}$, i.e. $C_{t}\left(y,z\right)=\left(c_{t}\left(y\right),\left(c_{t}^{-1}\right)^{*}z\right)$,
then $\beta_{0}\circ a_{t}=C_{t}\circ b_{0}\circ\beta_{0}$. Since~$c_{t}$ is a smooth curve in $\text{Dif\/f}\left(\text{Graph}\left(\beta_{0}\right)\right)$
passing through the identity map, the cotangent lift $C_{t}$ is a
smooth curve in $\text{Ham}\left(T^{*}\text{Graph}\left(\beta_{0}\right)\right)$.
Thus $j_{t}=\beta_{t}\circ a_{t}=\psi_{t}\circ C_{t}\circ b_{0}\circ\beta_{0}$
lies in the $\text{Ham}\left(T^{*}L_{0}\right)$ orbit through $b_{0}\circ\beta_{0}$,
which means the integral manifold $N$ is contained in a $\text{Ham}\left(T^{*}L_{0}\right)$
orbit. For the general case when $M$ is any symplectic manifold,
the previous discussion implies that the intersection of any integral
manifold with a chart neighbourhood~$U_{i}$ on $\text{Lag}\left(L_{0},M\right)$
lies in a $\text{Ham}\left(M\right)$ orbit. Thus any integral manifold
containing a point $i\in\text{Lag}\left(L_{0},M\right)$ is contained
in $\text{Ham}\left(M\right)\cdot i$, which means that such orbits
are maximal integral manifolds.

Finally, we will show that the atlas $\left(U_{i},\varphi_{i}\right)_{i\in\text{Lag}\left(L_{0},M\right)}$
consists of distinguished charts. The af\/f\/ine translations of $F$
in $Z^{1}\left(L_{0}\right)\oplus\mathfrak{X}\left(L_{0}\right)$
consist of elements $\left(\alpha,Y\right)$ which are pairwise cohomologous
in the f\/irst factor, i.e $\left(\alpha,Y\right)\in\left(\alpha_{0},Y_{0}\right)+F$
if and only if $\alpha-\alpha_{0}\in B^{1}\left(L_{0}\right)$. Let
$\left(U_{i},\varphi_{i}\right)$ be a chart on $\text{Lag}\left(L_{0},M\right)$
with corresponding symplectic embedding $\lambda_{i}:T^{*}L_{0}\supset W_{i}\to M$.
The zero section in $T^{*}L_{0}$ can be deformed to the graph of
any 1-form $\alpha\in\Omega^{1}\left(L_{0}\right)$ on $L_{0}$ by
taking the time 1 f\/low of the transformation $\left(x,p\right)\mapsto\left(x,p+t\alpha_{x}\right)$
of the cotangent bundle. When $\alpha$ is closed this transformation
is symplectic; when $\alpha$ is exact it is a Hamiltonian symplectomorphism.
So the graph of any exact form can be obtained by deforming the zero
section in $T^{*}L_{0}$ along a Hamiltonian vector f\/ield. Conversely,
suppose that $\phi\in\text{Ham}\left(T^{*}L_{0}\right)$ is a Hamiltonian
symplectomorphism and $\left\{ \psi_{t}\right\} $ is a collection
of symplectomorphisms satisfying $\psi_{0}=\text{Id}$, $\psi_{1}=\phi$,
and $\dot{\psi}_{t}=X_{H_{t}}\circ\psi_{t}$ for some family of Hamiltonian
vector f\/ields $X_{H_{t}}$ on $T^{*}L_{0}$. If $\text{\textbf{O}}$
denotes the zero section, then $j_{t}:=\varphi_{i}\left(\lambda_{i}\circ\psi_{t}\circ\text{\textbf{O}}\right)$
is a smooth curve in $Z^{1}\left(L_{0}\right)\oplus\mathfrak{X}\left(L_{0}\right)$
with time derivative equal to $\left(\alpha_{\text{d}\lambda_{i}\circ\dot{\psi}_{t}\circ\text{\textbf{O}}},\text{d}\pi\circ\dot{\psi}_{t}\circ\text{\textbf{O}}\right)$
for all $t$. Since $\alpha_{\text{d}\lambda_{i}\circ\dot{\psi}_{t}\circ\text{\textbf{O}}}=\left(\psi_{t}\circ\text{\textbf{O}}\right)^{*}dH_{t}$,
the curve $j_{t}$ must be contained in $B^{1}\left(L_{0}\right)\oplus\mathfrak{X}\left(L_{0}\right)$,
which means $\phi\circ\text{\textbf{O}}=\psi_{1}\circ\text{\textbf{O}}$
is the graph of an exact 1-form precomposed with a dif\/feomorphism
of $L_{0}$. It follows that two 1-forms are cohomologous if and only
if their graphs in $T^{*}L_{0}$ can be joined by f\/lowing along a
Hamiltonian vector f\/ield. So a curve in $U_{i}$ lies in a $\text{Ham}\left(M\right)$
orbit if and only if it is mapped into a translation of $F$ under
$\varphi_{i}$. Thus inverse images of translations of $F$ are integral
manifolds of $\mathcal{E}$ and intersections of integral manifolds
with each chart neighbourhood $U_{i}$ get mapped into translations
of $F$ under $\varphi_{i}$.
\end{proof}

\begin{definition}
In the spirit of Weinstein's terminology in \cite{Weinstein3}, we
will call the foliation $\mathcal{E}$ the \emph{isodrastic foliation}
of $\text{Lag}\left(L_{0},M\right)$. An individual leaf of $\mathcal{E}$
will be called an \emph{isodrast} in $\text{Lag}\left(L_{0},M\right)$.
\end{definition}

The group of orientation preserving dif\/feomorphisms $\text{Dif\/f}_{+}\left(L_{0}\right)$
acts freely on $\text{Lag}\left(L_{0},M\right)$ via
\[
a\cdot i:=i\circ a^{-1}.
\]
 The quotient $\text{Lag}\left(L_{0},M\right)/\text{Dif\/f}_{+}\left(L_{0}\right)$
is naturally identif\/ied with the set of oriented, compact Lagrangian
submanifolds in $M$ dif\/feomorphic to $L_{0}$.

\subsection*{Notation}

Set\[
\text{Lag}\left(M\right):=\text{Lag}\left(L_{0},M\right)/\text{Dif\/f}_{+}\left(L_{0}\right).\]

\begin{proposition}\label{sub:PropTLag(M)}
The set $\emph{Lag}\left(M\right)$ of oriented
Lagrangian submanifolds in~$M$ diffeomorphic to~$L_{0}$ is
a smooth manifold modeled on spaces $Z^{1}\left(L\right)$ for $L\in\emph{Lag}\left(M\right)$.
The tangent spaces to $\emph{Lag}\left(M\right)$ are
given by\[
T_{L}\emph{Lag}\left(M\right)=Z^{1}\left(L\right)\]
 and for each representative $i\in\emph{Lag}\left(L_{0},M\right)$
in the class $L\in\emph{Lag}\left(M\right)$,\[
T_{\left[i\right]}\emph{Lag}\left(M\right)=\Gamma_{\emph{closed}}\left(i^{*}TM\right)/\mathfrak{X}\left(L_{0}\right).\]
\end{proposition}

\begin{proof} We will f\/irst describe the manifold structure on $\text{Lag}\left(M\right)$.
For each Lagrangian $L\in\text{Lag}\left(M\right)$, by Theorem~\eqref{sub:TheoremWLN}
there exists a symplectic embedding $\lambda_{L}:W_{L}\to T^{*}L$
of an open neighbourhood of $L$ onto an open neighbourhood of the
zero section in the cotangent bundle. Def\/ine
\begin{gather*}
U_{L}:   =   \big\{ N\in\text{Lag}(M )\mid N\subset W_{L},\,\lambda_{L} (N )=\text{Graph} (\alpha ),\,\alpha\in Z^{1}(L)\big\} ,\\
\varphi_{L}   :  \  U_{L}\to Z^{1} (L ),\qquad
\varphi_{L} (N )   :=   \alpha.
\end{gather*}
 The mapping $\varphi_{L}$ maps $U_{L}$ onto the set $\left\{ \alpha\in Z^{1}\left(L\right)\mid\alpha\left(L\right)\subset\lambda_{L}\left(W_{L}\right)\right\} $
which is $c^{\infty}$-open \linebreak in~$Z^{1}\left(L\right)$. Thus the collection
$\left(U_{L},\varphi_{L}\right)_{L\in\text{Lag}\left(M\right)}$ def\/ines
a smooth atlas on $\text{Lag}\left(M\right)$ as the transition functions
$\varphi_{LN}$ are smooth by smoothness of the symplectic embeddings
$\lambda_{L}$.

As for the tangent spaces, suppose that $K_{t}$ is a smooth curve
in $\text{Lag}\left(M\right)$ such that $K_{0}=L$. If $\left(U_{L'},\varphi_{L'}\right)$
is a chart with $L\in U_{L'}$, and $\varphi_{L'}\left(K_{t}\right)=\text{Graph}\left(\alpha_{t}\right)$
for a smooth curve $\alpha_{t}$ in $Z^{1}\left(N\right)$, then the
derivative $\left.\frac{d}{dt}\right|_{t=0}\varphi_{N}\left(K_{t}\right)\in Z^{1}\left(L'\right)$
canonically def\/ines an element $\beta\in Z^{1}\left(L\right)$ via\[
\alpha_{0}^{*}\beta:=\left.\frac{d}{dt}\right|_{t=0}\varphi_{L'}\left(K_{t}\right).\]
 Thus, $T_{L}\text{Lag}\left(M\right)\subset Z^{1}\left(L\right)$.
Conversely, each $\beta\in Z^{1}\left(L\right)$ def\/ines a~smooth
curve in $\text{Symp}\left(T^{*}L\right)$ via $\psi_{t}\left(x,p\right):=\left(x,p+t\beta\left(x\right)\right)$.
If $\text{\textbf{O}}$ denotes the zero section of $T^{*}L$, then
$\text{Graph}\left(\psi_{t}\circ\text{\textbf{O}}\right)$ is a~smooth
curve in $\text{Lag}\left(T^{*}L\right)$. If $\lambda_{L}:M\supset W_{L}\to T^{*}L$
is the symplectic embedding associated to the chart $\left(U_{L},\varphi_{L}\right)$,
then $c\left(t\right):=\lambda_{L}^{-1}\left(\text{Graph}\left(\psi_{t}\circ\text{\textbf{O}}\right)\right)$
is a smooth curve in $\text{Lag}\left(M\right)$ such that $\left.\frac{d}{dt}\right|_{t=0}\varphi_{L}\circ c\left(t\right)=\beta$.
So $Z^{1}\left(L\right)\subset T_{L}\text{Lag}\left(M\right)$.

We will now describe the identif\/ication of tangent spaces of $\text{Lag}\left(M\right)$
with vector spaces $\Gamma_{\text{closed}}\left(i^{*}TM\right)/\mathfrak{X}\left(L_{0}\right)$.
Let $i\in\text{Lag}\left(L_{0},M\right)$ be a representative in the
class $L\in\text{Lag}\left(M\right)$. Let $\lambda_{i}:T^{*}L_{0}\supset W_{i}\to M$
be the symplectic embedding chosen in the def\/inition of the chart
$\left(U_{i},\varphi_{i}\right)$ on $\text{Lag}\left(L_{0},M\right)$,
and $s_{i}:Z^{1}\left(L_{0}\right)\to\Gamma_{\text{closed}}\left(i^{*}TM\right)$
the corresponding splitting map (see Remark~\ref{sub:SplitMap}).
Then the linear map
\begin{gather*}
Z^{1}\left(L\right)   \to   \Gamma_{\text{closed}} (i^{*}TM )/\mathfrak{X} (L_{0} ),\qquad
\alpha   \mapsto  [s_{i}(i^{*}\alpha)]
\end{gather*}
 is a vector space isomorphism.
\end{proof}

\begin{proposition}\label{sub:PropLag(L_0,M)Lag(M)}
The manifold $\emph{Lag}\left(L_{0},M\right)$ is the
total space of a principal $\emph{Dif\/f}_{+}\left(L_{0}\right)$ bundle
over $\emph{Lag}\left(M\right)$.
\end{proposition}

\begin{proof} We begin by describing
a f\/iber bundle atlas. Let $p:\text{Lag}\left(L_{0},M\right)\to\text{Lag}\left(M\right)$
denote the projection to the quotient. For $i\in\text{Lag}\left(L_{0},M\right)$
let $\lambda_{i}:T^{*}L_{0}\supset W_{i}\to M$ be the symplectic
embedding chosen in def\/ining the chart $\left(U_{i},\varphi_{i}\right)$
on $\text{Lag}\left(L_{0},M\right)$. Def\/ine
\[
U_{[i]}:=\big\{ N\in\text{Lag} (M )\mid N=\lambda_{i} (\text{Graph} (\alpha ) ),\,\alpha\in Z^{1} (L_{0} ),\,\alpha (L_{0} )\subset W_{i}\big\} .
\]
 Then $p^{-1}\left(U_{\left[i\right]}\right)$ consists of all $l\in\text{Lag}\left(L_{0},M\right)$
such that $l=\lambda_{i}\circ\alpha\circ a$ where $\alpha\in Z^{1}\left(L_{0}\right)$,
$\alpha\left(L_{0}\right)\subset W_{i}$ , and $a\in\text{Dif\/f}_{+}\left(L_{0}\right)$.
Def\/ine $\psi_{\left[i\right]}:p^{-1}\left(U_{\left[i\right]}\right)\to U_{\left[i\right]}\times\text{Dif\/f}_{+}\left(L_{0}\right)$
by{\samepage
\[
\psi_{[i]}(\lambda_{i}\circ\alpha\circ a):=(\lambda_{i}(\text{Graph}(\alpha)),a)\]
 so that the collection $\left(U_{\left[i\right]},\psi_{\left[i\right]}\right)_{i\in\text{Lag}\left(L_{0},M\right)}$
def\/ines a f\/iber bundle atlas.}

If $N\in U_{\left[i\right]\left[j\right]}:=U_{\left[i\right]}\cap U_{\left[j\right]}$,
and $\psi_{\left[j\right]}^{-1}\left(N,a\right)=\lambda_{j}\circ\alpha\circ a$
then\[
\psi_{\left[i\right]\left[j\right]}\left(N,a\right)=\left(N,\pi\circ\lambda_{i}^{-1}\circ\lambda_{j}\circ\alpha\circ a\right).\]
 So if $N\in U_{\left[i\right]\left[j\right]}$ and $N=\lambda_{j}\left(\text{Graph}\left(\alpha\right)\right)$
then def\/ine $u_{\left[i\right]\left[j\right]}:U_{\left[i\right]\left[j\right]}\to\text{Dif\/f}_{+}\left(L_{0}\right)$
by\[
u_{\left[i\right]\left[j\right]}\left(N\right):=\pi\circ\lambda_{i}^{-1}\circ\lambda_{j}\circ\alpha.\]
 It follows that if $N\in U_{\left[i\right]}\cap U_{\left[j\right]}\cap U_{\left[k\right]}$
and $N=\lambda_{j}\left(\text{Graph}\left(\alpha\right)\right)=\lambda_{k}\left(\text{Graph}\left(\beta\right)\right)$
then\begin{gather*}
u_{\left[i\right]\left[j\right]}\left(N\right)u_{\left[j\right]\left[k\right]}\left(N\right)   =  \pi\circ\lambda_{i}^{-1}\circ\lambda_{j}\circ\alpha\circ\pi\circ\lambda_{j}^{-1}\circ\lambda_{k}\circ\beta
   =   \pi\circ\lambda_{i}^{-1}\circ\lambda_{k}\circ\beta  =  u_{[i][k]}(N),\\
u_{\left[i\right]\left[i\right]}\left(N\right)=\text{Id}_{L_{0}},\\
\psi_{\left[i\right]\left[j\right]}\left(N,a\right)=u_{\left[i\right]\left[j\right]}\circ a.
\end{gather*}
 So the f\/iber bundle atlas $\left(U_{\left[i\right]},\psi_{\left[i\right]}\right)_{i\in\text{Lag}\left(L_{0},M\right)}$
and the collection of maps $u_{\left[i\right]\left[j\right]}:U_{\left[i\right]\left[j\right]}\to\text{Dif\/f}_{+}\left(L_{0}\right)$
def\/ine a principal $\text{Dif\/f}_{+}\left(L_{0}\right)$ structure.
\end{proof}

The group of Hamiltonian symplectomorphisms $\text{Ham}\left(M\right)$
acts on $\text{Lag}\left(M\right)$ via\[
\text{Ham}\left(M\right)\circlearrowright\text{Lag}\left(M\right):\,\phi\cdot L:=\phi\left(L\right).\]
 As in the proof of Proposition~\ref{sub:PropLag(L_0,M)Foliation}, the
family of subspaces $B^{1}\left(L\right)\subset T_{L}\text{Lag}\left(M\right)$
and charts $\left(U_{L},\varphi_{L}\right)$ for $L\in\text{Lag}\left(M\right)$
def\/ine a foliation \emph{$\mathcal{F}$} on $\text{Lag}\left(M\right)$.
This foliation has f\/inite codimension since the transverse space at
each Lagrangian submanifold $L$ is modeled on $H^{1}\left(L\right)$.

\begin{definition}
The foliation $\mathcal{F}$ will be called the \emph{isodrastic foliation}
of $\text{Lag}\left(M\right)$, and each leaf of $\mathcal{F}$ will
be called an \emph{isodrast} in $\text{Lag}\left(M\right)$.
\end{definition}

The foliation $\mathcal{E}$ gives the set $\text{Lag}\left(L_{0},M\right)$
the structure of a smooth manifold $\text{Lag}\left(L_{0},M\right)^{\mathcal{E}}$
modeled on spaces $\Gamma_{\text{exact}}\left(i^{*}TM\right)$ for
$i\in\text{Lag}\left(L_{0},M\right)$. Similarly, the foliation $\mathcal{F}$
gives $\text{Lag}\left(M\right)$ the structure of a smooth manifold
$\text{Lag}\left(M\right)^{\mathcal{F}}$ modeled on spaces $B^{1}\left(L\right)$
for $L\in\text{Lag}\left(M\right)$. As in Proposition~\ref{sub:PropLag(L_0,M)Lag(M)},
f\/iber bundle charts can be chosen to def\/ine a principal $\text{Dif\/f}_{+}\left(L_{0}\right)$
bundle $p:\text{Lag}\left(L_{0},M\right)^{\mathcal{E}}\to\text{Lag}\left(M\right)^{\mathcal{F}}$.
This bundle restricts to a principal $\text{Dif\/f}_{+}\left(L_{0}\right)$
bundle over each connected component of $\text{Lag}\left(M\right)^{\mathcal{F}}$,
i.e.\ over each isodrast in $\text{Lag}\left(M\right)$.

\section{Weighted Lagrangian submanifolds}

In this section we introduce the notion of weightings and weighted
submanifolds. The set $\text{Lag}_{w}\left(L_{0},M\right)$ of pairs
$\left(i,\eta\right)$ consisting of Lagrangian embeddings $i:L_{0}\hookrightarrow M$
and top degree forms $\eta$ that satisfy $\int_{L_{0}}\eta=1$ has
the smooth structure of the Cartesian product $\text{Lag}\left(L_{0},M\right)\times\big\{ \eta\in\Omega^{n}\left(L_{0}\right)\mid\int_{L_{0}}\eta=1\big\} $.
The foliation $\mathcal{E}$ of $\text{Lag}\left(L_{0},M\right)$
canonically induces a foliation $\mathcal{E}_{w}$ of $\text{Lag}_{w}\left(L_{0},M\right)$.
The space $\text{Lag}_{w}\left(L_{0},M\right)$ is the total space
of a principal $\text{Dif\/f}_{+}\left(L_{0}\right)$ bundle, whose
base can be identif\/ied with the set $\text{Lag}_{w}\left(M\right)$
of Lagrangian submanifolds in $M$ equipped with a~top degree form
of total measure 1. The foliation $\mathcal{E}_{w}$ descends to a
foliation $\mathcal{F}_{w}$ of the base, so that $\text{Lag}_{w}\left(L_{0},M\right)^{\mathcal{E}_{w}}$
(cf.\ Section~\ref{sub:Foliations}) is the total space of a principal
$\text{Dif\/f}_{+}\left(L_{0}\right)$ bundle over $\text{Lag}_{w}\left(M\right)^{\mathcal{F}_{w}}$.
On each leaf of $\mathcal{E}_{w}$ we def\/ine a 2-form $\Omega$, basic
with respect to this principal group action, which descends to a weak
symplectic structure on $\text{Lag}_{w}\left(M\right)^{\mathcal{F}_{w}}$.
Finally, we show that the tangent spaces of $\text{Lag}_{w}\left(M\right)$
and of leaves of $\mathcal{F}_{w}$ can be identif\/ied with the tangent
space descriptions in Weinstein's original construction, and that
Weinstein's symplectic structure $\Omega^{W}$ corresponds to our
symplectic structure $\Omega$.

\begin{definition}
A \emph{weighting} of a compact oriented manifold $L$ is a top degree
form $\rho$ on $L$ satisfying $\int_{L}\rho=1$. A pair $\left(L,\rho\right)$
will be called a \emph{weighted manifold}.
\end{definition}

\subsection*{Notation}

Let $\Omega_{1}^{n}\left(S_{0}\right)$ denote the set of $n$-forms
on a manifold $S_{0}$ that integrate to $1$ (where $n=\dim S_{0}$),
$\Omega_{0}^{n}\left(S_{0}\right)$ the set of $n$-forms on $S_{0}$
that integrate to 0, and $\text{Lag}_{w}\left(L_{0},M\right)$ the
product $\text{Lag}\left(L_{0},M\right)\times\Omega_{1}^{n}\left(L_{0}\right)$.
That is,\begin{gather*}
\Omega_{1}^{n}\left(S_{0}\right)  :=  \left\{ \eta\in\Omega^{n}\left(S_{0}\right)\mid\int_{S_{0}}\eta=1\right\} ,\\
\Omega_{0}^{n}\left(S_{0}\right)  :=  \left\{ \vartheta\in\Omega^{n}\left(S_{0}\right)\mid\int_{S_{0}}\vartheta=0\right\} ,\\
\text{Lag}_{w}\left(L_{0},M\right)  :=  \text{Lag}\left(L_{0},M\right)\times\Omega_{1}^{n}\left(L_{0}\right).
\end{gather*}

Integration along $L_{0}$ def\/ines a continuous linear functional
$\int_{L_{0}}:\Omega^{n}\left(L_{0}\right)\to\mathbb{R}$ on the convenient
vector space $\Omega^{n}\left(L_{0}\right)=\Gamma\left(\bigwedge^{n}T^{*}L_{0}\right)$.
So the kernel $\Omega_{0}^{n}\left(L_{0}\right)$ is a $c^{\infty}$-closed
(convenient) subspace. The space $\Omega_{1}^{n}\left(L_{0}\right)$
is an af\/f\/ine translation of $\Omega_{0}^{n}\left(L_{0}\right)$, which
means it is a smooth manifold modeled on $\Omega_{0}^{n}\left(L_{0}\right)$.
So $\text{Lag}_{w}\left(L_{0},M\right)$ is a smooth manifold modeled
on the space $Z^{1}\left(L_{0}\right)\oplus\mathfrak{X}\left(L_{0}\right)\oplus\Omega_{0}^{n}\left(L_{0}\right)$
with the product atlas $\left(U_{\left(i,\eta\right)},\varphi_{\left(i,\eta\right)}\right)_{\left(i,\eta\right)\in\text{Lag}_{w}\left(L_{0},M\right)}$.
That is, if $\left(U_{i},\varphi_{i}\right)_{i\in\text{Lag}\left(L_{0},M\right)}$
is the atlas on $\text{Lag}\left(L_{0},M\right)$ def\/ined in Proposition~\ref{sub:PropLag(L_0,M)} then the charts $\left(U_{\left(i,\eta\right)},\varphi_{\left(i,\eta\right)}\right)$
are def\/ined by
\begin{gather}
U_{(i,\eta)}  :=  U_{i}\times(\eta+\Omega_{0}^{n}(L_{0})),\label{eq:Lag_w(L_0,M)Chart1}\\
\varphi_{(i,\eta)}  : \  U_{(i,\eta)}\to Z^{1}(L_{0})\oplus\mathfrak{X}(L_{0})\oplus\Omega_{0}^{n}(L_{0}),\nonumber \\
\varphi_{(i,\eta)}(l,\eta+\vartheta)  := (\varphi_{i}(l),\vartheta)=(\alpha,Y,\vartheta).\nonumber
\end{gather}
 This atlas and the subspace $B^{1}(L_{0})\oplus\mathfrak{X}(L_{0})\oplus\Omega_{0}^{n}(L_{0})$
def\/ine a foliation $\mathcal{E}_{w}$ on $\text{Lag}_{w}(L_{0},M)$.

\begin{definition}
We will call the foliation $\mathcal{E}_{w}$ the \emph{isodrastic
foliation} of $\text{Lag}_{w}\left(L_{0},M\right)$ and a leaf of
$\mathcal{E}_{w}$ will be called an \emph{isodrast} in $\text{Lag}_{w}\left(L_{0},M\right)$.
\end{definition}

For each $\left(i,\eta\right)\in\text{Lag}_{w}\left(L_{0},M\right)$,
the tangent space $T_{\left(i,\eta\right)}\text{Lag}_{w}\left(L_{0},M\right)$
equals the vector space $\Gamma_{\text{closed}}\left(i^{*}TM\right)\oplus\Omega_{0}^{n}\left(L_{0}\right)$.
Meanwhile, the tangent space to an isodrast $\text{\textbf{H}}_{w}$
in $\text{Lag}_{w}\left(L_{0},M\right)$ at $\left(i,\eta\right)$
is given by $T_{\left(i,\eta\right)}\text{\textbf{H}}_{w}=\Gamma_{\text{exact}}\left(i^{*}TM\right)\oplus\Omega_{0}^{n}\left(L_{0}\right)$.
To each point $\left(i,\eta\right)\in\text{\textbf{H}}_{w}$ in an
isodrast we assign a skew-symmetric bilinear form on $T_{\left(i,\eta\right)}\text{\textbf{H}}_{w}$
via the expression\begin{equation}
\Omega_{\left(i,\eta\right)}\left(\left(X_{1},\vartheta_{1}\right),\left(X_{2},\vartheta_{2}\right)\right)
:=\int_{L_{0}}[\omega(X_{1},X_{2})\eta+h_{1}\vartheta_{2}-h_{2}\vartheta_{1}]\label{eq:OmegaLag_w(L_0,M)}\end{equation}
 where $\alpha_{X_{k}}=dh_{k}$ for some $h_{k}\in C^{\infty}\left(L_{0},\mathbb{R}\right)$.
This assignment does not depend on the choice of primitives $h_{k}$
since the top degree forms $\vartheta_{k}$ integrate to 0. Equivalently,
the pointwise assignment in~\eqref{eq:OmegaLag_w(L_0,M)} can be viewed
as an assignment on the smooth manifold $\text{Lag}_{w}\left(L_{0},M\right)^{\mathcal{E}_{w}}$.

\begin{proposition}\label{sub:PropOmegeH_w}
The pointwise assignment in \eqref{eq:OmegaLag_w(L_0,M)} on
$\emph{Lag}_{w}\left(L_{0},M\right)^{\mathcal{E}_{w}}$ defines
a basic differen\-tial $2$-form $\Omega$ with respect to the action
of $\emph{Dif\/f}_{+}(L_{0})$ on $\emph{Lag}_{w}(L_{0},M)^{\mathcal{E}_{w}}$
given by $a\cdot(i,\eta):=(i\circ a^{-1},(a^{-1})^{*}\eta)$.
\end{proposition}

\begin{proof}
We will f\/irst show that the assignment $\Omega$ def\/ines
a dif\/ferential 2-form on each leaf $\text{\textbf{H}}_{w}$ of~$\mathcal{E}_{w}$.
The assignment in \eqref{eq:OmegaLag_w(L_0,M)} def\/ines a map $\Omega:\text{\textbf{H}}_{w}\to L_{\text{alt}}^{2}(T\text{\textbf{H}}_{w},\text{\textbf{H}}_{w}\times\mathbb{R})$.
To check that this map is smooth, it is enough to check it in each
chart. If $\left(U_{\left(i,\eta\right)},\varphi_{\left(i,\eta\right)}\right)$
denotes a chart on $\text{\textbf{H}}_{w}$ then $\Omega$ def\/ines
a map from $U_{\left(i,\eta\right)}$ to $L_{\text{alt}}^{2}(\Gamma_{\text{exact}}\left(i^{*}TM\right)\oplus\Omega_{0}^{n}\left(L_{0}\right),\mathbb{R})$
(after $B^{1}\left(L_{0}\right)\oplus\mathfrak{X}\left(L_{0}\right)\oplus\Omega_{0}^{n}\left(L_{0}\right)$
has been identif\/ied with $\Gamma_{\text{exact}}\left(i^{*}TM\right)\times\Omega_{0}^{n}\left(L_{0}\right)$
via the splitting map $s_{i}:Z^{1}\left(L_{0}\right)\to\Gamma_{\text{closed}}\left(i^{*}TM\right)$
(see Remark \ref{sub:SplitMap})). This map is smooth if it maps smooth
curves in $U_{\left(i,\eta\right)}$ to smooth curves in $L_{\text{alt}}^{2}(\Gamma_{\text{exact}}\left(i^{*}TM\right)\oplus\Omega_{0}^{n}\left(L_{0}\right),\mathbb{R})$.
A curve in $L_{\text{alt}}^{2}(\Gamma_{\text{exact}}\left(i^{*}TM\right)\oplus\Omega_{0}^{n}\left(L_{0}\right),\mathbb{R})$
is smooth if the induced map $\mathbb{R}\times\left(\Gamma_{\text{exact}}\left(i^{*}TM\right)\oplus\Omega_{0}^{n}\left(L_{0}\right)\right)^{2}\to\mathbb{R}$
is smooth. Thus to verify that $\Omega$ is smooth, it is enough to
check the following statement: If $\left(M,\omega\right)=\left(T^{*}L_{0},\omega_{T^{*}L_{0}}\right)$,
$i:L_{0}\hookrightarrow T^{*}L_{0}$ denotes the zero section inclusion,
$\left(\alpha_{t}\circ a_{t},\eta_{t}\right)$ is a smooth curve in
$\text{\textbf{H}}_{w}$, $\left(X_{1}\left(t\right),\vartheta_{1}\left(t\right)\right)$
and $\left(X_{2}\left(t\right),\vartheta_{2}\left(t\right)\right)$
are smooth curves in $\Gamma_{\text{exact}}(i^{*}T(T^{*}L_{0}))\oplus\Omega_{0}^{n}\left(L_{0}\right)$
satisfying $\alpha_{X_{k}\left(t\right)}=dh_{k}\left(t\right)$ for
smooth curves $h_{k}\left(t\right)$ in $C^{\infty}\left(L_{0},\mathbb{R}\right)$,
$Z_{1}\left(t\right)$ and $Z_{2}\left(t\right)$ are the unique time
dependent vector f\/ields on $T^{*}L_{0}$ satisfying $\imath\left(Z_{k}\left(t\right)\right)\omega=\pi^{*}\alpha_{X_{k}\left(t\right)}$
for all $t$, $Y_{k}\left(t\right):=X_{k}\left(t\right)-\left.Z_{k}\left(t\right)\right|_{L_{0}}\in\mathfrak{X}\left(L_{0}\right)$,
and $s\in C^{\infty}\left(\mathbb{R},\mathbb{R}\right)$ is a smooth
map, then the map
\begin{gather*}
\mathbb{R}  \to  \mathbb{R},\\
t  \mapsto  \int_{L_{0}}\Big[\omega\Big(Z_{1}(t)\circ\alpha_{s(t)}\circ a_{s(t)}+\text{d}\alpha_{s(t)}\cdot Y_{1}(t)_{a_{s(t)}},  \\
\qquad\qquad Z_{2} (t )\circ\alpha_{s (t )}\circ a_{s (t )}+\text{d}\alpha_{s (t )}\cdot Y_{2} (t )_{a_{s (t )}}\Big)\eta_{s(t)}
+h_{1}(t)\vartheta_{2}(t)-h_{2}(t)\vartheta_{1}(t)\Big]
\end{gather*}
 is smooth. Since this statement follows from the smoothness of all
quantities in the integral, $\Omega$~is indeed a section of $L_{\text{alt}}^{2}\left(T\text{\textbf{H}}_{w},\text{\textbf{H}}_{w}\times\mathbb{R}\right)\to\text{\textbf{H}}_{w}$.

We will now show that $\Omega$ is basic with respect to the action
of $\text{Dif\/f}_{+}\left(L_{0}\right)$ on $\text{Lag}_{w}\left(L_{0},M\right)^{\mathcal{E}_{w}}$.
If $\left(i_{t},\eta_{t}\right)$ is a smooth curve in $\text{Lag}_{w}\left(L_{0},M\right)^{\mathcal{E}_{w}}$
with time derivative $\left(X,\vartheta\right)\in\Gamma_{\text{exact}}\left(i^{*}TM\right)\oplus\Omega_{0}^{n}\left(L_{0}\right)$
at $t=0$, then the tangent vector to the curve $l_{a}\left(i_{t},\eta_{t}\right)=\left(i_{t}\circ a^{-1},\left(a^{-1}\right)^{*}\eta_{t}\right)$
at $t=0$ is given by $\left(X\circ a^{-1},\left(a^{-1}\right)^{*}\vartheta\right)$.
Thus,\begin{gather*}
l_{a}^{*}\Omega_{\left(i,\eta\right)}\left(\left(X_{1},\vartheta_{1}\right),\left(X_{2},\vartheta_{2}\right)\right)   = \Bigg[\int_{L_{0}}\omega\left(X_{1}\circ a^{-1},X_{2}\circ a^{-1}\right)\left(a^{-1}\right)^{*}\eta \\
\phantom{l_{a}^{*}\Omega_{\left(i,\eta\right)}\left(\left(X_{1},\vartheta_{1}\right),\left(X_{2},\vartheta_{2}\right)\right)   =}{}
 +h_{1}\circ a^{-1}\left(a^{-1}\right)\vartheta_{2}-h_{2}\circ a^{-1}\left(a^{-1}\right)^{*}\vartheta_{1}\Bigg]\\
\phantom{l_{a}^{*}\Omega_{\left(i,\eta\right)}\left(\left(X_{1},\vartheta_{1}\right),\left(X_{2},\vartheta_{2}\right)\right)   }{}  =   \Omega_{\left(i,\eta\right)}\left(\left(X_{1},\vartheta_{1}\right),\left(X_{2},\vartheta_{2}\right)\right),
  \end{gather*}
 which means $\Omega$ is $\text{Dif\/f}_{+}\left(L_{0}\right)$-invariant.
To check that $\Omega$ is horizontal, let $Y\in\mathfrak{X}\left(L_{0}\right)$
be in the Lie algebra of $\text{Dif\/f}_{+}\left(L_{0}\right)$. If
$a_{t}$ is a smooth curve in $\text{Dif\/f}_{+}\left(L_{0}\right)$
through the identity map with time derivative $Y$ at $t=0$, then
the generating vector f\/ield at a point $\left(i,\eta\right)\in\text{Lag}_{w}\left(L_{0},M\right)$
is given by
\[
Y_{\text{Lag}_{w}\left(L_{0},M\right)}\left(i,\eta\right)=\left.\frac{d}{dt}\right|_{t=0}\big(i\circ a_{t}^{-1},\big(a_{t}^{-1}\big)^{*}\eta\big)=\left(-Y,-\mathcal{L}_{Y}\eta\right).
\]
 Since\[
\Omega_{\left(i,\eta\right)}\left(\left(Y,\mathcal{L}_{Y}\eta\right),
\left(X_{2},\vartheta_{2}\right)\right)=\int_{L_{0}}\left[\omega\left(Y,X_{2}\right)\eta-h_{2}\mathcal{L}_{Y}
\eta\right]=0\]
 (after integrating by parts, $\int_{L_{0}}\mathcal{L}_{Y}h_{2}\eta=-\int_{L_{0}}h_{2}\mathcal{L}_{Y}\eta$),
we conclude that $\Omega$ is also horizontal and thus basic.
\end{proof}

The quotient $\text{Lag}_{w}\left(L_{0},M\right)/\text{Dif\/f}_{+}\left(L_{0}\right)$
is naturally identif\/ied with the set of weighted Lagrangian submanifolds.
Explicitly, the identif\/ication of the quotient $\text{Lag}_{w}\left(L_{0},M\right)/\text{Dif\/f}_{+}\left(L_{0}\right)$
with the set of pairs $\left(L,\rho\right)$ is via the correspondence
$\left[\left(i,\eta\right)\right]\leftrightarrow\left(L,\rho\right)$
where $L=i\left(L_{0}\right)$ and $i^{*}\rho=\eta$.

\subsection*{Notation}

Set\[
\text{Lag}_{w}\left(M\right):=\text{Lag}_{w}\left(L_{0},M\right)/\text{Dif\/f}_{+}\left(L_{0}\right).\]

\begin{proposition}\label{sub:PropLag_w(M)SmoothStr}
$\emph{Lag}_{w}\left(M\right)$ is a smooth manifold modeled
on spaces $Z^{1}\left(L\right)\oplus\Omega_{0}^{n}\left(L\right)$
for $L\in\emph{Lag}\left(M\right)$. For each representative
$\left(i,\eta\right)$ in the class $\left(L,\rho\right)\in\emph{Lag}_{w}\left(M\right)$,\[
T_{\left(L,\rho\right)}\emph{Lag}_{w}\left(M\right)=\Gamma_{\emph{closed}}\left(i^{*}TM\right)\oplus\Omega_{0}^{n}\left(L_{0}\right)/\{\left(Y,\mathcal{L}_{Y}\eta\right)\mid Y\in\mathfrak{X}\left(L_{0}\right)\}.\]
\end{proposition}

\begin{proof} For each $\left(L,\rho\right)\in\text{Lag}_{w}\left(M\right)$,
by Theorem~\ref{sub:TheoremWLN} there exists a symplectic embedding $\lambda_{\left(L,\rho\right)}:M\supset W_{\left(L,\rho\right)}\to T^{*}L$
of a neighbourhood of $L$ onto a neighbourhood of the zero section
in the cotangent bundle. If $\pi_{T^{*}L}:T^{*}L\to L$ denotes the
cotangent bundle projection, then the restriction of $\pi_{T^{*}L}$
to the graph of any 1-form $\alpha\in\Omega^{1}\left(L\right)$ in
$T^{*}L$ def\/ines a dif\/feomorphism of that graph onto $L$. Def\/ine
\begin{gather*}
U_{\left(L,\rho\right)}:   =   \Big\{ \left(N,\sigma\right)\in\text{Lag}_{w}\left(M\right)\mid N\subset W_{\left(L,\rho\right)},\,\lambda_{\left(L,\rho\right)}\left(N\right)=\text{Graph}\left(\alpha\right), \\
\qquad\qquad \ \alpha\in Z^{1}\left(L\right),\,\sigma=\lambda_{\left(L,\rho\right)}^{*}(\left.\pi_{T^{*}L}
\right|_{\text{Graph}\left(\alpha\right)})^{*}\left(\rho+\theta\right),\,\theta\in\Omega_{0}^{n}\left(L\right)\Big\} ,\\
\varphi_{\left(L,\rho\right)}   :  \  U_{\left(L,\rho\right)}\to Z^{1}\left(L\right)\oplus\Omega_{0}^{n}\left(L\right),\\
\varphi_{\left(L,\rho\right)}\left(N,\sigma\right)   :=   \left(\alpha,\theta\right).
\end{gather*}
 All chart changings are smooth again by the smoothness of the symplectic
embeddings $\lambda_{L}$, so the collection $\left(U_{\left(L,\rho\right)},\varphi_{\left(L,\rho\right)}\right)_{\left(L,\rho\right)\in\text{Lag}_{w}\left(M\right)}$
def\/ines a smooth atlas on $\text{Lag}_{w}\left(M\right)$.

We will now describe the identif\/ication of tangent spaces of $\text{Lag}_{w}(M)\!$
with spaces $\Gamma_{\text{closed}}(i^{*}TM\!)\!\!$ $\oplus\Omega_{0}^{n}(L_{0})/\{(Y,\mathcal{L}_{Y}\eta)\mid Y\in\mathfrak{X}\left(L_{0}\right)\}$.
Let $\left(i,\eta\right)\in\text{Lag}_{w}\left(L_{0},M\right)$ be
a representative in the class $\left(L,\rho\right)\in\text{Lag}\left(M\right)$.
Let $\lambda_{i}:T^{*}L_{0}\supset W_{i}\to M$ denote the symplectic
embedding chosen in the def\/inition of the chart $\left(U_{i},\varphi_{i}\right)$
on $\text{Lag}_{w}\left(L_{0},M\right)$, and $s_{i}:Z^{1}\left(L_{0}\right)\to\Gamma_{\text{closed}}\left(i^{*}TM\right)$
the corresponding splitting map. Then the linear map
\begin{gather*}
T_{\left[\left(i,\eta\right)\right]}\text{Lag}_{w}\left(M\right)   \to   \Gamma_{\text{closed}}\left(i^{*}TM\right)\times\Omega_{0}^{n}\left(L_{0}\right)/\{\left(Y,\mathcal{L}_{Y}\eta\right)\mid Y\in\mathfrak{X}\left(L_{0}\right)\},\\
\left(\alpha,\theta\right)   \mapsto
 \left[\left(s_{i}\left(i^{*}\alpha\right),\theta\right)\right]
 \end{gather*}
 is a vector space isomorphism.
\end{proof}

The canonical projection
$\text{Lag}_{w}\left(M\right)\to\text{Lag}\left(M\right)$, which
forgets the weightings, pulls back $\mathcal{F}$ to a foliation $\mathcal{F}_{w}$
on $\text{Lag}_{w}\left(M\right)$. That is, the collection of subspaces
$\left\{ B^{1}\left(L\right)\oplus\Omega_{0}^{n}\left(L\right)\right\} $
and atlas $\left(U_{\left(L,\rho\right)},\varphi_{\left(L,\rho\right)}\right)$
indexed by ${\left(L,\rho\right)\in\text{Lag}_{w}\left(M\right)}$
def\/ine a foliation $\mathcal{F}_{w}$ on $\text{Lag}_{w}\left(M\right)$.

\begin{definition}
The foliation $\mathcal{F}_{w}$ will be called the \emph{isodrastic
foliation} of $\text{Lag}_{w}\left(M\right)$ and a leaf $\text{\textbf{I}}_{w}$
of $\mathcal{F}_{w}$ will be called an \emph{isodrast} in $\text{Lag}_{w}\left(M\right)$.
\end{definition}

\begin{proposition}\label{sub:PropOmegaSymplectic}
The smooth manifold $\emph{Lag}_{w}\left(L_{0},M\right)^{\mathcal{E}_{w}}$
is the total space of a principal $\emph{Dif\/f}_{+}\left(L_{0}\right)$
bundle over $\emph{Lag}_{w}\left(M\right)^{\mathcal{F}_{w}}$.
The basic $2$-form $\Omega$ on $\emph{Lag}_{w}\left(L_{0},M\right)^{\mathcal{E}_{w}}$
descends to a weak symplectic structure on $\emph{Lag}_{w}\left(M\right)^{\mathcal{F}_{w}}$.
Thus each isodrast in $\emph{Lag}_{w}\left(M\right)$ is
a weakly symplectic manifold.
\end{proposition}

\begin{proof}
We begin by describing
a f\/iber bundle atlas. For $\left(i,\eta\right)\in\text{Lag}_{w}\left(L_{0},M\right)^{\mathcal{E}_{w}}$,
let $\lambda_{\left(i,\eta\right)}:T^{*}L_{0}\supset W_{i}\to M$
denote the symplectic embedding chosen in def\/ining the chart $\left(U_{\left(i,\eta\right)},\varphi_{\left(i,\eta\right)}\right)$.
Def\/ine
\begin{gather*}
U_{\left[\left(i,\eta\right)\right]}:   =   \Big\{ \left(N,\sigma\right)\in\text{Lag}_{w}\left(M\right)\mid N=\lambda_{\left(i,\eta\right)}\left(\text{Graph}\left(\alpha\right)\right), \\
 \phantom{U_{\left[\left(i,\eta\right)\right]}:   =   \Big\{ }{} \alpha\in B^{1}\left(L_{0}\right),\,\sigma=\left.\lambda_{\left(i,\eta\right)}^{-1}\right|_{N}^{*}\pi^{*}
 \left(i^{*}\rho+\vartheta\right),
\; \vartheta\in\Omega_{0}^{n}\left(L_{0}\right)\Big\} .
 \end{gather*}
 If $p:\text{Lag}_{w}\left(L_{0},M\right)^{\mathcal{E}_{w}}\to\text{Lag}_{w}\left(M\right)^{\mathcal{F}_{w}}$
denotes the projection to the quotient, then
\begin{gather*}
p^{-1}\left(U_{\left[\left(i,\eta\right)\right]}\right)   =   \Big\{ \left(l,\kappa\right)\in\text{Lag}_{w}\left(L_{0},M\right)^{\mathcal{E}_{w}}\, |\, l=\lambda_{\left(i,\eta\right)}\circ\alpha\circ a,\,\alpha\in B^{1}\left(L_{0}\right), \\
\phantom{p^{-1}\left(U_{\left[\left(i,\eta\right)\right]}\right)   =   \Big\{}{}
  \alpha\left(L_{0}\right)\subset W_{i},\, a\in\text{Dif\/f}_{+}\left(L_{0}\right)\Big\} .\end{gather*}
 Def\/ine $\psi_{\left[\left(i,\eta\right)\right]}:p^{-1}\left(U_{\left[\left(i,\eta\right)\right]}\right)\to U_{\left[\left(i,\eta\right)\right]}\times\text{Dif\/f}_{+}\left(L_{0}\right)$
by\[
\psi_{\left[\left(i,\eta\right)\right]}\left(\left(\lambda_{\left(i,\eta\right)}\circ\alpha\circ a,\kappa\right)\right):=\left(\left(\lambda_{\left(i,\eta\right)}\left(\text{Graph}\left(\alpha\right)\right),\left.\lambda_{\left(i,\eta\right)}^{-1}\right|_{N}^{*}\pi^{*}\left(a^{-1}\right)^{*}\kappa\right),a\right)\]
 so that the collection $\left(U_{\left[\left(i,\eta\right)\right]},\psi_{\left[\left(i,\eta\right)\right]}\right)_{\left(i,\eta\right)\in\text{Lag}_{w}\left(L_{0},M\right)}$
def\/ines a f\/iber bundle atlas.

If $\left(N,\sigma\right)\in U_{\left[\left(i,\eta\right)\right]\left[\left(j,\nu\right)\right]}$,
and $\psi_{\left[\left(j,\nu\right)\right]}^{-1}\left(\left(N,\sigma\right),a\right)=\left(\lambda_{j}\circ\alpha\circ a,\kappa\right)$
then\[
\psi_{\left[\left(i,\eta\right)\right]\left[\left(j,\nu\right)\right]}\left(\left(N,\sigma\right),a\right)=\left(\left(N,\sigma\right),\pi\circ\lambda_{\left(i,\eta\right)}^{-1}\circ\lambda_{\left(j,\nu\right)}\circ\alpha\circ a\right).\]
 So if $\left(N,\sigma\right)\in U_{\left[\left(i,\eta\right)\right]\left[\left(j,\nu\right)\right]}$
and $N=\lambda_{j}\left(\text{Graph}\left(\alpha\right)\right)$ then
def\/ine $u_{\left[i\right]\left[j\right]}:U_{\left[\left(i,\eta\right)\right]\left[\left(j,\nu\right)\right]}\to\text{Dif\/f}_{+}\left(L_{0}\right)$
by\[
u_{\left[i\right]\left[j\right]}\left(\left(N,\sigma\right)\right):=\pi\circ\lambda_{\left(i,\eta\right)}^{-1}\circ\lambda_{\left(j,\nu\right)}\circ\alpha.\]
 It follows that if $\left(N,\sigma\right){\in} U_{\left[\left(i,\eta\right)\right]}\cap U_{\left[\left(j,\nu\right)\right]}\cap U_{\left[\left(k,\mu\right)\right]}$,
$N{=}\lambda_{\left(j,\nu\right)}\left(\text{Graph}\left(\alpha\right)\right)$,
and $N{=}\lambda_{\left(k,\mu\right)}\left(\text{Graph}\left(\beta\right)\right)$
then\begin{gather*}
u_{\left[i\right]\left[j\right]}
\left(\left(N,\sigma\right)\right)u_{\left[j\right]\left[k\right]}\left(\left(N,\sigma\right)\right)   =   \pi\circ\lambda_{i}^{-1}\circ\lambda_{j}\circ\alpha\circ\pi\circ\lambda_{j}^{-1}\circ\lambda_{k}\circ\beta\\
\phantom{u_{\left[i\right]\left[j\right]}
\left(\left(N,\sigma\right)\right)u_{\left[j\right]\left[k\right]}\left(\left(N,\sigma\right)\right)}{}   =  \pi\circ\lambda_{i}^{-1}\circ\lambda_{k}\circ\beta
   =   u_{\left[i\right]\left[k\right]}\left(\left(N,\sigma\right)\right),\\
u_{\left[i\right]\left[i\right]}\left(\left(N,\sigma\right)\right)=\text{Id}_{L_{0}},\\
\psi_{\left[\left(i,\eta\right)\right]\left[\left(j,\nu\right)\right]}\left(\left(N,\sigma\right),a\right)=u_{\left[\left(i,\eta\right)\right]\left[\left(j,\nu\right)\right]}\circ a.
\end{gather*}
 Hence, the f\/iber bundle atlas $(U_{[(i,\eta)]},\!\psi_{[(i,\eta)]})_{(i,\eta){\in}
 \text{Lag}_{w}(L_{0},M)}$
and the collection of maps $u_{[(i,\eta)][(j,\nu)]}{:}\!$
$U_{\left[\left(i,\eta\right)\right]\left[\left(j,\nu\right)\right]}\to\text{Dif\/f}_{+}\left(L_{0}\right)$
def\/ine a principal $\text{Dif\/f}_{+}\left(L_{0}\right)$ structure.

Since $\Omega$ def\/ines a basic 2-form on the total space $\text{Lag}_{w}\left(L_{0},M\right)^{\mathcal{E}_{w}}$,
it descends to a dif\/ferential 2-form (also denoted $\Omega$) on $\text{Lag}_{w}\left(M\right)^{\mathcal{F}_{w}}$
(see Section~\ref{sub:PrincipalBundles}).

We will now check closedness of $\Omega$ locally in a chart $\left(U_{\left(L,\rho\right)},\varphi_{\left(L,\rho\right)}\right)$.
On $U_{\left(L,\rho\right)}$ tangent vectors can be identif\/ied with
pairs $\left(Z,\vartheta\right)$ where $Z\in\mathfrak{X}\left(T^{*}L_{0}\right)$
is a vector f\/ield on the cotangent bundle satisfying $\imath\left(Z\right)\omega=\pi^{*}dh$
for $h\in C^{\infty}\left(L_{0},\mathbb{R}\right)$, and $\vartheta\in\Omega_{0}^{n}\left(L_{0}\right)$.
Under such an identif\/ication, if $\left(i,\eta\right)$ is a representative
in the class $\left(L,\rho\right)\in\text{Lag}_{w}\left(M\right)$,
and if we identify $M$ with $T^{*}L_{0}$ using the symplectic embedding
$\lambda_{i}:T^{*}L_{0}\supset W_{i}\to M$, then terms like $\omega\left(X_{1},X_{2}\right)\eta=i^{*}\left[\omega\left(Z_{1},Z_{2}\right)\right]\eta$
in the expression for $\Omega$ vanish since the $Z_{k}$'s are tangent
to the cotangent f\/ibers. So on $U_{\left(L,\rho\right)}$,\begin{equation}
\Omega_{\left(N,\sigma\right)}\left(\left(X_{1},\vartheta_{1}\right),\left(X_{2},\vartheta_{2}\right)\right)=\int_{L_{0}}\left(h_{1}\vartheta_{2}-h_{2}\vartheta_{1}\right).\label{eq:OmegaLocal}\end{equation}
 It follows that $d\Omega=0$ since locally $\Omega$ does not depend
on $\left(N,\sigma\right)$.

Finally, weak nondegeneracy follows from the local expression for
$\Omega$ in \eqref{eq:OmegaLocal} and the fact that the $h_{k}$'s
and $\vartheta_{k}$'s can be chosen independently. Indeed, if $h_{1}$
is nonzero on some open subset $V\subset L_{0}$, then we can take
$h_{2}$ to be zero and choose $\vartheta_{2}$ such that it is supported
on $V$ and $\int_{L_{0}}h_{1}\vartheta_{2}$ is nonzero. If $\vartheta_{1}$
is nonzero on an open subset $V\subset L_{0}$, then we can choose
$\vartheta_{2}$ to be zero and choose $h_{2}$ to be supported on
$V$ so that $\int_{L_{0}}h_{2}\vartheta_{1}$ is nonzero.
\end{proof}

$\text{Lag}_{w}\left(M\right)$ can also be described as the set of
equivalence classes $\left(L,\left[\widetilde{\rho}\right]\right)$
where $\widetilde{\rho}$ is an $n$-form on a neighbourhood of $L$
satisfying $\int_{L}\widetilde{\rho}=1$, and $\widetilde{\rho}_{1}\sim\widetilde{\rho}_{2}$
if and only if $\widetilde{\rho}_{1}$ and $\widetilde{\rho}_{2}$
have the same pullback to $L$. In \cite{Weinstein3} Weinstein used
this approach and heuristically viewed $\text{Lag}_{w}\left(M\right)$
and each leaf $\text{\textbf{I}}_{w}$ of $\mathcal{F}_{w}$ as inf\/inite
dimensional manifolds. He viewed tangent vectors as equivalence classes
of paths in $\text{Lag}_{w}(M)$ and $\text{\textbf{I}}_{w}$ to give
the following description of their tangent spaces and wrote down a
closed, weakly nondegenerate, skew-symmetric bilinear form $\Omega^{W}$
on each isodrast $\text{\textbf{I}}_{w}$:

\begin{theorem}[see Theorem~3.2 \& Lemma 3.3 in \cite{Weinstein3}]\label{Wtheorem} The tangent
space to $\emph{Lag}_{w}\left(M\right)$ at $\left(L,\rho\right)$
can be identified with the set of quadruples $\big(L,\widetilde{\rho},X,\widetilde{\theta}\,\big)$,
where $\widetilde{\rho}$ is an $n$-form on a~neighbourhood of $L$
satisfying $\int_{L}\widetilde{\rho}=1$, $X$ is a symplectic vector
field on a neighbourhood of~$L$, and~$\widetilde{\theta}$ is an
$n$-form on a neighbourhood of $L$ satisfying $\int_{L}L_{X}\widetilde{\rho}+\widetilde{\theta}=0$,
subject to the following equivalence relation. $\big(L,\widetilde{\rho}_{1},X_{1},\widetilde{\theta}_{1}\big)\sim
\big(L,\widetilde{\rho}_{2},X_{2},\widetilde{\theta}_{2}\big)$
if and only if the following conditions hold: $(1)$~$\widetilde{\rho}_{1}$
and $\widetilde{\rho}_{2}$ have the same pullback to $L$; $(2)$ $X_{1}-X_{2}$
is tangent to $L$; $(3)$ the pullbacks to $L$ of $L_{X_{1}}\widetilde{\rho}_{1}+\widetilde{\theta}_{1}$
and $L_{X_{1}}\widetilde{\rho}_{2}+\widetilde{\theta}_{2}$ are equal.

The tangent vectors to an isodrast $\emph{\textbf{I}}_{w}$
are represented by equivalence classes $\big[\big(L,\widetilde{\rho},X_{f},\widetilde{\theta}\,\big)\big]$
where $X_{f}$ is a Hamiltonian vector field on a neighbourhood of
$L$. $\emph{\textbf{I}}_{w}$ admits a closed, weakly nondegenerate,
skew-symmetric bilinear form $\Omega^{W}$ defined by\begin{equation}
\Omega_{\left(L,\rho\right)}^{W}\left(\xi_{1},\xi_{2}\right):=
\int_{L}\left[\left\{ f_{1},f_{2}\right\} \rho+f_{1}\big(L_{X_{f_{2}}}\widetilde{\rho}_{2}+\widetilde{\theta}_{2}\big)
-f_{2}\big(L_{X_{f_{1}}}\widetilde{\rho}_{1}+\widetilde{\theta}_{1}\big)\right].\label{eq:WeinsteinForm}
\end{equation}
\end{theorem}

We will show that this heuristic description of the tangent spaces
and bilinear structure $\Omega^{W}$ due to Weinstein can be derived
from the smooth structures on $\text{Lag}_{w}\left(M\right)$ and
$\text{\textbf{I}}_{w}$ def\/ined in Proposition~\ref{sub:PropLag_w(M)SmoothStr}
and from the weak symplectic structure $\Omega$ on $\text{\textbf{I}}_{w}$
(see \eqref{eq:OmegaLag_w(L_0,M)}).

\subsection*{Notation}

For $\left(L,\rho\right)\in\text{Lag}_{w}\left(M\right)$ let $Q_{\left(L,\rho\right)}^{\text{symp}}$
denote the space of equivalence classes $\big[\big(L,\widetilde{\rho},X,\widetilde{\theta}\,\big)\big]$
where $X$ is a symplectic vector f\/ield def\/ined on a neighbourhood
of $L$. Let $Q_{\left(L,\rho\right)}^{\text{ham}}$ denote the space
of equiva\-lence classes $\big[\big(L,\widetilde{\rho},X_{f},\widetilde{\theta}\,\big)\big]$
where $X_{f}$ is a Hamiltonian vector f\/ield def\/ined on a neighbourhood
of~$L$.

\begin{proposition}\label{sub:PropOmega^W=00003D00003DOmega}
For every representative $\left(i,\eta\right)\in\emph{Lag}_{w}\left(L_{0},M\right)$
in the class $\left(L,\rho\right)\in\emph{Lag}_{w}\left(M\right)$
there exists a vector space isomorphism
\begin{gather*}
\tau_{\left(i,\eta\right)}^{\emph{symp}}: \ T_{\left(L,\rho\right)}\emph{Lag}_{w}\left(M\right)   =  \Gamma_{\emph{closed}}\left(i^{*}TM\right)\times\Omega_{0}^{n}\left(L_{0}\right)/
 \left\{ \left(Y,\mathcal{L}_{Y}\eta\right)\mid Y\in\mathfrak{X}\left(L_{0}\right)\right\} \to Q_{\left(L,\rho\right)}^{\emph{symp}}.
 \end{gather*}
 For every representative $\left(i,\eta\right)\in\emph{Lag}_{w}\left(L_{0},M\right)$
in the class $\left(L,\rho\right)$ in an isodrast $\emph{\textbf{I}}_{w}$,
this restricts to a vector space isomorphism
\begin{gather*}
\tau_{\left(i,\eta\right)}^{\emph{ham}}: \ T_{\left(L,\rho\right)}\emph{\textbf{I}}_{w}   =  \Gamma_{\emph{exact}}\left(i^{*}TM\right)\times\Omega_{0}^{n}\left(L_{0}\right)/
 \left\{ \left(Y,\mathcal{L}_{Y}\eta\right)\mid Y\in\mathfrak{X}\left(L_{0}\right)\right\} \to Q_{\left(L,\rho\right)}^{\emph{ham}}.
 \end{gather*}
 Under the induced vector space identification $T_{\left(L,\rho\right)}\emph{\textbf{I}}_{w}\simeq Q_{\left(L,\rho\right)}^{\emph{ham}}$,
if $\zeta_{1},\zeta_{2}\in T_{\left(L,\rho\right)}\emph{\textbf{I}}_{w}$
correspond to $\xi_{1},\xi_{2}\in Q_{\left(L,\rho\right)}^{\emph{ham}}$
then $\Omega_{\left(L,\rho\right)}\left(\zeta_{1},\zeta_{2}\right)=\Omega_{\left(L,\rho\right)}^{W}\left(\xi_{1},\xi_{2}\right)$.
\end{proposition}

\begin{proof}
For $\left(L,\rho\right)=\left[\left(i,\eta\right)\right]\in\text{Lag}_{w}\left(M\right)$
let $\lambda_{\left(i,\eta\right)}:T^{*}L_{0}\supset W_{\left(i,\eta\right)}\to M$
denote the symplectic embedding chosen for the chart $\left(U_{\left(i,\eta\right)},\varphi_{\left(i,\eta\right)}\right)$
on $\text{Lag}_{w}\left(L_{0},M\right)$. Let $\left[\left(X,\vartheta\right)\right]\in T_{\left(L,\rho\right)}\text{Lag}_{w}\left(M\right)$
be a tangent vector, with representative $\left(X,\vartheta\right)$
such that $\text{d}\lambda_{\left(i,\eta\right)}^{-1}\circ X$ is
tangent to the cotangent f\/ibers in $T^{*}L_{0}$. Let $Z_{\alpha_{X}}$
be the unique vector f\/ield on $T^{*}L_{0}$ satisfying $\imath\left(Z_{\alpha_{X}}\right)\omega_{T^{*}L_{0}}=\pi^{*}\alpha_{X}$.
Set $\widetilde{\rho}:=\big(\lambda_{\left(i,\eta\right)}^{-1}\big)^{*}\pi^{*}i^{*}\rho$,
$\widetilde{Z}_{\alpha_{X}}:=\text{d}\lambda_{\left(i,\eta\right)}\left(Z_{\alpha_{X}}\right)$
and $\widetilde{\vartheta}:=\big(\lambda_{\left(i,\eta\right)}^{-1}\big)^{*}\pi^{*}\vartheta$,
and def\/ine
\begin{gather*}
\tau_{\left(L,\rho\right)}^{\text{symp}}: \ T_{\left(L,\rho\right)}\text{Lag}_{w}\left(M\right)   \to   Q_{\left(L,\rho\right)}^{\text{symp}}, \qquad
\tau_{\left(L,\rho\right)}^{\text{symp}}\left(\left[\left(X,\vartheta\right)\right]\right)   :=   \big[\big(L,\widetilde{\rho},\widetilde{Z}_{\alpha_{X}},\widetilde{\vartheta}\, \big)\big].
\end{gather*}
 The linear map $\tau_{\left(L,\rho\right)}^{\text{symp}}$ has an
inverse given by{\samepage \begin{gather*}
Q_{\left(L,\rho\right)}^{\text{symp}}   \to   T_{(L,\rho)}\text{Lag}_{w}(M), \qquad
\big[\big(L,\widetilde{\rho},\widetilde{X},\widetilde{\vartheta}\, \big)\big]   \mapsto   \big[\big(\widetilde{X}\circ i,i^{*}\big(L_{\widetilde{X}}\widetilde{\rho}+\widetilde{\vartheta}\, \big)\big)\big].
\end{gather*}
 The isomorphism $\tau_{\left(L,\rho\right)}^{\text{ham}}:T_{\left(L,\rho\right)}\text{\textbf{I}}_{w}\to Q_{\left(L,\rho\right)}^{\text{ham}}$
is described similarly.}

Finally if $\zeta_{1},\zeta_{2}\in T_{\left(L,\rho\right)}\text{\textbf{I}}_{w}$
with $\zeta_{k}=\left[\left(X_{k},\vartheta_{k}\right)\right]$, with
representatives $\left(X_{k},\vartheta_{k}\right)$ such that $\text{d}\lambda_{\left(i,\eta\right)}^{-1}\circ X_{k}$
is tangent to the cotangent f\/ibers in $T^{*}L_{0}$, then $\omega(X_{1},X_{2})=0$.
So if $\alpha_{X_{k}}=dh_{k}$, $\widetilde{h}_{k}:=\big(\lambda_{\left(i,\eta\right)}^{-1}\big)^{*}\pi^{*}h_{k}$,
and $\widetilde{\vartheta}_{k}:=\big(\lambda_{\left(i,\eta\right)}^{-1}\big)^{*}\pi^{*}\vartheta_{k}$
then\begin{gather*}
\Omega_{\left(L,\rho\right)}^{W}\big(\tau_{\left(L,\rho\right)}^{\text{ham}}
\left(\zeta_{1}\right),\tau_{\left(L,\rho\right)}^{\text{ham}}\left(\zeta_{2}\right)\big)   =   \int_{L}\big(\widetilde{h}_{1}\widetilde{\vartheta}_{2}-\widetilde{h}_{2}\widetilde{\vartheta}_{1}\big)
   =   \int_{L_{0}}\left(h_{1}\vartheta_{2}-h_{2}\vartheta_{1}\right)
   =   \Omega_{\left(L,\rho\right)}\left(\zeta_{1},\zeta_{2}\right).\!\!\!\!\!\!\!\tag*{\qed}
   \end{gather*}\renewcommand{\qed}{}
\end{proof}

\begin{example}\label{sub:ExampleSphere}
Let $M=S^{2}$ and $L_{0}=S^{1}$. Since $S^{1}$ is one dimensional,
all embeddings are Lagrangian and all 1-forms on $S^{1}$ are closed.
So $\text{Lag}\left(S^{1},S^{2}\right)=\text{Emb}\left(S^{1},S^{2}\right)$
and for every embedding $i$ we have that $T_{i}\text{Lag}\left(S^{1},S^{2}\right)=\Gamma\left(i^{*}TS^{2}\right)$.

For any point $\left[\left(i,\eta\right)\right]=\left(L,\rho\right)$
in a leaf $\text{\textbf{I}}_{w}$, if $j$ is a compatible almost
complex structure on~$S^{2}$, i.e. $g\left(\cdot,\cdot\right):=\omega_{S^{2}}\left(\cdot,j\cdot\right)$
def\/ines a Riemannian metric on $S^{2}$, then for every tangent vector
$\xi\in T_{\left[\left(i,\eta\right)\right]}\text{\textbf{I}}_{w}$
there exists a unique representative $\left(X,\vartheta\right)\in\xi$
with $X\left(x\right)\in jT_{i\left(x\right)}L$ for every $x\in L_{0}$.
For such choices of representatives the expression for $\Omega$ becomes\[
\Omega\left(\left[\left(X_{1},\vartheta_{1}\right)\right],\left[\left(X_{2},\vartheta_{2}\right)\right]\right)=\int_{S^{1}}\left(h_{1}\vartheta_{2}-h_{2}\vartheta_{1}\right)\]
 where $\alpha_{X_{k}}=dh_{k}$ for $h_{k}\in C^{\infty}\left(S^{1},\mathbb{R}\right)$
(because the subspaces $jT_{i\left(x\right)}L\subset T_{i\left(x\right)}S^{2}$
are Lagrangian).

Each $\vartheta_{k}$ can be written as $r_{k}\left(x\right)dx$ for
some function $r_{k}$ on $S^{1}$. Meanwhile, any function $f$ on
$S^{1}$ has a Fourier series expansion\[
f\left(x\right)=\sum_{n=-\infty}^{\infty}\widehat{f}\left(n\right)e^{in\pi x}\]
 which reduces the expression for $\Omega$ to\[
\Omega\left(\xi_{1},\xi_{2}\right)=\sum_{n=-\infty}^{\infty}\big(\widehat{h}_{1}\left(n\right)
\widehat{\vartheta}_{2}\left(n\right)-\widehat{h}_{2}\left(n\right)\widehat{\vartheta}_{1}\left(n\right)\big).\]
 This expression is a countably inf\/inite version of the standard symplectic
vector space structure.
\end{example}

\begin{remark}
Weinstein's original construction was more general than we have described
so far. It included the case of Lagrangian submanifolds which are
neither compact nor oriented. In this case Weinstein used compactly
supported densities instead of volume forms. All of our constructions also carry through for non-oriented Lagrangian submanifolds equipped with compactly supported densities.
\end{remark}

\begin{example}
Let $\left(M,\omega\right)=\left(\mathbb{R}^{2},dq\wedge dp\right)$
and $L_{0}=\mathbb{R}$. As in Ex. \ref{sub:ExampleSphere}, $\mathbb{R}$
is one dimensional which means $\text{Lag}\left(\mathbb{R},\mathbb{R}^{2}\right)=\text{Emb}\left(\mathbb{R},\mathbb{R}^{2}\right)$
and for every embedding $i$ we have that $T_{i}\text{Lag}\left(L_{0},M\right)=\Gamma\left(i^{*}TM\right)$.
Moreover, since $H^{1}\left(L_{0}\right)=0$ the leaves of $\mathcal{E}$
consist of path connected components in $\text{Emb}\left(\mathbb{R},\mathbb{R}^{2}\right)$.
Thus the leaves of $\mathcal{F}$ consist of oriented one dimensional
submanifolds in $\mathbb{R}^{2}$ dif\/feomorphic to $\mathbb{R}$.

Though $L_{0}$ is not compact, we can still use compactly supported
1-forms as weightings. A~leaf $\text{\textbf{I}}_{w}$ then of $\mathcal{F}_{w}$
consists of isotopic one dimensional submanifolds in $\mathbb{R}^{2}$
dif\/feomorphic to~$\mathbb{R}$, equipped with compactly supported
1-forms. Any tangent vector $X\in T_{i}\text{Lag}\left(\mathbb{R},\mathbb{R}^{2}\right)$
can be written in components as $X=\widetilde{q}\frac{\partial}{\partial q}+\widetilde{p}\frac{\partial}{\partial p}$.
Since any 1-form $\eta$ on $\mathbb{R}$ can be written as $\widetilde{\eta}\left(x\right)dx$
for some function $\widetilde{\eta}$, the expression for $\Omega$
on such a leaf $\text{\textbf{I}}_{w}$ becomes\[
\Omega\left(\xi_{1},\xi_{2}\right)=\int_{-\infty}^{\infty}\big(\widetilde{q}_{1}\left(x\right)\widetilde{p}_{2}
\left(x\right)\widetilde{\eta}-\widetilde{q}_{2}\left(x\right)\widetilde{p}_{1}\left(x\right)
\widetilde{\eta}+\widetilde{h}_{1}\widetilde{\vartheta}_{2}-\widetilde{h}_{2}\widetilde{\vartheta}_{1}\big)dx.\]
\end{example}

\section{Positive weighted Lagrangian submanifolds}

In this section we will consider an open subset $\text{Lag}_{pw}\left(M\right)$
of $\text{Lag}_{w}(M)$ consisting of Lagrangian submanifolds weighted
with volume forms. All constructions involving not necessarily positive
weightings from before carry over to this case. In particular there
is a foliation $\mathcal{F}_{pw}$ of~$\text{Lag}_{pw}\left(M\right)$
whose leaves have a weak symplectic structure. The space $\text{Lag}_{pw}\left(M\right)$
also has a dif\/ferent description. By f\/ixing a positive weighting $\eta_{0}$,
the space of positive weighted Lagrangian submanifolds can be identif\/ied
with the quotient of $\text{Lag}\left(L_{0},M\right)$ by the group
of dif\/feomorphisms that preserve~$\eta_{0}$. This identif\/ication
is $\text{Ham}\left(M\right)$ equivariant and makes use of Moser's
Lemma~\cite{Moser}.

Fix $L_{0}$ to be a compact oriented manifold and $\left(M,\omega\right)$
a symplectic manifold with $\dim L_{0}=\frac{1}{2}\dim M$ as before.

\subsection*{Notation}

Let $\text{Vol}_{1}\left(S_{0}\right)$ denote the set of volume forms
on a compact oriented manifold $S_{0}$ that integrate to 1 and $\text{Lag}_{pw}\left(L_{0},M\right)$
the product $\text{Lag}\left(L_{0},M\right)\times\text{Vol}_{1}\left(L_{0}\right)$.
That is,\begin{gather*}
\text{Vol}_{1}\left(S_{0}\right)  :=  \left\{ \eta\in\Omega^{n}\left(S_{0}\right)\mid\eta\text{ nowhere vanishing, }\int_{S_{0}}\eta=1\right\} ,\\
\text{Lag}_{pw}\left(L_{0},M\right)  :=  \text{Lag}\left(L_{0},M\right)\times\text{Vol}_{1}\left(L_{0}\right).
\end{gather*}
 For every $\eta\in\text{Vol}_{1}\left(L_{0}\right)$, there exists
a $c^{\infty}$-open neighbourhood $U_{\eta}$ of $0\in\Omega_{0}^{n}\left(L_{0}\right)$
such that $\eta+U_{\eta}\subset\text{Vol}_{1}\left(L_{0}\right)$.
Thus $\text{Vol}_{1}\left(L_{0}\right)$ is a $c^{\infty}$-open subset
of the af\/f\/ine space $\Omega_{1}^{n}\left(L_{0}\right)$. It follows
that $\text{Lag}\left(L_{0},M\right)\times\text{Vol}_{1}\left(L_{0}\right)$
is a smooth manifold with an atlas $\left(U_{\left(i,\eta\right)},\varphi_{\left(i,\eta\right)}\right)_{\left(i,\eta\right)\in\text{Lag}_{pw}\left(L_{0},M\right)}$
given by \eqref{eq:Lag_w(L_0,M)Chart1} except that now\[
U_{\left(i,\eta\right)}=U_{i}\times\left(\eta+U_{\eta}\right).\]
 The atlas $\left(U_{\left(i,\eta\right)},\varphi_{\left(i,\eta\right)}\right)_{\left(i,\eta\right)\in\text{Lag}_{pw}\left(L_{0},M\right)}$
and the subspace $B^{1}\left(L_{0}\right)\oplus\mathfrak{X}\left(L_{0}\right)\oplus\Omega_{0}^{n}\left(L_{0}\right)$
def\/ine a~foliation $\mathcal{E}_{pw}$ on $\text{Lag}_{pw}\left(L_{0},M\right)$.

\begin{definition}
We will call $\mathcal{E}_{pw}$ the \emph{isodrastic foliation} of
$\text{Lag}_{pw}\left(L_{0},M\right)$ and a leaf of $\mathcal{E}_{pw}$
will be called an \emph{isodrast} in $\text{Lag}_{pw}\left(L_{0},M\right)$.
\end{definition}

The quotient of $\text{Lag}_{pw}\left(L_{0},M\right)$ by the $\text{Dif\/f}_{+}\left(L_{0}\right)$
action\[
a\cdot\left(i,\eta\right):=\big(i\circ a^{-1},\left(a^{-1}\right)^{*}\eta\big)\]
 can be identif\/ied with the set of oriented Lagrangian submanifolds
in $M$ equipped with positive weightings.

\subsection*{Notation}

Set\[
\text{Lag}_{pw}\left(M\right):=\text{Lag}_{pw}\left(L_{0},M\right)/\text{Dif\/f}_{+}\left(L_{0}\right).\]

As in Proposition~\ref{sub:PropLag_w(M)SmoothStr}, the set $\text{Lag}_{pw}\left(M\right)$
admits an atlas of charts $\left(U_{\left(L,\rho\right)},\varphi_{\left(L,\rho\right)}\right)$
indexed by $\left(L,\rho\right)\in\text{Lag}_{pw}\left(M\right)$
modeled on spaces $Z^{1}\left(L\right)\oplus\Omega_{0}^{n}\left(L\right)$
for $L\in\text{Lag}\left(M\right)$. For each representative $\left(i,\eta\right)$
in the class $\left(L,\rho\right)\in\text{Lag}_{pw}\left(M\right)$
the tangent space to $\text{Lag}_{pw}\left(M\right)$ at $\left(L,\rho\right)$
is equal to $\Gamma_{\text{closed}}\left(i^{*}TM\right)\oplus\Omega_{0}^{n}\left(L_{0}\right)/\{\left(Y,\mathcal{L}_{Y}\eta\right)\mid Y\in\mathfrak{X}\left(L_{0}\right)\}$.
The canonical projection $\text{Lag}_{pw}\left(M\right)\to\text{Lag}\left(M\right):\,\left(L,\rho\right)\mapsto L$
pulls back the foliation $\mathcal{F}$ on $\text{Lag}\left(M\right)$
to a foliation $\mathcal{F}_{pw}$ on $\text{Lag}_{pw}\left(M\right)$.

\begin{definition}
We will call $\mathcal{F}_{pw}$ the \emph{isodrastic foliation} of
$\text{Lag}_{pw}\left(M\right)$ and a leaf of $\mathcal{F}_{pw}$
will be called an \emph{isodrast} in $\text{Lag}_{pw}\left(M\right)$.
\end{definition}

Using a result of Moser \cite{Moser}, we can describe isodrasts in
$\text{Lag}_{pw}\left(M\right)$ more explicitly. Moser's Lemma states
that if $\Lambda_{0}$ and $\Lambda_{1}$ are two volume forms on
a~compact manifold $N$, such that $\int_{N}\Lambda_{0}=\int_{N}\Lambda_{1}$,
then there exists an isotopy $\psi_{t}\in\text{Dif\/f}_{+}\left(N\right)$
satisfying $\psi_{1}^{*}\Lambda_{0}=\Lambda_{1}$. Thus a~positive
weighting on a Lagrangian submanifold~$L$ can be moved to any other
positive weighting via an isotopy of~$L$. Any such isotopy $\psi_{t}$
can be lifted to a~Hamiltonian isotopy of~$M$ in the following way.
Choose a symplectic embedding $\lambda:M\supset U\to T^{*}L$ of a
neighbourhood~$U$ of~$L$ onto a neighbourhood of the zero section
in $T^{*}L$. If $\psi_{t}^{\sharp}$ denotes the cotangent lift of
$\psi_{t}$, i.e. $\psi_{t}^{\sharp}\left(x,p\right):=\big(\psi_{t}\left(x\right),\left(\psi_{t}^{-1}\right)^{*}p\big)$,
then $\tau_{t}:=\lambda^{-1}\circ\psi_{t}^{\sharp}\circ\lambda$ extends
$\psi_{t}$ to a Hamiltonian isotopy of a~neighbourhood of~$L$. Let
$H_{t}$ be a time dependent family of Hamiltonians corresponding
to~$\tau_{t}$. If $\chi$ is a cutof\/f function supported near $L$,
and equal to 1 near $L$, then the family $\chi H_{t}$ generates
a~Hamiltonian isotopy of~$M$ which restricts to~$\psi_{t}$ on~$L$.
It follows that the isodrasts in $\text{Lag}_{pw}\left(M\right)$
are equal to the $\text{Ham}\left(M\right)$ orbits in $\text{Lag}_{pw}\left(M\right)$
under the action
\[
\phi\cdot\left[\left(i,\eta\right)\right]:=\left[\left(\phi\circ i,\eta\right)\right]\leftrightarrow\phi\cdot\left(L,\rho\right):=
\big(\phi\left(L\right),\left(\phi^{-1}\right)^{*}\rho\big).\]

We can also describe $\text{Lag}_{pw}\left(M\right)$ and each leaf
$\text{\textbf{I}}_{w}\subset\text{Lag}_{pw}\left(M\right)$ in a
dif\/ferent way. Suppose that $L_{0}$ is equipped with a f\/ixed volume
form $\eta_{0}$ that integrates to 1.

\subsection*{Notation}

Let $\text{Dif\/f}\left(S_{0},\eta_{0}\right)$ denote the group of
dif\/feomorphisms of a manifold $S_{0}$ that preserve a f\/ixed volume
form $\eta_{0}$, and $\mathfrak{X}\left(S_{0},\eta_{0}\right)$ the
set of divergence free vector f\/ields on $S_{0}$. That is,\begin{gather*}
\text{Dif\/f}\left(S_{0},\eta_{0}\right)   :=   \left\{ a\in\text{Dif\/f}\left(S_{0}\right)\mid a^{*}\eta_{0}=\eta_{0}\right\} ,\\
\mathfrak{X}\left(S_{0},\eta_{0}\right)   :=   \left\{ \xi\in\mathfrak{X}\left(S_{0}\right)\mid\mathcal{L}_{\xi}\eta_{0}=0\right\} .\end{gather*}

The space $\mathfrak{X}\left(S_{0},\eta_{0}\right)$ is the Lie algebra
of $\text{Dif\/f}\left(S_{0},\eta_{0}\right)$ (see Theorem~43.7 in \cite{KrieglMichor}).
Meanwhile, $\text{Dif\/f}\left(L_{0},\eta_{0}\right)$ acts freely on
$\text{Lag}\left(L_{0},M\right)$ via\[
a\cdot i:=i\circ a^{-1}.\]

\begin{proposition}\label{PWQuotient}
The quotient space $\emph{Lag}\left(L_{0},M\right)/\emph{Dif\/f}\left(L_{0},\eta_{0}\right)$
is a smooth manifold modeled on spaces $\Gamma_{\emph{closed}}\left(i^{*}TM\right)/\mathfrak{X}\left(L_{0},\eta_{0}\right)$
for $i\in\emph{Lag}\left(L_{0},M\right)$. Moreover,
there exists a $\emph{Ham}\left(M\right)$ equivariant diffeomorphism
$\upsilon:\emph{Lag}\left(L_{0},M\right)/\emph{Dif\/f}\left(L_{0},\eta_{0}\right)\to\emph{Lag}_{pw}\left(M\right)$.
\end{proposition}

\begin{proof}
Def\/ine\begin{gather*}
\upsilon: \ \text{Lag}\left(L_{0},M\right)/\text{Dif\/f}\left(L_{0},\eta_{0}\right)   \to  \text{Lag}_{pw}\left(M\right),
\qquad
\left[i\right]   \mapsto   \left[\left(i,\eta_{0}\right)\right].\end{gather*}
 Then $\upsilon$ is injective since\[
\left[\left(i_{1},\eta_{0}\right)\right]=\left[\left(i_{2},\eta_{0}\right)\right]\ \Rightarrow \ i_{2}=i_{1}\circ a^{-1},\quad \left(a^{-1}\right)^{*}\eta_{0}=\eta_{0}\]
 for some $a\in\text{Dif\/f}_{+}\left(L_{0}\right)$. To check surjectivity
suppose that $\left(L,\rho\right)=\left[\left(i,\eta\right)\right]\in\text{Lag}_{pw}\left(M\right)$
is a positive weighted Lagrangian submanifold. By Moser's Lemma, since
$\eta$ and $\eta_{0}$ are volume forms on $L_{0}$ that both induce
the orientation of $L_{0}$ and integrate to 1, there exists an isotopy
$\psi_{t}\in\text{Dif\/f}_{+}\left(L_{0}\right)$ such that $\psi_{1}^{*}\eta_{0}=\eta$.
Thus\[
\upsilon\left(\left[i\circ\psi_{1}^{-1}\right]\right)=\left[\left(i\circ\psi_{1}^{-1},\eta_{0}\right)\right]=\left[\left(i,\eta\right)\right].\]

We will now describe charts into spaces $\Gamma_{\text{closed}}\left(i^{*}TM\right)/\mathfrak{X}\left(L_{0},\eta_{0}\right)$.
Let $i\in\text{Lag}(L_{0},M)$ and let $\lambda_{i}:T^{*}L_{0}\supset W_{i}\to M$
be the symplectic embedding chosen in def\/ining the chart $\left(U_{i},\varphi_{i}\right)$
on $\text{Lag}\left(L_{0},M\right)$. Given a representative $X$
of a class $\left[X\right]\in\Gamma_{\text{closed}}\left(i^{*}TM\right)/\mathfrak{X}\left(L_{0},\eta_{0}\right)$,
the section $\text{d}\lambda_{i}^{-1}\circ X\in\Gamma\left(\left.T\left(T^{*}L_{0}\right)\right|_{L_{0}}\right)$
can be decomposed as $\text{d}\lambda_{i}^{-1}\circ X=\left.Z_{\alpha_{X}}\right|_{L_{0}}+Y$
where $\imath\left(Z_{\alpha_{X}}\right)\omega_{T^{*}L_{0}}=\pi^{*}\alpha_{X}$
and $Y\in\mathfrak{X}\left(L_{0}\right)$. For a dif\/ferent choice
of representative, this decomposition changes only in the component
$Y$ tangent to $L_{0}$. Thus this decomposition def\/ines a vector
space isomorphism\begin{gather*}
\zeta_{1}: \ \Gamma_{\text{closed}}\left(i^{*}TM\right)/\mathfrak{X}\left(L_{0},\eta_{0}\right)   \to   Z^{1}\left(L_{0}\right)\oplus\left(\mathfrak{X}\left(L_{0}\right)/\mathfrak{X}\left(L_{0},\eta_{0}\right)\right), \qquad
\left[X\right]   \mapsto   \left(\alpha_{X},\left[Y\right]\right),\end{gather*}
 where $\text{d}\lambda_{i}^{-1}\circ X=\left.Z_{\alpha_{X}}\right|_{L_{0}}+Y$.

Consider the map\begin{gather*}
\zeta_{2}: \ \mathfrak{X}\left(L_{0}\right)/\mathfrak{X}\left(L_{0},\eta_{0}\right)   \to   \Omega_{0}^{n}\left(L_{0}\right),\qquad
\left[Y\right]   \mapsto   \mathcal{L}_{Y}\eta_{0}.
\end{gather*}
 We claim that $\zeta_{2}$ is a vector space isomorphism. It is injective
since $\mathcal{L}_{Y_{1}}\eta_{0}=\mathcal{L}_{Y_{2}}\eta_{0}$ implies
that $Y_{1}-Y_{2}\in\mathfrak{X}\left(L_{0},\eta_{0}\right)$. To
check surjectivity, choose a metric $g_{0}$ on $L_{0}$ such that
the induced volume form $\mu\left(g_{0}\right)$ equals $\eta_{0}$.
Suppose that $\vartheta=\vartheta'\cdot\eta_{0}\in\Omega_{0}^{n}\left(L_{0}\right)$
for $\vartheta'\in C^{\infty}\left(L_{0},\mathbb{R}\right)$. By the
Hodge Decomposition Theorem (see e.g.~\cite{Warner}), there exists
a function $h'\in C^{\infty}\left(L_{0},\mathbb{R}\right)$ (unique
up to constants) such that $\triangle h'=\vartheta'$. For such an
$h'$, it follows that $\mathcal{L}_{\nabla h'}\eta_{0}=\triangle h'\cdot\eta_{0}=\vartheta$.

The isomorphisms $\zeta_{1}$ and $\zeta_{2}$ combine to def\/ine a
vector space isomorphism $\zeta$ from the quotient $\Gamma_{\text{closed}}\left(i^{*}TM\right)/\mathfrak{X}\left(L_{0},\eta_{0}\right)$
to $Z^{1}\left(L_{0}\right)\oplus\Omega_{0}^{n}\left(L_{0}\right)$.
Let $E_{U_{\left[i\right]}}$ denote $\Gamma_{\text{closed}}\left(i^{*}TM\right)/\mathfrak{X}\left(L_{0},\eta_{0}\right)$
equipped with the convenient structure induced by this isomorphism.
Def\/ine\begin{gather*}
U_{\left[i\right]}   :=   \Big\{ \left[l\right]\in\text{Lag}\left(L_{0},M\right)/\mathfrak{X}\left(L_{0},\eta_{0}\right)\mid\left[\left(l,\eta_{0}\right)\right]=\left(N,\sigma\right)\in U_{\left[\left(i,\eta_{0}\right)\right]} \\
\hphantom{U_{\left[i\right]}   :=   \Big\{ }{} N=\lambda_{\left(i,\eta\right)}\left(\text{Graph}\left(\alpha\right)\right),\,\alpha\in Z^{1}\left(L_{0}\right),\, \sigma=\left.\lambda_{\left(i,\eta\right)}^{-1}\right|_{N}^{*}\pi^{*}\left(i^{*}\rho+\vartheta\right),
 \,\vartheta\in\Omega_{0}^{n}\left(L_{0}\right)\Big\} ,\\
\varphi_{\left[i\right]}   : \  U_{\left[i\right]}\to E_{U_{\left[i\right]}},\qquad
\varphi_{\left[i\right]}\left(\left[l\right]\right)   :=   \zeta^{-1}\left(\alpha,\vartheta\right),
\end{gather*}
 where $\left(U_{\left[\left(i,\eta_{0}\right)\right]},\varphi_{\left[\left(i,\eta_{0}\right)\right]}\right)$
is one of the charts on $\text{Lag}_{pw}\left(M\right)$. The charts
$\left(U_{\left[i\right]},\varphi_{\left[i\right]}\right)$ indexed
by $i\in\text{Lag}\left(L_{0},M\right)$ def\/ine a smooth atlas on
$\text{Lag}\left(L_{0},M\right)/\mathfrak{X}\left(L_{0},\eta_{0}\right)$.

From the def\/inition of the smooth structure on $\text{Lag}\left(L_{0},M\right)/\mathfrak{X}\left(L_{0},\eta_{0}\right)$,
it follows that the map $\upsilon$ is a dif\/feomorphism. Since the
$\text{Ham}\left(M\right)$ action commutes with the $\text{Dif\/f}\left(L_{0},\eta_{0}\right)$
action on $\text{Lag}\left(L_{0},M\right)$ and with the $\text{Dif\/f}_{+}\left(L_{0}\right)$
action on $\text{Lag}_{pw}\left(L_{0},M\right)$, we have well def\/ined
$\text{Ham}\left(M\right)$ actions on the quotients. Thus for $\phi\in\text{Ham}\left(M\right)$,\[
\upsilon\left(\phi\cdot\left[i\right]\right)=\upsilon\left(\left[\phi\circ i\right]\right)=\left[\left(\phi\circ i,\eta_{0}\right)\right]=\phi\cdot\left[\left(i,\eta_{0}\right)\right]=\phi\cdot\upsilon\left(\left[i\right]\right)\]
 which verif\/ies the $\text{Ham}\left(M\right)$ equivariance of $\upsilon$.
\end{proof}

As in the case of not necessarily positive weighted Lagrangian
submanifolds, the smooth manifold $\text{Lag}_{pw}\left(L_{0},M\right)^{\mathcal{E}_{pw}}$
is the total space of a $\text{Dif\/f}_{+}\left(L_{0}\right)$ bundle
over $\text{Lag}_{pw}\left(M\right)^{\mathcal{F}_{pw}}$. We can def\/ine
a basic 2-form $\Omega$ on $\text{Lag}_{pw}\left(L_{0},M\right)^{\mathcal{E}_{pw}}$
by the expression\begin{equation}
\Omega_{\left(i,\eta\right)}\left(\left(X_{1},\vartheta_{1}\right),\left(X_{2},\vartheta_{2}\right)\right):=\int_{L_{0}}\left[\omega\left(X_{1},X_{2}\right)\eta+h_{1}\vartheta_{2}-h_{2}\vartheta_{1}\right]\label{eq:PWOmega}\end{equation}
 where $\alpha_{X_{k}}=dh_{k}$ for $h_{k}\in C^{\infty}\left(L_{0},\mathbb{R}\right)$.
This then descends to a weak symplectic structure (also labeled $\Omega$)
on $\text{Lag}_{pw}\left(M\right)^{\mathcal{F}_{pw}}$. So in particular
the isodrasts in $\text{Lag}_{pw}\left(M\right)$ are weakly symplectic
manifolds.

\section{Embeddings into a symplectic manifold}

In this section we will make precise a heuristic construction by Donaldson
\cite{Donaldson1} of a symplectic structure and moment map for a
dif\/feomorphism group action restricted to the space of embeddings.

Let $S_{0}$ be a f\/ixed f\/inite dimensional, compact, and oriented
manifold equipped with a volume form $\eta_{0}$, and let $\left(M,\omega\right)$
be a f\/inite dimensional symplectic manifold. The set of embeddings
$\text{Emb}\left(S_{0},M\right)$ of $S_{0}$ into $M$ is an open
subset of the space $C^{\infty}\left(S_{0},M\right)$ of all smooth
maps. Thus $\text{Emb}\left(S_{0},M\right)$ is a smooth manifold
modeled on spaces $\Gamma\left(i^{*}TM\right)$ for $i\in\text{Emb}\left(S_{0},M\right)$.
Assign to each point $i\in\text{Emb}\left(S_{0},M\right)$ a skew
symmetric bilinear form on $T_{i}\text{Emb}\left(S_{0},M\right)$
via the expression\begin{equation}
\Omega_{i}^{D}\left(X_{1},X_{2}\right):=\int_{S_{0}}\omega\left(X_{1},X_{2}\right)\eta_{0}\label{eq:Omega^D}\end{equation}
 for $X_{1},X_{2}\in T_{i}\text{Emb}\left(S_{0},M\right)$.

\begin{proposition}
The pointwise assignment $\Omega^{D}$ in \eqref{eq:Omega^D}
defines a weak symplectic structure on $\emph{Emb}\left(S_{0},M\right)$.
\end{proposition}

\begin{proof} Checking smoothness amounts to checking the following
statement: If $X_{1}\left(t\right)$ and $X_{2}\left(t\right)$ are
smooth time dependent vector f\/ields on $M$, $i_{t}$ is a smooth
curve in $\text{Emb}\left(S_{0},M\right)$, and $s:\mathbb{R}\to\mathbb{R}$
is a smooth function, then the map
\begin{gather*}
\mathbb{R}  \to  \mathbb{R},\\
t  \mapsto  \int_{S_{0}}\omega_{i_{s\left(t\right)}}\left(X_{1}\left(t\right)\circ i_{s\left(t\right)},X_{2}\left(t\right)\circ i_{s\left(t\right)}\right)\eta_{0}
\end{gather*}
 is smooth. This statement follows from the smoothness of all functions
in the integrand.

We will now prove closedness by choosing special extensions of tangent
vectors to vector f\/ields on $\text{Emb}\left(S_{0},M\right)$. Let
$X_{1},X_{2},X_{3}\in T_{i}\text{Emb}\left(S_{0},M\right)$ be tangent
vectors. Let $Z_{1}$, $Z_{2}$ and~$Z_{3}$ be vector f\/ields def\/ined
on a neighbourhood of $i\left(S_{0}\right)$ in $M$ such that $Z_{k}\circ i=X_{k}$.
Let $\xi_{1}$, $\xi_{2}$ and $\xi_{3}$ be vector f\/ields def\/ined
on the chart neighbourhood $U_{i}\subset\text{Emb}\left(S_{0},M\right)$
by the expression $\xi_{k}\left(j\right):=Z_{k}\circ j$. For these
particular vector f\/ields, Lie brackets like $\left[\xi_{1},\xi_{2}\right]$
at a point $i\in\text{Emb}\left(S_{0},M\right)$ can be written in
terms of the Lie bracket $\left[Z_{1},Z_{2}\right]$:
\begin{gather*}
\left[\xi_{1},\xi_{2}\right]\left(i\right)  =  \text{d}\xi_{2}\left(i\right)\cdot\xi_{1}\left(i\right)-\text{d}\xi_{1}\left(i\right)\cdot\xi_{2}\left(i\right)
  =  \text{d}Z_{2}\cdot Z_{1}\circ i-\text{d}Z_{1}\cdot Z_{2}\circ i
  =  \left[Z_{1},Z_{2}\right]\circ i.
 \end{gather*}
 Thus,
\begin{gather*}
d\Omega_{i}^{D}\left(X_{1},X_{2},X_{3}\right)   =  \left(d\Omega^{D}\left(\xi_{1},\xi_{2},\xi_{3}\right)\right)_{i}
   =   \left(\xi_{1}\Omega^{D}\left(\xi_{2},\xi_{3}\right)+\xi_{2}\Omega^{D}\left(\xi_{3},\xi_{1}\right)
   +\xi_{3}\Omega^{D}\left(\xi_{1},\xi_{2}\right)\right.\\
\left.\phantom{d\Omega_{i}^{D}\left(X_{1},X_{2},X_{3}\right)   =}{} -\Omega^{D}\left(\left[\xi_{1},\xi_{2}\right],\xi_{3}\right)-\Omega^{D}\left(\left[\xi_{2},\xi_{3}\right],\xi_{1}\right)-\Omega^{D}\left(\left[\xi_{3},\xi_{2}\right],\xi_{1}\right)\right)_{i}\\
\phantom{d\Omega_{i}^{D}\left(X_{1},X_{2},X_{3}\right)  }{}   =   \int_{S_{0}}d\omega\left(Z_{1},Z_{2},Z_{3}\right)\eta_{0}  =  0
\end{gather*}
 by closedness of $\omega$.

As for weak nondegeneracy, suppose that $X_{1}\in\Gamma\left(i^{*}TM\right)$
is nonzero on a neighbourhood~$W$ of $x\in S_{0}$. Let $j$ be a
compatible almost complex structure on $M$ (i.e.\ $g\left(\cdot,\cdot\right):=\omega\left(\cdot,j\cdot\right)$
is a~Riemannian metric on $M$). Let $\chi$ be a positive function
supported on $W$. Def\/ine $X_{2}\in\Gamma\left(i^{*}TM\right)$ by
$X_{2}\left(x\right):=\chi\left(i\left(x\right)\right)\cdot jX_{1}\left(x\right)$.
It follows that
\begin{gather*}
\int_{S_{0}}\omega\left(X_{1},X_{2}\right)\eta_{0}=\int_{S_{0}}\left(\chi\circ i\right)\cdot g\left(X_{1},X_{1}\right)\eta_{0}>0.\tag*{\qed}
\end{gather*}
\renewcommand{\qed}{}
\end{proof}

\begin{remark}
Donaldson originally def\/ined the above weakly symplectic structure
$\Omega^{D}$ on the space of smooth mappings $C^{\infty}(S_{0},M)$.
When $S_{0}=M$, the 2-form $\Omega^{D}$ restricts to a symplectic
structure def\/ined by Khesin and Lee on the open subset of orientation
preserving dif\/feo\-mor\-phisms of $M$ (relative to the orientation induced
by the symplectic volume form, which is taken to be~$\eta_{0}$; see
Section~3 in \cite{KhesinLee}).
\end{remark}

The Lie group $\text{Dif\/f}\left(S_{0},\eta_{0}\right)$ acts freely
on $\text{Emb}\left(S_{0},M\right)$ via\[
a\cdot i:=i\circ a^{-1}.\]

\begin{definition}
Let $G$ be a Lie group with Lie algebra $\mathfrak{g}$. By $\mathfrak{g}^{*}$
we will mean all bounded linear functionals on the convenient vector
space $\mathfrak{g}$. Let $\left\langle \cdot,\cdot\right\rangle :\mathfrak{g}^{*}\times\mathfrak{g}\to\mathbb{R}$
denote the canonical pairing between $\mathfrak{g}^{*}$ and $\mathfrak{g}$.
The \emph{coadjoint representation of} $G$, $\text{Ad}^{*}:G\to GL(\mathfrak{g}^{*})\subset L(\mathfrak{g}^{*},\mathfrak{g}^{*})$,
is def\/ined by\[
\left\langle \text{Ad}_{g}^{*}\zeta,\xi\right\rangle :=\left\langle \zeta,\text{Ad}_{g^{-1}}\xi\right\rangle \qquad\text{for any} \ \xi\in\mathfrak{g}.\]
\end{definition}

\begin{definition}
Let $\left(M,\sigma\right)$ be a weakly symplectic smooth manifold.
Let $G\times M\to M$ be a~smooth action of a Lie group $G$ on $M$,
such that $l_{g}^{*}\sigma=\sigma$ for all $g\in G$. This $G$ action
is called \emph{Hamiltonian} if there exists a $G$ equivariant map
(called the \emph{moment map}) $\mu\in C^{\infty}\left(M,\mathfrak{g}^{*}\right)$
 such that for all $\xi\in\mathfrak{g}$, the function $\left\langle \mu,\xi\right\rangle \in C^{\infty}\left(M,\mathbb{R}\right)$
is a Hamiltonian for $\xi_{M}$:\[
\text{d}\left\langle \mu,\xi\right\rangle =\imath\left(\xi_{M}\right)\sigma.\]
\end{definition}

\begin{proposition}\label{sub:PropDMomentMap}
If $i^{*}\left[\omega\right]$ is the zero class in
$H^{2}\left(S_{0}\right)$ for all $i\in\emph{Emb}\left(S_{0},M\right)$
and $H^{1}\left(S_{0}\right)=0$, then the $\emph{Dif\/f}\left(S_{0},\eta_{0}\right)$
action on $\emph{Emb}\left(S_{0},M\right)$ is Hamiltonian.
\end{proposition}

\begin{proof}
We f\/irst note that the $\text{Dif\/f}\left(S_{0},\eta_{0}\right)$
action on $\text{Emb}\left(S_{0},M\right)$ is symplectic:\begin{gather*}
\left(l_{a}^{*}\Omega^{D}\right)_{i}\left(X_{1},X_{2}\right)   =   \int_{S_{0}}\omega_{i\circ a^{-1}\left(x\right)}\left(X_{1}\circ a^{-1}\left(x\right),X_{2}\circ a^{-1}\left(x\right)\right)\eta_{0}\left(x\right)\\
 \phantom{\left(l_{a}^{*}\Omega^{D}\right)_{i}\left(X_{1},X_{2}\right)}{}  =   \int_{S_{0}}\omega_{i\left(y\right)}\left(X_{1}\left(y\right),X_{2}\left(y\right)\right)
 \left(a^{-1}\right)^{*}\eta_{0}\left(y\right)
   =   \Omega_{i}^{D}\left(X_{1},X_{2}\right).
   \end{gather*}
 If $i^{*}\omega=dA,$ def\/ine $\mu:\text{Emb}\left(S_{0},M\right)\to\mathfrak{X}\left(S_{0},\eta_{0}\right)^{*}$
by\[
\left\langle \mu,\xi\right\rangle :=\int_{S_{0}}A\left(\xi\right)\eta_{0}.\]
 This def\/inition is independent of the choice of $A$ since $H^{1}\left(S_{0}\right)=0$.
The map $\mu$ is smooth by the usual local arguments. To check that
$\mu$ is a moment map, let $X\in\Gamma\left(i^{*}TM\right)$ be a
tangent vector at $i\in\text{Emb}\left(S_{0},M\right)$. Let $Z$
be a vector f\/ield on a neighbourhood of $i\left(S_{0}\right)$ satisfying
$Z\circ i=X$, and suppose $Z$ generates a f\/low $\tau_{t}$ on $M$.
Let $A_{t}$ be a smooth curve in $\Omega^{1}\left(S_{0}\right)$
satisfying $dA_{t}:=\left(\tau_{t}\circ i\right)^{*}\omega$. Then\[
\left.\frac{d}{dt}\right|_{t=0}dA_{t}=i^{*}\left.\frac{d}{dt}\right|_{t=0}\tau_{t}^{*}\omega=i^{*}\mathcal{L}_{Z}\omega=i^{*}\left(d\circ\imath\left(Z\right)\omega\right).\]
 This together with $H^{1}\left(S_{0}\right)=0$ imply that $\left.\frac{d}{dt}\right|_{t=0}A_{t}=i^{*}\imath\left(Z\right)\omega+dh$
for some $h\in C^{\infty}\left(S_{0},\mathbb{R}\right)$. It follows
that\[
\left\langle d\mu\left(X\right),\xi\right\rangle =\int_{S_{0}}\left[\omega\left(X,i_{*}\xi\right)+dh\left(\xi\right)\right]\eta_{0}=\Omega^{D}\left(X,i_{*}\xi\right)\]
 which verif\/ies the moment map condition.

Finally, $\mu$ is $\text{Dif\/f}\left(S_{0},\eta_{0}\right)$ equivariant:
\begin{gather*}
\left\langle \mu\left(l_{a}\left(i\right)\right),\xi\right\rangle    =   \int_{S_{0}}\left(a^{-1}\right)^{*}A\left(\xi\right)\eta_{0}
  =   \int_{S_{0}}A\left(\text{d}a^{-1}\cdot\xi\right)\eta_{0} \\
\phantom{\left\langle \mu\left(l_{a}\left(i\right)\right),\xi\right\rangle}{}
  =   \left\langle \mu\left(i\right),\text{Ad}_{a^{-1}}\xi\right\rangle
  =   \left\langle \text{Ad}_{a}^{*}\mu\left(i\right),\xi\right\rangle .\tag*{\qed}
  \end{gather*}\renewcommand{\qed}{}
\end{proof}

Let us now consider the special case when $\int_{S_{0}}\eta_{0}=1$,
the manifold $S_{0}$ is half the dimension of $M$, and assume that
the topological conditions $H^{1}(S_{0})=0$ and $[i^{*}\omega]=0\in H^{2}(S_{0})$
in Proposition~\ref{sub:PropDMomentMap} hold so that we have a moment map
$\mu$ on $\text{Emb}\left(S_{0},M\right)$. The level surface $\mu^{-1}\left\{ 0\right\} $
is given by\[
\mu^{-1}\left\{ 0\right\} =\left\{ i\in\text{Emb}\left(S_{0},M\right)\mid i^{*}\omega=0\right\} =\text{Lag}\left(S_{0},M\right).\]
 The group $\text{Dif\/f}\left(S_{0},\eta_{0}\right)$ acts freely on
$\mu^{-1}\left\{ 0\right\} $ in the usual way, with the quotient
given by\[
\mu^{-1}\left\{ 0\right\} /\text{Dif\/f}\left(S_{0},\eta_{0}\right)=\text{Lag}\left(S_{0},M\right)/\text{Dif\/f}\left(S_{0},\eta_{0}\right).\]
 By Proposition~\ref{PWQuotient}, the set $\text{Lag}\left(S_{0},M\right)/\text{Dif\/f}\left(S_{0},\eta_{0}\right)$
is a smooth manifold modeled on convenient spaces $\Gamma_{\text{closed}}\left(i^{*}TM\right)/\mathfrak{X}\left(S_{0},\eta_{0}\right)$
for $i\in\text{Lag}\left(S_{0},M\right)$. In fact, the manifold $\text{Lag}\left(S_{0},M\right)$
is the total space of a principal $\text{Dif\/f}(S_{0},\eta_{0})$ bundle
over the quotient $\text{Lag}\left(S_{0},M\right)/\text{Dif\/f}\left(S_{0},\eta_{0}\right)$.
Since the 2-form $\left.\Omega^{D}\right|_{\mu^{-1}\left\{ 0\right\} }$
is basic, it descends to a unique 2-form $\Omega_{\text{red}}^{D}$
on $\mu^{-1}\left\{ 0\right\} /\text{Dif\/f}\left(S_{0},\eta_{0}\right)$.

Under the topological assumption that $H^{1}\left(S_{0}\right)=0$,
the subspace given by the isodrastic foliation $\mathcal{F}_{pw}$
at each point equals the entire tangent space to $\text{Lag}_{pw}\left(M\right)$
at that point. The weak symplectic structure $\Omega$ on isodrasts
in this case becomes well def\/ined on all of $\text{Lag}_{pw}\left(M\right)$.
In fact, the pull back of $\Omega$ under the dif\/feomorphism $\upsilon$
in Proposition~\ref{PWQuotient} is given by\[
\left(\upsilon^{*}\Omega\right)\left(\left[X_{1}\right],\left[X_{2}\right]\right)=\int_{L_{0}}\omega\left(X_{1},X_{2}\right)\eta_{0}=\Omega_{\text{red}}^{D}\left(\left[X_{1}\right],\left[X_{2}\right]\right).\]
 In other words, the ``symplectic quotient'' $\left(\mu^{-1}\left\{ 0\right\} /\text{Dif\/f}_{+}\left(S_{0},\eta_{0}\right),\Omega_{\text{red}}^{D}\right)$
is symplectomorphic to $\left(\text{Lag}_{pw}\left(M\right),\Omega\right)$.

This last result is suggestive, leading one to wonder if the symplectic
structure $\Omega^{D}$ on $\text{Emb}\left(S_{0},M\right)$ might
be related to the symplectic structure $\Omega$ on isodrasts in $\text{Lag}_{pw}\left(M\right)$
via some sort of reduction procedure. In the next section we will
make this relationship clear.

\section[Reduction of $\left(\text{Emb}\left(S_{0},M\right),\Omega^{D}\right)$]{Reduction of $\boldsymbol{\left(\text{Emb}\left(S_{0},M\right),\Omega^{D}\right)}$}

In this section we will def\/ine a notion of reduction of weakly symplectic
spaces with respect to a~group action. We will then show that the
image of $\text{Ham}\left(M\right)$ orbits through isotropic embeddings
in $\text{Emb}\left(S_{0},M\right)$ under the projection to the quotient
$\text{Emb}\left(S_{0},M\right)/\text{Dif\/f}\left(S_{0},\eta_{0}\right)$
are reductions. Moreover, these reductions are naturally symplectomorphic
to the spaces of positive weighted isotropic submanifolds of $M$.
In particular, when $\dim S_{0}=\frac{1}{2}\dim M$, the
particular reductions of $\text{Emb}\left(S_{0},M\right)$ consisting
of positive weighted Lagrangian submanifolds are symplectomorphic
to isodrasts in $\text{Lag}_{pw}\left(M\right)$.

In Proposition~\ref{sub:PropDMomentMap}, the topological assumption $H^{1}\left(S_{0}\right)=0$
was essential in def\/ining a moment map. Since the transverse spaces
to the leaves of an isodrastic foliation become nontrivial exactly
when $H^{1}\left(L_{0}\right)$ is nontrivial, we would like to remove
such a topological condition on $S_{0}$. This means we must use a
notion of reduction that does not depend on having a moment map.

Let us begin by looking at the standard reduction of a f\/inite dimensional
symplectic manifold $\left(P,\sigma\right)$ with respect to a Hamiltonian
$G$ action using a moment map $\mu$. Suppose $r$ is a regular value
of $\mu$. The tangent space at $p$ to the level surface $\mu^{-1}\left\{ r\right\} $
is equal to the set $\mathcal{D}_{p}$ of all vectors $X\in T_{p}P$
satisfying $\sigma\left(X,\xi_{P}\left(p\right)\right)=0$ for all
$\xi\in\mathfrak{g}$. These subspaces $\mathcal{D}_{p}$ are def\/ined
for any symplectic action on $P$, even in the absence of a moment
map, and def\/ine a distribution $\mathcal{D}$ on $P$. If $G\circlearrowright P$
is a free symplectic action, then this distribution can be taken as
the starting point of the ``optimal reduction method'' of Juan-Pablo
Ortega and Tudor S.~Ratiu \cite{OrtegaRatiu}. We will not describe
the details here, but simply note that for a symplectic $G$ action
$G\circlearrowright(P,\sigma)$
\begin{itemize}\itemsep=0pt
\item the optimal reduction method yields symplectic spaces $\left(P_{\rho},\sigma_{\rho}\right)$
where $\rho\in P/G_{\mathcal{D}}$ and $G_{\mathcal{D}}$ is the pseudogroup
of f\/lows generated by Hamiltonian vector f\/ields in $\mathfrak{X}_{\mathcal{D}}\left(P\right)$
corresponding to $G$-invariant Hamiltonian functions;
\item the ``optimal momentum map'' is given by the projection $\mathcal{J}:P\to P/G_{\mathcal{D}}$;
\item each reduced symplectic space $\left(P_{\rho},\sigma_{\rho}\right)$
is the quotient of an integral manifold of $\mathcal{D}$ (the level
surface $\mathcal{J}^{-1}\left\{ \rho\right\} $) by the stabilizer
$G_{\rho}$ at $\rho$ under the action $G\circlearrowright P/G_{\mathcal{D}}$:
$g\cdot\left[p\right]:=\left[g\cdot p\right]$;
\item if $G\circlearrowright P$ is a free Hamiltonian action with moment
map $\mu$, and the point $r\in\mathfrak{g}^{*}$ is a~regular value
of $\mu$, and $\mu^{-1}\left\{ r\right\} $ is connected, then $\mu^{-1}\left\{ r\right\} $
is a $G_{\mathcal{D}}$ orbit $\rho$ and the reduced symplectic space
$P_{\rho}$ coincides with the symplectic quotient $\mu^{-1}\left\{ r\right\} /G_{r}$
(here $G_{r}$ denotes the stabilizer of $r\in\mathfrak{g}^{*}$ with
respect to the coadjoint action of $G$).
\end{itemize}
This suggests a way to def\/ine reduction in the inf\/inite dimensional
case, and motivates the following def\/inition:

\begin{definition}
Let $\left(P,\sigma\right)$ be a weakly symplectic smooth manifold.
Let $G\circlearrowright P$ be a smooth free action of a Lie group
$G$ on $P$, such that $l_{g}^{*}\sigma=\sigma$ for all $g\in G$.
The collection of subspaces\[
\mathcal{D}_{x}:=\left\{ X\in T_{x}M\mid\sigma_{x}\left(X,\xi_{M}\left(x\right)\right)=0\ \forall\, \xi\in\mathfrak{g}\right\} \]
 for $x\in P$ def\/ines a distribution $\mathcal{D}$ on $P$. Let
$i_{N}:N\hookrightarrow P$ be a maximal integral manifold of $\mathcal{D}$
and let $q:P\to P/G$ denote the projection to the orbit space. Suppose
that
\begin{enumerate}\itemsep=0pt
\item[1)] $q\left(N\right)$ is a smooth manifold,
\item[2)] there exists a unique weak symplectic structure $\sigma_{\text{red}}$
on $q\left(N\right)$ such that $\left(\left.q\right|_{N}\right)^{*}\sigma_{\text{red}}=i_{N}^{*}\sigma$.
\end{enumerate}
Then the weakly symplectic manifold $\left(q\left(N\right),\sigma_{\text{red}}\right)$
will be called a \emph{reduction} or \emph{symplectic quotient} of
$\left(P,\sigma\right)$ with respect to the $G$ action.
\end{definition}

\subsection*{Convention}

From now on, we will make no topological assumptions on $i^{*}\omega$
or $H^{1}\left(S_{0}\right)$.

For the $\text{Dif\/f}\left(S_{0},\eta_{0}\right)$ action on the symplectic
manifold $\left(\text{Emb}\left(S_{0},M\right),\Omega^{D}\right)$,
the subspaces $\mathcal{D}_{i}$ can be described in very familiar
terms:

\begin{proposition}
For every $i\in\emph{Emb}\left(S_{0},M\right)$,\begin{equation}
\mathcal{D}_{i}=\Gamma_{\emph{exact}}\left(i^{*}TM\right)=\left\{ X\in\Gamma\left(i^{*}TM\right)\mid\alpha_{X}\in B^{1}\left(S_{0}\right)\right\} .\label{eq:DistributionD}\end{equation}
\end{proposition}

\begin{proof}
The distribution $\mathcal{D}$ on $\text{Emb}\left(S_{0},M\right)$
is def\/ined by\[
\mathcal{D}_{i}:=\left\{ X\in\Gamma\left(i^{*}TM\right)\mid\ \Omega_{i}^{D}\left(X,\xi_{\text{Emb}\left(S_{0},M\right)}\left(i\right)\right)=0,\ \forall\, \xi\in\mathfrak{X}\left(S_{0},\eta_{0}\right)\right\} \]
 for $i\in\text{Emb}(S_{0},M)$. Since $\int_{S_{0}}dh\left(\xi\right)\eta_{0}=0$
for any function $h$ on $S_{0}$ and all $\xi\in\mathfrak{X}\left(S_{0},\eta_{0}\right)$,
it follows that $\left\{ X\in\Gamma\left(i^{*}TM\right)\mid\alpha_{X}\in B^{1}\left(S_{0}\right)\right\} \subset\mathcal{D}_{i}$.

Let $X\in\mathcal{D}_{i}$, i.e. $\int_{S_{0}}\alpha_{X}\left(\xi\right)\eta_{0}=\int_{S_{0}}\alpha_{X}\wedge\imath_{\xi}\eta_{0}=0$
for all divergence free $\xi$. If $U$ is a coordinate neighbourhood
in $S_{0}$, $\eta_{0}=\eta dx_{1}\wedge\cdots\wedge dx_{n}$ on $U$,
and $f$ a function with $\text{supp}\left(f\right)\subset U$, then
def\/ine the divergence free vector f\/ield $Y_{12}:=\frac{1}{\eta}\left(\frac{\partial f}{\partial x_{2}}\frac{\partial}{\partial x_{1}}-\frac{\partial f}{\partial x_{1}}\frac{\partial}{\partial x_{2}}\right)$.
If $\alpha_{X}=a_{i}dx_{i}$ on $U$, then\begin{gather*}
\int_{S_{0}}\alpha_{X}\wedge\imath_{Y_{12}}\eta_{0}   =   \int_{U}\left(a_{1}\frac{\partial f}{\partial x_{2}}-a_{2}\frac{\partial f}{\partial x_{1}}\right)dx_{1}\wedge\cdots\wedge dx_{n}\\
\phantom{\int_{S_{0}}\alpha_{X}\wedge\imath_{Y_{12}}\eta_{0}}{}
   =   \int_{U}\left(\frac{\partial a_{2}}{\partial x_{1}}-\frac{\partial a_{1}}{\partial x_{2}}\right)fdx_{1}\wedge\cdots\wedge dx_{n}
   \end{gather*}
 where we have used integration by parts. If this is to vanish for
all $f$ then $\frac{\partial a_{2}}{\partial x_{1}}=\frac{\partial a_{1}}{\partial x_{2}}$.
Similarly, by considering vector f\/ields like $Y_{13}:=\frac{1}{\eta}\left(\frac{\partial f}{\partial x_{3}}\frac{\partial}{\partial x_{1}}-\frac{\partial f}{\partial x_{1}}\frac{\partial}{\partial x_{3}}\right)$
it follows that $\frac{\partial a_{3}}{\partial x_{1}}=\frac{\partial a_{1}}{\partial x_{3}}$,
etc., which means $\alpha_{X}$ is closed.

Let $g_{0}$ be a Riemannian metric on $S_{0}$ whose volume form
equals $\eta_{0}$. For every $X\in\mathcal{D}_{i}$, since $\alpha_{X}$
is closed, there exists a function $h$ on $S_{0}$ such that $\beta_{X}:=\alpha_{X}-dh$
is harmonic. Moreover, $\int_{S_{0}}\beta_{X}\left(\xi\right)\eta_{0}=0$
for every $\xi\in\mathfrak{X}\left(S_{0},\eta_{0}\right)$. Def\/ine
the vector f\/ield $Y_{\beta_{X}}$ on $S_{0}$ by $\beta_{X}=g_{0}\left(Y_{\beta_{X}},\cdot\right)$.
Let $V$ be a coordinate neighbourhood in $S_{0}$, and suppose $\beta_{X}=b_{i}dx_{i}$
on $V$. On $V$,\begin{gather*}
\mathcal{L}_{Y_{\beta_{X}}}\eta_{0}  =   d\circ\imath_{Y_{\beta_{X}}}\left(\sqrt{\det g_{0}}dx_{1}\wedge\cdots\wedge dx_{n}\right)\\
   \phantom{\mathcal{L}_{Y_{\beta_{X}}}\eta_{0}}{}
   =   d\left(\sqrt{\det g_{0}}g_{0}^{ij}b_{j}\left(-1\right)^{i}dx_{1}\wedge\cdots\wedge\hat{dx}_{i}\wedge\cdots\wedge dx_{n}\right)
   =   d*\beta_{X}=*\delta\beta_{X}=0
   \end{gather*}
 which means $Y_{\beta_{X}}\in\mathfrak{X}\left(S_{0},\eta_{0}\right)$.
So $0=\int_{S_{0}}\beta_{X}\left(Y_{\beta_{X}}\right)\eta_{0}=\int_{S_{0}}g_{0}\left(Y_{\beta_{X}},Y_{\beta_{X}}\right)\eta_{0}$
implies that $\beta_{X}=0$, which means $\alpha_{X}=dh$.
\end{proof}

The group $\text{Ham}\left(M\right)$ acts freely on $\text{Emb}\left(S_{0},M\right)$
under the action\[
\text{Ham}\left(M\right)\circlearrowright\text{Emb}\left(S_{0},M\right):\quad\phi\cdot i:=\phi\circ i.\]
 In what follows we will show that $\text{Ham}\left(M\right)$ orbits
through isotropic embeddings are maximal integral manifolds of $\mathcal{D}$.
For this purpose, we will need to make use of the following isotropic
embedding theorem of Weinstein:

\begin{theorem}[see \cite{Weinstein2}]\label{sub:ThWIsotropicET} Let $\left(M,\omega\right)$ be
a symplectic manifold and $i:I\hookrightarrow M$ be an isotropic
submanifold, i.e.~$i^{*}\omega=0$. The vector bundle $T^{*}I\oplus\left(TI^{\omega}/TI\right)$
admits a symplectic structure on a~neighbourhood of the zero
section, which is given by $\omega_{T^{*}I}+\omega_{\mathbb{R}^{2n}}$
on the zero section. Furthermore, there exists a neighbourhood
$U_{1}$ of $I$ in $M$, a neighbourhood $U_{2}$
of $I$ in $T^{*}I\oplus\left(TI^{\omega}/TI\right)$,
and a~symplectomorphism from $U_{1}$ to $U_{2}$ fixing
$I$.
\end{theorem}

\begin{proposition}\label{sub:PropReductions}
If $i\in\emph{Emb}\left(S_{0},M\right)$ is an isotropic
embedding, then the $\emph{Ham}\left(M\right)$ orbit~$\mathcal{N}$
through $i$ is a maximal integral manifold of the distribution
$\mathcal{D}$ $($see \eqref{eq:DistributionD}$)$. The restriction
of~$\Omega^{D}$ to~$\mathcal{N}$ in~$\emph{Emb}\left(S_{0},M\right)$
descends to a unique weak symplectic structure $\Omega_{\emph{red}}^{D}$
on the image  \mbox{$\mathcal{O}{:=}q\left(\mathcal{N}\right)$} in the
orbit space under the projection $q:\emph{Emb}\left(S_{0},M\right)\to\emph{Emb}\left(S_{0},M\right)/\emph{Dif\/f}\left(S_{0},\eta_{0}\right)$.
Thus the symplectic spaces $\left(\mathcal{O},\Omega_{\emph{red}}^{D}\right)$
are reductions of $\emph{Emb}\left(S_{0},M\right)$ with
respect to the $\emph{Dif\/f}\left(S_{0},\eta_{0}\right)$~action.
\end{proposition}

\begin{proof}
We will f\/irst show that a $\text{Ham}\left(M\right)$
orbit through an isotropic embedding is an integral manifold of $\mathcal{D}$.
Suppose that $i:S_{0}\to M$ is an isotropic embedding. By Theorem~\ref{sub:ThWIsotropicET},
if $S=i\left(S_{0}\right)$ and $N=TS^{\omega}/TS$ denotes the symplectic
normal bundle, then by choosing a symplectic embedding $\lambda:M\supset U\to T^{*}S\oplus N$
we can assume that $M=T^{*}S\oplus N$ and that $i$ is the zero section
inclusion. Let $\text{pr}_{1}:T^{*}S\oplus N\to T^{*}S$ and $\text{pr}_{2}:T^{*}S\oplus N\to N$
be the natural maps which forget the points in the symplectic normal
and cotangent f\/ibers respectively. That is, $\text{pr}_{1}(x,p,v):=(x,p)$
and $\text{pr}_{2}(x,p,v):=(x,v)$. Given $X\in\mathcal{D}_{i}$ with
$\alpha_{X}=dh$ for some $h\in C^{\infty}\left(S_{0},\mathbb{R}\right)$,
let $X=X_{\text{f\/ib}}+X_{\text{tan}}$ denote the decomposition of
$X$ into components tangent to the f\/ibers and tangent to $S$. Extend
$X_{\text{f\/ib}}$ constantly along the f\/ibers in $T^{*}S\oplus N$
to a vector f\/ield $Z$ def\/ined on a neighbourhood of the zero section.
It follows that $Z$ is a Hamiltonian vector f\/ield satisfying $\imath\left(Z\right)\omega=d(\text{pr}_{1}^{*}\pi_{T^{*}S}^{*}f+\text{pr}_{2}^{*}H_{N})$
where $f\in C^{\infty}\left(S,\mathbb{R}\right)$ satisf\/ies $i^{*}f=h$,
and $H_{N}\in C^{\infty}(N,\mathbb{R})$ is def\/ined by $H_{N}(x,v):=\omega_{\mathbb{R}^{2n}}(\text{d}\text{pr}_{2}\circ X(x),v)$.
Since $\omega=\omega_{T^{*}S}+\omega_{\mathbb{R}^{2n}}$ along the
zero section, $\left.\imath(X_{\text{tan}})\omega\right|_{TS}=0$.
Thus for each point $x$ in the zero section, we have that $\imath(X_{\text{tan}})\omega$
def\/ines a~linear functional $(\imath(X_{\text{tan}})\omega)_{x}$
on the f\/iber $T_{x}^{*}S\oplus N_{x}$. The smooth function $H_{\text{tan}}:T^{*}S\oplus N\to\mathbb{R}$
def\/ined by $H_{\text{tan}}(x,p,v):=(\imath(X_{\text{tan}})\omega)_{x}(p,v)$
is the primitive of a Hamiltonian vector f\/ield $Z_{\text{tan}}$ satisfying
$Z_{\text{tan}}\circ i=X_{\text{tan}}$. It follows that $\mathcal{D}_{i}\subset T_{i}(\text{Ham}(M)\cdot i)$.
The converse inclusion follows from the fact that $\alpha_{X_{H}\circ i}=di^{*}H$
for any Hamiltonian vector f\/ield $X_{H}$. So $\text{Ham}\left(M\right)$
orbits must be integral manifolds of $\mathcal{D}$. The fact that
they must be maximal can be shown as in the proof of Proposition~\ref{sub:PropLag(L_0,M)Foliation}.

Let $\mathcal{N}$ be a $\text{Ham}\left(M\right)$ orbit through
an isotropic embedding $i$. The image $\mathcal{O}:=q\left(\mathcal{N}\right)$
is equal to the $\text{Ham}\left(M\right)$ orbit in the quotient
$\text{Emb}\left(S_{0},M\right)/\text{Dif\/f}\left(S_{0},\eta_{0}\right)$
under the action $\phi\cdot\left[i\right]:=\left[\phi\circ i\right]$,
it consists of positive weighted isotropic submanifolds of $M$ dif\/feomorphic
to $S_{0}$, and can be given a smooth manifold structure similar
to that on the space of positive weighted Lagrangian submanifolds.
The space $\text{Dif\/f}\left(S_{0},\eta_{0}\right)\cdot\mathcal{N}$
is the total space of a principal $\text{Dif\/f}\left(S_{0},\eta_{0}\right)$
bundle over $\mathcal{O}$. The pullback of $\Omega^{D}$ to $\text{Dif\/f}\left(S_{0},\eta_{0}\right)\cdot\mathcal{N}$
def\/ines a closed basic 2-form which descends to a closed 2-form $\Omega_{\text{red}}^{D}$
on $\mathcal{O}$. To check weak nondegeneracy, we f\/irst note that
the 2-form $\Omega_{\text{red}}^{D}$ is given by the expression\[
\Omega_{\text{red}}^{D}\left(\left[X_{1}\right],\left[X_{2}\right]\right)=\int_{S_{0}}\omega\left(X_{1},X_{2}\right)\eta_{0}.\]
 Again, by choosing a symplectic embedding $\lambda:M\supset U\to T^{*}S\oplus N$
we can assume that $M=T^{*}S\oplus N$ and that $i$ is the zero section
inclusion. Given $\left[X_{2}\right]\in\text{Ker}\left(\Omega_{\text{red}}^{D}\right)$,
let $X_{2}=X_{\text{f\/ib}}+X_{\text{tan}}$ be the decomposition of
$X_{2}$ into components tangent to the f\/ibers and to the zero section
respectively. Extend $X_{\text{f\/ib}}$ and $X_{\text{tan}}$ to Hamiltonian
vector f\/ields $Z_{\text{f\/ib}}$ and $Z_{\text{tan}}$ respectively
as before. For $\left[X_{1}\right]\in T_{\left[i\right]}\mathcal{O}$,
let $Z_{f_{1}}$be a Hamiltonian vector f\/ield on $T^{*}S\oplus N$
such that $X_{1}=Z_{f_{1}}\circ i$. It follows that\begin{gather*}
0=\Omega_{\text{red}}^{D}\left(\left[X_{1}\right],\left[X_{2}\right]\right)   =   \int_{S_{0}}\!\omega\left(Z_{f_{1}},X_{2}\right)\eta_{0}
   =   \int_{S_{0}}\!i^{*}\left[\omega\left(Z_{f_{1}},Z_{\text{f\/ib}}\right)\right]\eta_{0}
   +\!\int_{S_{0}}\!i^{*}\left[\omega\left(Z_{f_{1}},Z_{\text{tan}}\right)\right]\eta_{0}\!\\
\phantom{0=\Omega_{\text{red}}^{D}\left(\left[X_{1}\right],\left[X_{2}\right]\right)}{}  =  \!\int_{S_{0}}\!i^{*}\mathcal{L}_{Z_{\text{f\/ib}}}f_{1}\eta_{0}+\!\int_{S_{0}}\!i^{*}\mathcal{L}_{Z_{\text{tan}}}f_{1}\eta_{0}
  =  \int_{S_{0}}\!i^{*}\mathcal{L}_{Z_{\text{f\/ib}}}f_{1}\eta_{0}+\!\int_{S_{0}}\!\mathcal{L}_{Y}\left(i^{*}f_{1}\right)\eta_{0},\!
  \end{gather*}
 where $Y\in\mathfrak{X}\left(S_{0}\right)$ and $i_{*}Y=Z_{\text{tan}}\circ i$.
If $f_{1}$ is the pullback through the bundle projection of a function
on the base $S$, then the f\/irst term vanishes and the second term
vanishes if and only if $Y\in\mathfrak{X}\left(S_{0},\eta_{0}\right)$.
It follows that $Z_{\text{f\/ib}}$ must vanish which means $\left[X_{2}\right]=0$.

Finally, uniqueness of the 2-form $\Omega_{\text{red}}^{D}$ follows from the
fact that the principal $\text{Dif\/f}\left(S_{0},\eta_{0}\right)$
bundle $\text{Emb}\left(S_{0},\eta_{0}\right)\to\text{Emb}\left(S_{0},\eta_{0}\right)/\text{Dif\/f}\left(S_{0},\eta_{0}\right)$
restricted to $\mathcal{O}$ admits local sections.
\end{proof}

Thus
the image $\mathcal{O}:=q\left(\mathcal{N}\right)$ of a $\text{Ham}\left(M\right)$
orbit in $\text{Emb}\left(S_{0},M\right)$ through an isotropic embedding
is a symplectic quotient of $\left(\text{Emb}\left(S_{0},M\right),\Omega^{D}\right)$
with respect to the $\text{Dif\/f}\left(S_{0},M\right)$ action.

\subsection*{Convention}

From now on we will assume that $\int_{S_{0}}\eta_{0}=1$.

\begin{theorem}\label{sub:ThWisD}
Suppose $\emph{dim}\left(S_{0}\right)=\frac{1}{2}\emph{dim}\left(M\right)=n$
and that \emph{$\mathcal{N}$} is a $\emph{Ham}\left(M\right)$
orbit through a Lagrangian embedding, i.e.\ $\mathcal{N}$ is
an isodrast in $\emph{Lag}\left(S_{0},M\right)$. If $\emph{\textbf{I}}_{w}:=p\left(\mathcal{N}\times\emph{Vol}_{1}\left(S_{0}\right)\right)$,
where $p:\emph{Lag}_{pw}\left(S_{0},M\right)\to\emph{Lag}_{pw}\left(M\right)$,
and $\mathcal{O}:=q\left(\mathcal{N}\right)$, then the reduction
$\left(\mathcal{O},\Omega_{\emph{red}}^{D}\right)$ of $\left(\emph{Emb}\left(S_{0},M\right),\Omega^{D}\right)$
with respect to the $\emph{Dif\/f}\left(S_{0},\eta_{0}\right)$
action is symplectomorphic to the isodrast $\left(\emph{\textbf{I}}_{w},\Omega\right)$
in $\emph{Lag}_{pw}\left(M\right)$.
\end{theorem}

\begin{proof}
By Proposition~\ref{PWQuotient}, we have a $\text{Ham}\left(M\right)$ equivariant
dif\/feomorphism\[
\upsilon: \ \text{Lag}\left(S_{0},M\right)/\text{Dif\/f}\left(S_{0},\eta_{0}\right)\to\text{Lag}_{pw}\left(M\right)\]
 which induces a dif\/feomorphism (also labeled by $\upsilon$)\[
\upsilon: \ \mathcal{O}\to\text{\textbf{I}}_{w}.\]
 The pull back of $\Omega$ under this map is given by\begin{gather*}
\left(\upsilon^{*}\Omega\right)_{\left[i\right]}\left(\left[X_{1}\right],\left[X_{2}\right]\right)    =  \Omega_{\left[\left(i,\eta_{0}\right)\right]}\left(\left[\left(X_{1},0\right)\right],
\left[\left(X_{2},0\right)\right]\right) \\
\phantom{\left(\upsilon^{*}\Omega\right)_{\left[i\right]}\left(\left[X_{1}\right],\left[X_{2}\right]\right) }{}
  =   \int_{S_{0}}\omega\left(X_{1},X_{2}\right)\eta_{0}
  =  \left(\Omega_{\text{red}}^{D}\right)_{\left[i\right]}\left(\left[X_{1}\right],\left[X_{2}\right]\right).
  \end{gather*}
 So indeed, the symplectic quotient $\left(\mathcal{O},\Omega_{\text{red}}^{D}\right)$
is symplectomorphic to the isodrast $\left(\text{\textbf{I}}_{w},\Omega\right)$.
\end{proof}

\begin{remark}
This last theorem clarif\/ies the relationship between Donaldson's symplectic
structure $\Omega^{D}$ on $\text{Emb}\left(S_{0},M\right)$ and Weinstein's
symplectic structure $\Omega^{W}$ on isodrasts in $\text{Lag}_{pw}\left(M\right)$.
Namely, the isodrasts $\left(\text{\textbf{I}}_{w},\Omega^{W}\right)$
can be viewed as symplectic quotients of $\left(\text{Emb}\left(S_{0},M\right),\Omega^{D}\right)$.
And while Donaldson's and Weinstein's constructions were done heuristically,
our constructions of the symplectic spaces $\left(\text{Emb}\left(S_{0},M\right),\Omega^{D}\right)$
and $\left(\text{\textbf{I}}_{w},\Omega\right)$ as well as the description
of the relationship between them are rigourous.
\end{remark}

\begin{example}
{\sloppy Let $\left(M,\omega\right)=\left(\mathbb{R}^{2},dq\wedge dp\right)$
and $S_{0}=S^{1}$. Since $S^{1}$ is one dimensional $\text{Lag}\left(S^{1},\mathbb{R}^{2}\right)$ \linebreak $=\text{Emb}\left(S^{1},\mathbb{R}^{2}\right)$.
However, since $H^{1}\left(S^{1}\right)$ is nontrivial, we have nontrivial
isodrasts in $\text{Emb}\left(S^{1},\mathbb{R}^{2}\right)$ that can
be described as follows.

}

Let $\beta=pdq$ denote the canonical 1-form on the plane. Given a
map $i\in\text{Emb}\left(S^{1},\mathbb{R}^{2}\right)$, the action
integral $A\left(i\right)$ of $i$ is def\/ined as the integral of
$\beta$ around $i\left(S^{1}\right)$:\[
A\left(i\right):=\int_{S^{1}}i^{*}\beta.\]
 An isotopic deformation is Hamiltonian if and only if the action
integrals are constant along the deformation (see Proposition~2.1 in \cite{Weinstein3})%
\footnote{In fact, this was Weinstein's original motivation for the terminology
``isodrast''.}. The idea is as follows. Given two nearby loops in a symplectic manifold
$\left(M,\omega\right)$, we can def\/ine the dif\/ference in their action
integrals to be the integral of $-\omega$ over a cylindrical surface
joining the two loops. This is well def\/ined even when $\omega$ is
not exact. Lagrangian submanifolds near a given $L\in\text{Lag}\left(M\right)$
can be identif\/ied with graphs of 1-forms in $T^{*}L$ by Theorem~\ref{sub:TheoremWLN}.
Two such graphs can be joined by a Hamiltonian deformation if and
only if their corresponding 1-forms are cohomologous. If $\gamma'$
is a small deformation of a~loop $\gamma$ in the zero section corresponding
to a deformation of the zero section to a Lagrangian submanifold $L'=\text{Graph}\left(\alpha\right)$,
$C$ denotes a cylinder joining $\gamma$ and $\gamma'$, and if $\beta_{T^{*}L}$
denotes the canonical 1-form of the $T^{*}L$, then
\begin{gather*}
\int_{C}-\omega_{T^{*}L}   =   \int_{\gamma'}\beta_{T^{*}L}  =  \int_{\gamma}\alpha.
\end{gather*}
 So a small deformation of the zero section is the graph of an exact
1-form if and only if the dif\/ference in action integrals is 0 for
all loops $\gamma$ and $\gamma'$ in the zero section and the deformed
image respectively. It follows that two graphs of 1-forms can be joined
by a Hamiltonian deformation if and only if the dif\/ference in action
integrals remains constant for all loops in these Lagrangian submanifolds.
Returning to our example, this means that the isodrasts in $\text{Emb}\left(S^{1},\mathbb{R}^{2}\right)$
consist of circle embeddings that can be joined by an isotopy that
preserves action integrals, i.e. $\text{\textbf{H}}\subset\text{Emb}\left(S^{1},\mathbb{R}^{2}\right)$
is an isodrast if and only if it consists of isotopic circle embeddings
and the map $A:\text{Emb}\left(S^{1},\mathbb{R}^{2}\right)\to\mathbb{R}$
sending $i$ to its action integral $A\left(i\right)$ is constant
on $\text{\textbf{H}}$.

On the circle we can take $\eta_{0}=\frac{dt}{2\pi}$ so that $\text{Dif\/f}\left(S^{1},\frac{dt}{2\pi}\right)$
consists of rigid rotations of the circle. Let $\mathcal{O}:=q\left(\text{\textbf{H}}\right)$
be the image of an isodrast $\text{\textbf{H}}$ in the orbit space.
Each representative $X$ for a~tangent vector $\left[X\right]\in T_{\left[i\right]}\mathcal{O}$
can be decomposed as $X=Z+Y$ where $Z$ is a normal to $i\left(S^{1}\right)$
and~$Y$ is tangent to $i\left(S^{1}\right)$. It follows that if
$\alpha_{X_{k}}=dh_{k}$ then the reduced symplectic structure is
given by\begin{gather*}
\Omega_{\text{red}}^{D}\left(\left[X_{1}\right],\left[X_{2}\right]\right)   =  \frac{1}{2\pi}\left[\int_{0}^{2\pi}\omega\left(Z_{1},Y_{2}\right)dt+\omega\left(Y_{1},Z_{2}\right)\right]dt
   =   \frac{1}{2\pi}\int_{0}^{2\pi}\left[\mathcal{L}_{Y_{2}}h_{1}dt-\mathcal{L}_{Y_{1}}h_{2}\right]dt\\
\phantom{\Omega_{\text{red}}^{D}\left(\left[X_{1}\right],\left[X_{2}\right]\right)}{}   =   \Omega\left(\left[\left(Z_{1},-\mathcal{L}_{Y_{1}}\frac{dt}{2\pi}\right)\right],\left[\left(Z_{2},
-\mathcal{L}_{Y_{2}}\frac{dt}{2\pi}\right)\right]\right)
   =   \Omega\left(\left[\left(X_{1},0\right)\right],\left[\left(X_{2},0\right)\right]\right).
   \end{gather*}
\end{example}

\section{Weighted isotropic submanifolds}

The results of the last section suggest a way to generalize Weinstein's
construction of a symplectic structure on spaces of weighted Lagrangian
submanifolds to spaces of weighted isotropic submanifolds. In this
section we will pursue this idea and thereby obtain a generalization
of Theorem~\ref{sub:ThWisD} suggested by Proposition~\ref{sub:PropReductions}.

Let $\left(M,\omega\right)$ be a symplectic manifold and $I_{0}$
a f\/ixed compact oriented manifold with $\dim I_{0}\leq\frac{1}{2}\dim \left(M\right)$.

\subsection*{Notation}

Let $\text{Iso}\left(I_{0},M\right)$ denote the set of isotropic
embeddings of $I_{0}$ into $\left(M,\omega\right)$. That is,\[
\text{Iso}\left(I_{0},M\right):=\left\{ i\in\text{Emb}\left(I_{0},M\right)\mid i^{*}\omega=0\right\} .\]

Similar to the discussion in Proposition~\ref{sub:PropLag(L_0,M)}, Theorem~\ref{sub:ThWIsotropicET} can be used to equip $\text{Iso}\left(I_{0},M\right)$
with a smooth manifold structure locally modeled on spaces $\Gamma_{\text{closed}}\left(i^{*}TM\right)$
for $i\in\text{Iso}\left(I_{0},M\right)$. For each $i\in\text{Iso}\left(I_{0},M\right)$,
if $I:=i\left(I_{0}\right)$ then the sequence\[
\xymatrix{0\ar[r] & \Gamma\left(\left(TI\right)^{\omega}\right)\ar[r]^{f_{1}} & T_{i}\text{Iso}\left(I_{0},M\right)\ar[r]^{f_{2}} & Z^{1}\left(I_{0}\right)\ar[r] & 0}
\]
 where $f_{1}\left(Z\right)=Z\circ i$ and $f_{2}\left(X\right)=\alpha_{X}$,
is an exact sequence.

Meanwhile, the atlas on $\text{Iso}\left(I_{0},M\right)$ and the
spaces $\Gamma_{\text{exact}}\left(i^{*}TM\right)$ for $i\in\text{Iso}\left(I_{0},M\right)$
def\/ine a~foliation $\mathcal{E}$ of $\text{Iso}\left(I_{0},M\right)$,
whose leaves consist of $\text{Ham}\left(M\right)$ orbits under the
action\[
\text{Ham}\left(M\right)\circlearrowright\text{Iso}\left(I_{0},M\right):\ \ \phi\cdot i:=\phi\circ i.\]

\begin{definition}
The foliation $\mathcal{E}$ will be called the \emph{isodrastic foliation}
of $\text{Iso}\left(I_{0},M\right)$, and each leaf of $\mathcal{E}$
will be called an \emph{isodrast} in $\text{Iso}\left(I_{0},M\right)$.
\end{definition}

The group $\text{Dif\/f}_{+}\left(I_{0}\right)$ acts freely on $\text{Iso}\left(I_{0},M\right)$
via\[
a\cdot i:=i\circ a^{-1}.\]

\subsection*{Notation}

Set\[
\text{Iso}\left(M\right):=\text{Iso}\left(I_{0},M\right)/\text{Dif\/f}_{+}\left(I_{0}\right.)\]

If $p:\text{Iso}\left(I_{0},M\right)\to\text{Iso}\left(M\right)$
denotes the projection to the quotient, then $\text{Iso}\left(I_{0},M\right)$
is the total space of a principal $\text{Dif\/f}_{+}\left(I_{0}\right)$
bundle whose base $\text{Iso}\left(M\right)$ is identif\/ied with the
set of oriented isotropic submanifolds in $M$ dif\/feomorphic to $I_{0}$.

\subsection*{Notation}

Let $\text{Iso}_{w}\left(I_{0},M\right)$ denote the product $\text{Iso}\left(I_{0},M\right)\times\Omega_{1}^{k}\left(I_{0}\right)$
(where $k=\dim I_{0}$).

The group $\text{Dif\/f}_{+}\left(I_{0}\right)$ acts freely on $\text{Iso}_{w}\left(I_{0},M\right)$
by\[
a\cdot\left(i,\eta\right):=\big(i\circ a^{-1},\left(a^{-1}\right)^{*}\eta\big).\]

\subsection*{Notation}

Set\[
\text{Iso}_{w}\left(M\right):=\text{Iso}_{w}\left(I_{0},M\right)/\text{Dif\/f}_{+}\left(I_{0}\right).\]

The space $\text{Iso}_{w}\left(I_{0},M\right)$ is the total space
of a principal $\text{Dif\/f}_{+}\left(I_{0}\right)$ bundle, whose
base $\text{Iso}_{w}\left(M\right)$ can be identif\/ied with the space
of weighted isotropic submanifolds in $M$. The foliation $\mathcal{E}$
of $\text{Iso}\left(I_{0},M\right)$ induces a foliation $\mathcal{E}_{w}$
of $\text{Iso}_{w}\left(I_{0},M\right)$ which descends to a foliation
$\mathcal{F}_{w}$ of $\text{Iso}_{w}\left(M\right)$.

\begin{definition}
We will call $\mathcal{E}_{w}$ the \emph{isodrastic foliation} of
$\text{Iso}_{w}\left(I_{0},M\right)$, and each leaf of $\mathcal{E}_{w}$
will be called an \emph{isodrast} in $\text{Iso}_{w}\left(I_{0},M\right)$.
Similarly, $\mathcal{F}_{w}$ will be called the \emph{isodrastic
foliation} of $\text{Iso}_{w}\left(M\right)$, and each leaf of $\mathcal{F}_{w}$
will be called an \emph{isodrast} in $\text{Iso}_{w}\left(M\right)$.
\end{definition}

The pointwise assignment\[
\Omega_{\left(i,\eta\right)}\left(\left(X_{1},\vartheta_{1}\right),\left(X_{2},\vartheta_{2}\right)\right):=\int_{I_{0}}\left[\omega\left(X_{1},X_{2}\right)\eta+h_{1}\vartheta_{2}-h_{2}\vartheta_{1}\right]\]
 where $\alpha_{X_{k}}=dh_{k}$ for $h_{k}\in C^{\infty}\left(I_{0},\mathbb{R}\right)$
def\/ines a basic 2-form on $\text{Iso}_{w}\left(I_{0},M\right)^{\mathcal{E}_{w}}$.
It descends to a weak symplectic structure (also labeled $\Omega$)
on $\text{Iso}_{w}\left(I_{0},M\right)^{\mathcal{F}_{w}}$. Here closedness
and nondegeneracy of $\Omega$ on $\text{\textbf{I}}_{w}$ follow
from using the local model of isotropic submanifolds af\/forded by Theorem~\ref{sub:ThWIsotropicET}, and the fact that in such a model the symplectic
form along the zero section in $T^{*}I\oplus\Gamma\left(TI^{\omega}/TI\right)$
is given by $\omega_{T^{*}I}+\omega_{\mathbb{R}^{2n}}$ where $\omega_{\mathbb{R}^{2n}}$
is the standard symplectic vector space structure on the f\/ibers of
$\Gamma\left(TI^{\omega}/TI\right)$.

\subsection*{Notation}

Let $\text{Iso}_{pw}\left(I_{0},M\right)$ denote the product $\text{Iso}\left(I_{0},M\right)\times\text{Vol}_{1}\left(I_{0}\right)$
and set\[
\text{Iso}_{pw}\left(M\right):=\text{Iso}_{pw}\left(I_{0},M\right)/\text{Dif\/f}_{+}\left(I_{0}\right).\]

$\text{Iso}_{pw}\left(M\right)$ is an open subset of $\text{Iso}_{w}\left(M\right)$,
and so is a smooth manifold locally modeled on spaces $\Gamma_{\text{closed}}\left(i^{*}TM\right)\oplus\Omega_{0}^{k}\left(I_{0}\right)/\left\{ \left(Y,\mathcal{L}_{Y}\eta\right)\mid Y\in\mathfrak{X}\left(I_{0}\right)\right\} $
for $\left(i,\eta\right)\in\text{Iso}_{pw}\left(I_{0},M\right)$.
If $I_{0}$ is equipped with a f\/ixed volume form $\eta_{0}$ satisfying
$\int_{I_{0}}\eta_{0}=1$, then as in Proposition~\ref{PWQuotient} there
is a natural $\text{Ham}\left(M\right)$ equivariant dif\/feomorphism
$\upsilon:\text{Iso}\left(I_{0},M\right)/\text{Dif\/f}\left(I_{0},\eta_{0}\right)\to\text{Iso}_{pw}\left(M\right)$.
This leads to the following generalization of Theorem~\ref{sub:ThWisD}:

\begin{theorem}
Suppose that $\mathcal{N}$ is a $\emph{Ham}\left(M\right)$
orbit through an isotropic embedding, i.e. $\mathcal{N}$ is
an isodrast in $\emph{Iso}\left(S_{0},M\right)$. If $\emph{\textbf{I}}_{w}:=p\left(\mathcal{N}\times\emph{Vol}_{1}\left(S_{0}\right)\right)$
and $\mathcal{O}:=q\left(\mathcal{N}\right)$, where
$p:\emph{Iso}_{pw}\left(L_{0},M\right)\to\emph{Iso}_{pw}\left(M\right)$
and $q:\emph{Emb}\left(S_{0},M\right)\to\emph{Emb}\left(S_{0},M\right)/\emph{Dif\/f}\left(S_{0},\eta_{0}\right)$
are the projections to the quotients, then the reduction $\left(\mathcal{O},\Omega_{\emph{red}}^{D}\right)$
of $\left(\emph{Emb}\left(S_{0},M\right),\Omega^{D}\right)$ with
respect to the $\emph{Dif\/f}\left(S_{0},\eta_{0}\right)$ action
is symplectomorphic to the isodrast $\left(\emph{\textbf{I}}_{w},\Omega\right)$
in $\emph{Iso}_{pw}\left(M\right)$.
\end{theorem}

The proof is completely analogous to that of Theorem~\ref{sub:ThWisD}.

\section[Coadjoint orbits of $\text{Ham}\left(M\right)$]{Coadjoint orbits of $\boldsymbol{\text{Ham}\left(M\right)}$}

In this section we will extend a $\text{Ham}(M)$ moment map written
down by Weinstein (see Theorem~5.1 in \cite{Weinstein3}), which maps
each symplectic leaf $\left(\text{\textbf{I}}_{w},\Omega\right)\subset\text{Lag}_{pw}\left(M\right)$
onto a coadjoint orbit of a central extension of $\text{Ham}\left(M\right)$.
The result will be an identif\/ication of isodrasts consisting of positive
weighted isotropic submanifolds with coadjoint orbits of $\text{Ham}\left(M\right)$
(when $M$ is compact).

Suppose that $\left(M,\omega\right)$ is compact. The Lie algebra
of $G=\text{Ham}\left(M\right)$ can be identif\/ied with the space\[
\mathfrak{g}=\left\{ f\in C^{\infty}\left(M,\mathbb{R}\right)\mid\int_{M}f\frac{\omega^{n}}{n!}=0\right\} \simeq C^{\infty}\left(M,\mathbb{R}\right)/\mathbb{R}.\]
 Each positive weighted isotropic submanifold $\left(I,\rho\right)\in\text{Iso}_{pw}\left(M\right)$
def\/ines a element of the dual $\mathfrak{g}^{*}$ by integration:\[
f\mapsto\int_{I}f\rho.\]

\begin{theorem}
Let $\emph{\textbf{I}}_{w}$ be an isodrast in $\emph{Iso}_{pw}\left(M\right)$.
Then the $\emph{Ham}\left(M\right)$ action on $\left(\emph{\textbf{I}}_{w},\Omega\right)$
defined by $\phi\cdot\left(I,\rho\right):=\left(\phi\left(I\right),\left(\phi^{-1}\right)^{*}\rho\right)$
is Hamiltonian with moment map $\Phi:\emph{\textbf{I}}_{w}\to\mathfrak{g}^{*}$
given by\[
\left(I,\rho\right)\mapsto\left(f\mapsto\int_{I}f\rho\right).\]
 Moreover, the map $\Phi$ is a symplectomorphism onto
a coadjoint orbit in $\mathfrak{g}^{*}$ equipped with the
usual Kirillov--Kostant--Souriau form.
\end{theorem}

\begin{proof}
We will f\/irst
show that the $\text{Ham}\left(M\right)$ action is symplectic. If
$\left(I,\rho\right)=\left[\left(i,\eta\right)\right]$, then the
$\text{Ham}\left(M\right)$ action on $\text{\textbf{I}}_{w}$ can
be written as\[
\phi\cdot\left[\left(i,\eta\right)\right]:=\left[\left(\phi\circ i,\eta\right)\right].\]
 It follows that\begin{gather*}
\left(l_{\phi}^{*}\Omega\right)_{\left[\left(\phi\cdot i,\eta\right)\right]}\left(\left[\left(X_{1},\vartheta_{1}\right)\right],
\left[\left(X_{2},\vartheta_{2}\right)\right]\right)   =   \int_{I_{0}}\left[\omega\left(\text{d}\phi\cdot X_{1},\text{d}\phi\cdot X_{2}\right)\eta+h_{1}\vartheta_{2}-h_{2}\vartheta_{1}\right]\\
 \phantom{\left(l_{\phi}^{*}\Omega\right)_{\left[\left(\phi\cdot i,\eta\right)\right]}\left(\left[\left(X_{1},\vartheta_{1}\right)\right],
\left[\left(X_{2},\vartheta_{2}\right)\right]\right)}{}
   =   \int_{I_{0}}\left[\omega\left(X_{1},X_{2}\right)\eta+h_{1}\vartheta_{2}-h_{2}\vartheta_{1}\right] \\
 \phantom{\left(l_{\phi}^{*}\Omega\right)_{\left[\left(\phi\cdot i,\eta\right)\right]}\left(\left[\left(X_{1},\vartheta_{1}\right)\right],
\left[\left(X_{2},\vartheta_{2}\right)\right]\right)}{}
   =   \Omega_{\left[\left(i,\eta\right)\right]}\left(\left[\left(X_{1},\vartheta_{1}\right)\right],
   \left[\left(X_{2},\vartheta_{2}\right)\right]\right).
\end{gather*}

To check the moment map condition, we f\/irst note that each tangent
vector $\xi\in T_{\left[\left(i,\eta\right)\right]}\text{\textbf{I}}_{w}$
has a unique representative of the form $\left(X,0\right)$. If $f\in\mathfrak{g}$
then the generating vector f\/ield $f_{\text{\textbf{I}}_{w}}$ assigns
to each point $\left[\left(i,\eta\right)\right]\in\text{\textbf{I}}_{w}$
the tangent vector $\left[\left(X_{f}\circ i,0\right)\right]$ where
$X_{f}$ is the Hamiltonian vector f\/ield on $M$ corresponding to
the Hamiltonian $f$. If $\left[\left(X,0\right)\right]\in T_{\left[\left(i,\eta\right)\right]}\text{\textbf{I}}_{w}$
is a tangent vector to the isodrast, then choose an extension of $X$
to a Hamiltonian vector f\/ield $Z$ def\/ined on a~neighbourhood of $I=i\left(I_{0}\right)$.
It follows that\begin{gather*}
\text{d}\left\langle \Phi,f\right\rangle _{\left[\left(i,\eta\right)\right]}\left(\left[\left(X,0\right)\right]\right) =   \int_{I_{0}}i^{*}\mathcal{L}_{Z}f\eta
  =   \int_{I_{0}}i^{*}\omega\left(X_{f},Z\right)\eta\\
 \phantom{\text{d}\left\langle \Phi,f\right\rangle _{\left[\left(i,\eta\right)\right]}\left(\left[\left(X,0\right)\right]\right) }{}  =   \int_{I_{0}}\omega\left(X_{f}\circ i,X\right)\eta
  =   \Omega_{\left[\left(i,\eta\right)\right]}\left(f_{\text{\textbf{I}}_{w}},
  \left[\left(X,0\right)\right]\right).\end{gather*}

To check that the map $\Phi:\text{\textbf{I}}_{w}\to\mathfrak{g}^{*}$
is $\text{Ham}\left(M\right)$ equivariant we observe that\begin{gather*}
\left\langle \Phi\left(\phi\cdot\left[\left(i,\eta\right)\right]\right),f\right\rangle    =   \int_{I_{0}}i^{*}\phi^{*}f\eta
  =  \int_{I_{0}}i^{*}\left(f\circ\phi\right)\eta \\
  \phantom{\left\langle \Phi\left(\phi\cdot\left[\left(i,\eta\right)\right]\right),f\right\rangle}{}
   =   \left\langle \Phi\left(\left[\left(i,\eta\right)\right]\right),\text{Ad}_{\phi^{-1}}f\right\rangle
   =   \left\langle \text{Ad}_{\phi}^{*}\Phi\left(\left[\left(i,\eta\right)\right]\right),f\right\rangle .
   \end{gather*}

As for the image of $\Phi$, since the stabilizers at points $\left(I,\rho\right)\in\text{\textbf{I}}_{w}$
and $\Phi\left(I,\rho\right)\in\mathfrak{g}^{*}$ are given by\begin{gather*}
\text{Stab}_{\left(I,\rho\right)}\text{Ham}\left(M\right)   =   \big\{ \phi\in\text{Ham}\left(M\right)\mid\phi\left(I\right)=I,\:\left(\phi^{-1}\right)^{*}\rho=\rho\big\}
  =   \text{Stab}_{\Phi\left(I,\rho\right)}\text{Ham}\left(M\right),
  \end{gather*}
 the map $\Phi$ maps the isodrast $\text{\textbf{I}}_{w}$ dif\/feomorphically
onto a coadjoint orbit. Meanwhile, if $\left[\left(X_{1},0\right)\right]$
and $\left[\left(X_{2},0\right)\right]$ are tangent vectors in $T_{\left[\left(i,\eta\right)\right]}\text{\textbf{I}}_{w}$,
and $Z_{f_{1}}$ and $Z_{f_{2}}$ are extensions of $X_{1}$ and $X_{2}$
respectively to Hamiltonian vector f\/ields def\/ined on a neighbourhood
of $I=i\left(I_{0}\right)$, then the pullback under $\Phi$ of the
usual Kirillov--Kostant--Souriau form $\Omega_{\rm KKS}$ is given by
\begin{gather*}
\left(\Phi^{*}\Omega_{\rm KKS}\right)_{\left[\left(i,\eta\right)\right]}\left(\left[\left(X_{1},0\right)\right],\left[\left(X_{2},0\right)\right]\right) =-  \int_{I_{0}}i^{*}\left\{ f_{1},f_{2}\right\} \eta\\
\phantom{\left(\Phi^{*}\Omega_{\rm KKS}\right)_{\left[\left(i,\eta\right)\right]}\left(\left[\left(X_{1},0\right)\right],
\left[\left(X_{2},0\right)\right]\right)}{}
  =  \int_{I_{0}}i^{*}\omega\left(Z_{f_{1}},Z_{f_{2}}\right)\eta
  =  \Omega_{\left[\left(i,\eta\right)\right]}\left(\left[\left(X_{1},0\right)\right],
 \left[\left(X_{2},0\right)\right]\right).\!\!\!\!\tag*{\qed}
 \end{gather*}\renewcommand{\qed}{}
\end{proof}

\begin{remark}
This last theorem makes rigourous the heuristic moment map written
down by Weinstein. Also, it extends it to a map identifying isodrasts
in the space of positive weighted isotropic submanifolds with coadjoint
orbits in the dual of the Lie algebra of the group of Hamiltonian
symplectomorphisms.
\end{remark}

\section{Poisson structures}

In this section we will def\/ine a Poisson algebra $\left(\mathcal{A},\left\{ \cdot,\cdot\right\} \right)$
for a subalgebra $\mathcal{A}\subset C^{\infty}\left(M\right)$ of
smooth functions on a smooth manifold $M$. Such a manifold will be
called an $\mathcal{A}$-Poisson manifold if there are enough Hamiltonian
vector f\/ields in a sense we will specify. We def\/ine what a leaf of
such an $\mathcal{A}$-Poisson manifold is. We then show that reductions
of $\text{Emb}\left(S_{0},M\right)$ consisting of positive weighted
isotropic submanifolds are symplectic leaves of a Poisson structure.

\subsection*{Notation}

Let $C_{\int}^{\infty}\left(\text{Emb}\left(S_{0},M\right)\right)$
denote the set of maps $F:\text{Emb}\left(S_{0},M\right)\to\mathbb{R}$
such that for all $i_{0}\in\text{Emb}\left(S_{0},M\right)$ there
exists a $c^{\infty}$-open neighbourhood $U$ of $i_{0}$, a map
$A\in C^{\infty}\left(\mathbb{R}^{n},\mathbb{R}\right)$, and smooth
functions $h_{1},\ldots,h_{n}\in C^{\infty}\left(W,\mathbb{R}\right)$
def\/ined on a neighbourhood $W$ of $i_{0}\left(S_{0}\right)$ so that\begin{equation}
F\left(i\right)=A\left(\int_{S_{0}}\left(i^{*}h_{1}\right)\eta_{0},\ldots,\int_{S_{0}}\left(i^{*}h_{n}\right)\eta_{0}\right)\label{eq:FunctionsasIntegrals}\end{equation}
 for all $i\in U$.

Since functions in $C_{\int}^{\infty}\left(\text{Emb}\left(S_{0},M\right)\right)$
locally amount to integrating functions against $\eta_{0}$, the set
$C_{\int}^{\infty}\left(\text{Emb}\left(S_{0},M\right)\right)$ is
a subalgebra of $C^{\infty}(\text{Emb}(S_{0},M),\mathbb{R})$.

\begin{remark}
The algebra of functions $C_{\int}^{\infty}\left(\text{Emb}\left(S_{0},M\right)\right)$
was chosen with the axioms of dif\/fe\-ren\-tial structures in mind (see
\cite{Sikorski}; also cf.~Section~3 in \cite{CushmanSniatycki}).
A \emph{differential structure} on a~topological space $Q$ is a set
$C^{\infty}\left(Q\right)$ of continuous functions on $Q$ with the
following properties:
\begin{enumerate}\itemsep=0pt
\item The topology of $Q$ is generated by sets of the form $F^{-1}\left(V\right)$
where $F\in C^{\infty}\left(Q\right)$ and $V$ is an open subset
of $\mathbb{R}$.
\item If $B\in C^{\infty}\left(\mathbb{R}^{m},\mathbb{R}\right)$ and $F_{1},\ldots,F_{m}\in C^{\infty}\left(Q\right)$
then $B\left(F_{1},\ldots,F_{m}\right)\in C^{\infty}\left(Q\right)$.
\item If $F:Q\to\mathbb{R}$ is a function such that, for every $x\in Q$
there is an open neighbourhood $U$ of $x$ in $Q$ and a function
$F_{U}\in C^{\infty}\left(Q\right)$ such that $\left.F\right|_{U}=F_{U}$,
then $F\in C^{\infty}\left(Q\right)$.
\end{enumerate}
A topological space $Q$ together with a dif\/ferential structure $C^{\infty}\left(Q\right)$
is called a \emph{differential space}. When $\text{Emb}\left(S_{0},M\right)$
is equipped with the topology $\mathcal{T}$ generated by sets of
the form $F^{-1}\left(V\right)$ where $F\in C_{\int}^{\infty}\left(\text{Emb}\left(S_{0},M\right)\right)$
and $V$ is open in $\mathbb{R}$, then the topological space $\left(\text{Emb}\left(S_{0},M\right),\mathcal{T}\right)$
together with $C_{\int}^{\infty}\left(\text{Emb}\left(S_{0},M\right)\right)$
def\/ines a dif\/ferential space.
\end{remark}

\begin{proposition}\label{sub:Propv_FVectorfield}
For every $F\in C_{\int}^{\infty}\left(\emph{Emb}\left(S_{0},M\right)\right)$,
the local assignments given by $v_{F}\left(i\right):=\sum_{j=1}^{n}\frac{\partial A}{\partial y_{j}}\left.X_{h_{j}}\right|_{i\left(S_{0}\right)}\circ i$
on each neighbourhood $U$ define a unique vector field
$v_{F}$ on $\emph{Emb}\left(S_{0},M\right)$ satisfying
$dF=\imath\left(v_{F}\right)\Omega^{D}$.
\end{proposition}

\begin{proof} We will f\/irst
compute the exterior derivative of a map in $C_{\int}^{\infty}\left(\text{Emb}\left(S_{0},M\right)\right)$
locally. Suppose that $F\in C_{\int}^{\infty}\left(\text{Emb}\left(S_{0},M\right)\right)$
can be written as in \eqref{eq:FunctionsasIntegrals} on a $c^{\infty}$-open
neighbourhood $U$ of $i_{0}\in\text{Emb}\left(S_{0},M\right)$. Let
$X\in T_{i_{0}}\text{Emb}\left(S_{0},M\right)$ be a tangent vector
at $i_{0}$. Choose a vector f\/ield~$Z$ def\/ined on a neighbourhood
of $i\left(S_{0}\right)$ such that $Z\circ i=X$. Such a $Z$ can
be viewed as a vector f\/ield $v$ on $U$ (or perhaps on a $c^{\infty}$-open
subset $V\subset U$), assigning $Z\circ i$ to each $i\in U$. If\[
F\left(i\right)=A\left(\int_{S_{0}}\left(i^{*}h_{1}\right)\eta_{0},\ldots,\int_{S_{0}}\left(i^{*}h_{n}\right)\eta_{0}\right)\]
 on $U$ for some $A\in C^{\infty}\left(\mathbb{R}^{n},\mathbb{R}\right)$
and smooth functions $h_{1},\ldots,h_{n}\in C^{\infty}\left(W,\mathbb{R}\right)$
def\/ined on a~neighbourhood $W$ of $i_{0}\left(S_{0}\right)$, then
the pointwise exterior derivative of $F$ in the direction of $v$
is given by\[
dF_{i}\left(v_{i}\right)=\sum_{j=1}^{n}\frac{\partial A}{\partial y_{j}}\int_{S_{0}}\left(i^{*}\mathcal{L}_{Z}h_{j}\right)\eta_{0}.\]
 It follows that the vector f\/ield $v_{F}\in\mathfrak{X}\left(U\right)$,
def\/ined by $v_{F}\left(i\right):=\sum_{j=1}^{n}\frac{\partial A}{\partial y_{j}}X_{h_{j}}\circ i$
(where $X_{h_{j}}$ is the Hamiltonian vector f\/ield with Hamiltonian
$h_{j}$), satisf\/ies $dF=\imath\left(v_{F}\right)\Omega^{D}$ on $U$.
In fact, if $U_{1}$ and $U_{2}$ are $c^{\infty}$-open neighbourhoods
in $\text{Emb}\left(S_{0},M\right)$ with nonempty intersection, and\begin{gather*}
F\left(i\right)  =  A\left(\int_{S_{0}}\left(i^{*}h_{1}\right)\eta_{0},\ldots,\int_{S_{0}}\left(i^{*}h_{n}\right)\eta_{0}\right)\text{ on }U_{1},\\
F\left(i\right)  =  B\left(\int_{S_{0}}\left(i^{*}g_{1}\right)\eta_{0},\ldots,\int_{S_{0}}\left(i^{*}g_{m}\right)\eta_{0}\right)\text{ on }U_{2},\end{gather*}
 then on $U_{1}\cap U_{2}$
 \begin{gather*}
0   =   d\left(A\left(\int_{S_{0}}\left(i^{*}h_{1}\right)\eta_{0},\ldots,\int_{S_{0}}
\left(i^{*}h_{n}\right)\eta_{0}\right)
-B\left(\int_{S_{0}}\left(i^{*}g_{1}\right)\eta_{0},\ldots,\int_{S_{0}}\left(i^{*}g_{m}\right)\eta_{0}\right)\right)
\left(v\right)\\
 \phantom{0}{} =  \Omega^{D}\left(\sum_{j=1}^{n}\frac{\partial A}{\partial y_{j}}X_{h_{j}}\circ i-\sum_{l=1}^{m}\frac{\partial B}{\partial z_{l}}X_{g_{l}}\circ i,v\right).
 \end{gather*}
 So by the nondegeneracy of $\Omega^{D}$, the local assignments $v_{F}\left(i\right):=\sum_{j=1}^{n}\frac{\partial A}{\partial y_{j}}X_{h_{j}}\circ i$
coincide on overlapping regions and so def\/ine a vector f\/ield $v_{F}$
on $\text{Emb}\left(S_{0},M\right)$ satisfying $dF=\imath\left(v_{F}\right)\Omega^{D}$.
Uniqueness of the Hamiltonian vector f\/ield $v_{F}$ also follows from
nondegeneracy of $\Omega^{D}$.
\end{proof}

\begin{definition}
Let $M$ be a smooth manifold. A subalgebra $\mathcal{A}\subset C^{\infty}\left(M,\mathbb{R}\right)$
of smooth functions together with a Lie structure $\left[\cdot,\cdot\right]$
will be called a \emph{Poisson algebra} if\[
\left[FG,H\right]=F\left[G,H\right]+\left[F,H\right]G.\]
 In this case the bracket $\left[\cdot,\cdot\right]$ will be called
a \emph{Poisson bracket}. If $\left(\mathcal{A},\left[\cdot,\cdot\right]\right)$
is a Poisson algebra, we will say that $M$ is an $\mathcal{A}$-\emph{Poisson
manifold}.
\end{definition}

Def\/ine a skew-symmetric bilinear structure $\left\{ \cdot,\cdot\right\} _{\text{Emb}\left(S_{0},M\right)}$
on $C_{\int}^{\infty}\left(\text{Emb}\left(S_{0},M\right)\right)$
via\[
\left\{ F,G\right\} _{\text{Emb}\left(S_{0},M\right)}:=-\Omega^{D}\left(v_{F},v_{G}\right).\]

\begin{proposition}
The subalgebra $C_{\int}^{\infty}\left(\emph{Emb}\left(S_{0},M\right)\right)$
together with $\left\{ \cdot,\cdot\right\} _{\emph{Emb}\left(S_{0},M\right)}$
is a Poisson algebra.
\end{proposition}

\begin{proof}
$\left\{ \cdot,\cdot\right\} _{\text{Emb}\left(S_{0},M\right)}$
satisf\/ies Jacobi's identity pointwise since the Poisson bracket $\left\{ \cdot,\cdot\right\} _{M}$
on $M$ does, and is a derivation in the f\/irst argument for the same
reason. (Here $\left\{ \cdot,\cdot\right\} _{M}$ is def\/ined by $\left\{ f,g\right\} _{M}:=-\omega\left(X_{f},X_{g}\right)$
for $f,g\in C^{\infty}\left(M,\mathbb{R}\right)$.)
\end{proof}

By
taking restrictions of functions in $C_{\int}^{\infty}\left(\text{Emb}\left(S_{0},M\right)\right)$
to $\text{Iso}\left(S_{0},M\right)$ we obtain an algebra $C_{\int}^{\infty}\left(\text{Iso}\left(S_{0},M\right)\right)$
of smooth functions on $\text{Iso}\left(S_{0},M\right)$. For each
$F\in C_{\int}^{\infty}\left(\text{Iso}\left(S_{0},M\right)\right)$
and $i\in\text{Iso}\left(S_{0},M\right)$ the vector $v_{F}\left(i\right)$
is tangent to $\text{Iso}\left(S_{0},M\right)$ at $i$. Since the
functions in $C_{\int}^{\infty}\left(\text{Iso}\left(S_{0},M\right)\right)$
are $\text{Dif\/f}\left(S_{0},\eta_{0}\right)$ invariant, the algebra
$C_{\int}^{\infty}\left(\text{Iso}\left(S_{0},M\right)\right)$ descends
to an algebra of functions $C_{\int}^{\infty}\left(\text{Iso}_{pw}\left(M\right)\right)$
on $\text{Iso}_{pw}\left(M\right)$.

\subsection*{Convention}

In what follows, we will use the same letter $F$ to denote a function
in $C_{\int}^{\infty}\left(\text{Emb}\left(S_{0},M\right)\right)$,
its restriction to $\text{Iso}\left(S_{0},M\right)$, as well as the
corresponding map on the quotient $\text{Iso}_{pw}\left(M\right)$
to avoid introducing more notation.

For $F,G\in C_{\int}^{\infty}\left(\text{Iso}_{pw}\left(M\right)\right)$
and $\left[i\right]\in\text{Iso}_{pw}\left(M\right)$ def\/ine\begin{gather*}
\left\{ F,G\right\} _{\text{Iso}_{pw}\left(M\right)}\left(\left[i\right]\right):   =   -\Omega_{i}^{D}\left(v_{F}\left(i\right),v_{G}\left(i\right)\right)
   =   -\left(\Omega_{\text{red}}^{D}\right)_{\left[i\right]}\left(\text{d}r
   \left(v_{F}\right),\text{d}r\left(v_{G}\right)\right),\end{gather*}
 where $r:\text{Iso}\left(S_{0},M\right)\to\text{Iso}_{pw}\left(M\right)$
is the projection to the quotient. Then the subalgebra $C_{\int}^{\infty}\left(\text{Iso}_{pw}\left(M\right)\right)$
together with $\left\{ \cdot,\cdot\right\} _{\text{Iso}_{pw}\left(M\right)}$
is a Poisson algebra, which means $\text{Iso}_{pw}\left(M\right)$
is a~$C_{\int}^{\infty}\left(\text{Iso}_{pw}\left(M\right)\right)$-Poisson
manifold.

\begin{definition}
Let $\left(M,\mathcal{A},\left\{ \cdot,\cdot\right\} \right)$ be
an $\mathcal{A}$-Poisson manifold. Suppose that for every $F\in\mathcal{A}$
there exists a unique vector f\/ield $X_{F}$ on $M$ such that $\text{d}G\left(X_{F}\right)=\left\{ F,G\right\} $
for every $G\in\mathcal{A}$. The collection of vectors\[
\mathcal{C}_{x}:=\left\{ X_{F}\left(x\right)\in T_{x}M\mid F\in\mathcal{A}\right\} \]
 at each $x\in M$ def\/ines a distribution $\mathcal{C}$ on $M$.
A maximal integral manifold $\mathcal{M}$ of $\mathcal{C}$ will
be called a \emph{leaf} of the $\mathcal{A}$-Poisson manifold $M$.
\end{definition}

\begin{remark}
In f\/inite dimensions the distribution $\mathcal{C}$ is integrable.
In inf\/inite dimensions it might be possible to prove an analogue of
Frobenius' theorem for some spaces of mappings and apply it to our
spaces.
\end{remark}

\begin{remark}
The f\/irst condition in the previous def\/inition is satisf\/ied within the
subal\-geb\-ra $C_{\int}^{\infty}\left(\text{Iso}_{pw}\left(M\right)\right)$
but not in $C_{\int}^{\infty}\left(\text{Emb}\left(S_{0},M\right)\right)$
because of the directions generated by the $\text{Dif\/f}(S_{0},\eta_{0})$
action.
\end{remark}

For each $F\in C_{\int}^{\infty}\left(\text{Iso}_{pw}\left(M\right)\right)$
we have that $\text{d}r\left(v_{F}\right)$ is the unique vector f\/ield
on $\text{Iso}_{pw}\left(M\right)$ satis\-fying $dF=\imath\left(\text{d}r\left(v_{F}\right)\right)\Omega_{\text{red}}^{D}$
(uniqueness follows from the nondegeneracy of $\Omega_{\text{red}}^{D}$
\linebreak on~$\text{Iso}_{pw}\left(M\right)^{\mathcal{F}_{pw}}$). It follows
that if $F\in C_{\int}^{\infty}\left(\text{Iso}_{pw}\left(M\right)\right)$
then $\text{d}G\left(X_{F}\right)=\left\{ F,G\right\} _{\text{Iso}_{pw}\left(M\right)}$
for all $G\in C_{\int}^{\infty}\left(\text{Iso}_{pw}\left(M\right)\right)$.
On the $C_{\int}^{\infty}\left(\text{Iso}_{pw}\left(M\right)\right)$-Poisson
manifold $\text{Iso}_{pw}\left(M\right)$, the distribu\-tion~$\mathcal{C}$
is given by the collection of vectors
\[
\mathcal{C}_{\left[i\right]}:=\left\{ \text{d}r\left(v_{F}\right)\left(\left[i\right]\right)\in T_{\left[i\right]}\text{Iso}_{pw}\left(M\right)\mid F\in C_{\int}^{\infty}\left(\text{Iso}_{pw}\left(S_{0},M\right)\right)\right\} \]
 at each $\left[i\right]\in\text{Iso}_{pw}\left(M\right)$.

\begin{theorem}
Suppose that $\mathcal{N}$ is a $\emph{Ham}(M)$
orbit through an isotropic embedding in $\emph{Emb}(S_{0},\eta_{0})$,
i.e., $\mathcal{N}$ is an isodrast in $\emph{Iso}\left(S_{0},M\right)$,
with $\mathcal{O}:=q\left(\mathcal{N}\right)$ where the mapping
$q:\emph{Emb}\left(S_{0},M\right)\to\emph{Emb}\left(S_{0},M\right)/\emph{Dif\/f}\left(S_{0},\eta_{0}\right)$
is the projection to the quotient. Then the reduction
$\left(\mathcal{O},\Omega_{\emph{red}}^{D}\right)$ is a~symplectic
leaf of the $C_{\int}^{\infty}\left(\emph{Iso}_{pw}\left(M\right)\right)$-Poisson
manifold $\emph{Iso}_{pw}\left(M\right)$.
\end{theorem}

\begin{proof} We will
f\/irst show that any maximal integral manifold of $\mathcal{C}$ is
a $\text{Ham}\left(M\right)$ orbit in $\text{Iso}_{pw}\left(M\right)$
under the action $\phi\cdot\left[i\right]:=\left[\phi\circ i\right]$.
The manifold $\text{Iso}\left(S_{0},M\right)$ is the total space
of a principal $\text{Dif\/f}\left(S_{0},M\right)$ bundle over $\text{Iso}_{pw}\left(M\right)$.
Each vector f\/ield $v_{F}$ on $\text{Iso}\left(S_{0},M\right)$ with
Hamiltonian $F\in C_{\int}^{\infty}\left(\text{Iso}\left(S_{0},M\right)\right)$
is locally of the form $v_{F}\left(i\right)=\sum_{j=1}^{n}\frac{\partial A}{\partial y_{j}}X_{h_{j}}\circ i$.
It follows that curves in $\text{Iso}_{pw}\left(M\right)$ are everywhere
tangent to the distribution $\mathcal{C}$ if and only if they belong
to a $\text{Ham}\left(M\right)$ orbit. So the symplectic quotient
$\left(\mathcal{O}:=q\left(\mathcal{N}\right),\Omega_{\text{red}}^{D}\right)$
is a leaf of the $C_{\int}^{\infty}\left(\text{Iso}_{pw}\left(M\right)\right)$-Poisson
manifold $\text{Iso}_{pw}\left(M\right)$.
\end{proof}

\begin{remark}
By now we have revealed three dif\/ferent faces of isodrasts in $\text{Iso}_{pw}\left(M\right)$.
Namely, isodrasts consisting of positively weighted isotropic submanifolds
can be identif\/ied with symplectic reductions of $\left(\text{Emb}\left(S_{0},M\right),\Omega^{D}\right)$,
with coadjoint orbits of $\text{Ham}\left(M\right)$, and with symplectic
leaves of the $C_{\int}^{\infty}\left(\text{Iso}_{pw}\left(M\right)\right)$-Poisson
manifold $\text{Iso}_{pw}\left(M\right)$.
\end{remark}

\section{Momentum weighted Lagrangian submanifolds\label{sec:Momentum-Weighted-Lagrangians}}

In this section we will discuss a kinematic interpretation of isodrasts
in $\text{Lag}\left(M\right)$ to motivate a~dif\/ferent choice of ``weightings''
to obtain a symplectic structure.

Let us view points in an isodrast $\text{\textbf{I}}\subset\text{Lag}\left(M\right)$
as conf\/igurations of a submanifold constrained to move in $\text{\textbf{I}}$.
What are the possible velocities? By Proposition~\ref{sub:PropTLag(M)}
and the description of the isodrastic foliation $\mathcal{F}$, we
know that the tangent bundle of $\text{\textbf{I}}$ can be described
by\[
T\text{\textbf{I}}=\left\{ \left(L,dh\right)\in\text{\textbf{I}}\times B^{1}\left(L\right)\mid L\in\text{\textbf{I}}\right\} .\]
 Thus the velocities at a conf\/iguration $L\in\text{\textbf{I}}$ correspond
to functions on $L$ modulo constants, i.e.\ to elements of $C^{\infty}\left(L,\mathbb{R}\right)/\mathbb{R}$.
The conjugate momenta to conf\/igurations in $\text{\textbf{I}}$ (i.e.\
cotangent vectors) should be linear functionals from $C^{\infty}\left(L,\mathbb{R}\right)/\mathbb{R}$
to $\mathbb{R}$, that are in 1-1 correspondence with $C^{\infty}\left(L,\mathbb{R}\right)/\mathbb{R}$.
This expectation that the cotangent f\/ibers should be the ``same
size'' as $C^{\infty}\left(L,\mathbb{R}\right)/\mathbb{R}$
ref\/lects a physical expectation that all momenta should be accessible
by motions of particles in the conf\/iguration space, and that all such
motions can be assigned momenta. In any case, integration against
$n$-forms in $\Omega_{0}^{n}\left(L\right)$ certainly f\/its the above
description. This motivates the following def\/inition:

\begin{definition}
Let $\text{\textbf{I}}$ be an isodrast in $\text{Lag}\left(M\right)$.
A \emph{momentum weighting} of a Lagrangian submanifold $L\in\text{\textbf{I}}$
is a top degree form $\chi$ on $L$ satisfying $\int_{L}\chi=0$.
Pairs $\left(L,\chi\right)$ will be called \emph{momentum weighted
Lagrangian submanifolds}.
\end{definition}

\subsection*{Notation}

Let $\text{Lag}_{mw}\left(L_{0},M\right)$ denote the product $\text{Lag}\left(L_{0},M\right)\times\Omega_{0}^{n}\left(L_{0}\right)$
and $\text{Lag}_{mw}\left(M\right)$ the quotient of $\text{Lag}_{mw}\left(L_{0},M\right)$
by $\text{Dif\/f}_{+}\left(L_{0}\right)$ under the action $a\cdot\left(i,\chi\right):=\left(i\circ a^{-1},\left(a^{-1}\right)^{*}\chi\right)$.

As in the case of weighted Lagrangian submanifolds, the canonical
projection from $\text{Lag}_{mw}\left(M\right)$ to $\text{Lag}\left(M\right)$
pulls back the foliation $\mathcal{F}$ to a foliation $\mathcal{F}_{mw}$
on $\text{Lag}_{mw}\left(M\right)$, whose leaves $\text{\textbf{I}}_{mw}$
are $\text{Ham}\left(M\right)$ orbits under the action\[
\text{Ham}\left(M\right)\circlearrowright\text{Lag}_{mw}\left(M\right):\ \phi\cdot\left[\left(i,\chi\right)\right]:=\left[\left(\phi\circ i,\chi\right)\right].\]

\begin{definition}
We will call $\mathcal{F}_{mw}$ the \emph{isodrastic foliation} of
$\text{Lag}_{mw}\left(M\right)$ and each leaf of $\mathcal{F}_{mw}$
will be called an \emph{isodrast} in $\text{Lag}_{mw}\left(M\right)$.
\end{definition}

The tangent space to an isodrast $\text{\textbf{I}}_{mw}$ at a point
$\left[\left(i,\chi\right)\right]$ is given by\begin{gather*}
T_{\left[\left(i,\nu\right)\right]}\text{\textbf{I}}_{mw}   =   \left\{ \left(X,\vartheta\right)\in\Gamma\left(i^{*}TM\right)\times\Omega_{0}^{n}\left(L_{0}\right)\mid\alpha_{X}\in B^{1}\left(L_{0}\right)\right\}
 /\left\{ \left(Y,\mathcal{L}_{Y}\eta\right)\mid Y\in\mathfrak{X}\left(L_{0}\right)\right\} .\end{gather*}

\begin{theorem}
The bilinear map\[
\Omega_{\left[\left(i,\chi\right)\right]}\left(\left[\left(X_{1},\vartheta_{1}\right)\right],\left[\left(X_{2},\vartheta_{2}\right)\right]\right):=\int_{L_{0}}\left[\omega\left(X_{1},X_{2}\right)\chi+h_{1}\vartheta_{2}-h_{2}\vartheta_{1}\right]\]
 where $\alpha_{X_{k}}=dh_{k}$, defines an exact weak symplectic
structure on $\emph{\textbf{I}}_{mw}$ satisfying $\Omega=-d\Theta$
where\[
\Theta_{\left[\left(i,\chi\right)\right]}\left(\left[\left(X,\vartheta\right)\right]\right):=\int_{L_{0}}h\chi\]
 and $\alpha_{X}=dh$.
\end{theorem}

\begin{proof} We will compute the exterior
derivative of $\Theta$ locally in charts. That is, be means of a
symplectic embedding $\lambda_{\left(i,\chi\right)}:T^{*}L_{0}\supset W_{\left(i,\chi\right)}\to M$
chosen in def\/ining a chart $\left(U_{\left(i,\chi\right)},\varphi_{\left(i,\chi\right)}\right)$
on $\text{Lag}_{mw}\left(M\right)$, we can assume that $M=T^{*}L_{0}$
and that each tangent vector in $T_{\left[\left(i,\chi\right)\right]}\text{\textbf{I}}_{mw}$
is represented by a pair $\left(X,\vartheta\right)$ where $\alpha_{X}\in B^{1}\left(L_{0}\right)$
and $X$ is tangent to the cotangent f\/ibers. Given a tangent vector
$\xi=\left(X,\vartheta\right)\in T_{\left[\left(i,\chi\right)\right]}\text{\textbf{I}}_{mw}$,
we can extend it to a vector f\/ield on $\text{\textbf{I}}_{mw}$ (also
labeled $\xi$) in the following way. Let $Z_{f}$ denote the Hamiltonian
vector f\/ield def\/ined on $T^{*}L_{0}$ satisfying $\imath\left(Z_{f}\right)\omega=\pi^{*}\alpha_{X}$
where $\alpha_{X}=dh$ and $f=\pi^{*}h$. Def\/ine $\xi$ to be the
vector f\/ield on $\text{\textbf{I}}_{mw}$ that assigns $\xi\left(\left[i',\nu'\right]\right)=\left(Z_{f}\circ i',\vartheta\right)$
to nearby points $\left[\left(i',\chi'\right)\right]$. So given tangent
vectors~$\xi_{1}$ and~$\xi_{2}$ in $T_{\left[\left(i,\chi\right)\right]}\text{\textbf{I}}_{mw}$,
extend them to vector f\/ields (also labeled $\xi_{1}$ and~$\xi_{2}$)
as just described. Then\begin{gather*}
d\Theta_{\left[\left(i,\nu\right)\right]}\left(\xi_{1},\xi_{2}\right)   =   \left[\xi_{1}\Theta\left(\xi_{2}\right)-\xi_{2}\Theta\left(\xi_{1}\right)
-\Theta\left(\left[\xi_{1},\xi_{2}\right]\right)\right]_{\left[\left(i,\nu\right)\right]}\\
 \phantom{d\Theta_{\left[\left(i,\nu\right)\right]}\left(\xi_{1},\xi_{2}\right)}{}
 =   \int_{L_{0}}\Big[\big(i^{*}L_{X_{f_{1}}}f_{2}\big)\chi+\left(i^{*}f_{2}\right)
  \vartheta_{1}-\big(i^{*}L_{X_{f_{2}}}f_{1}\big)\chi-\left(i^{*}f_{1}\right)\vartheta_{2} \\
\phantom{d\Theta_{\left[\left(i,\nu\right)\right]}\left(\xi_{1},\xi_{2}\right)=}{}
 +\left(i^{*}\omega\left(X_{f_{1}},X_{f_{2}}\right)\right)\chi\Big]
   =   -\int_{L_{0}}\left(h_{1}\vartheta_{2}-h_{2}\vartheta_{1}\right)
  =   -\Omega_{\left[\left(i,\chi\right)\right]}\left(\xi_{1},\xi_{2}\right).
  \end{gather*}

Nondegeneracy then follows from the fact that the $h$'s and $\vartheta$'s
are independent of one \linebreak another.
\end{proof}

\section{Momentum weighted metrics}

In this section we apply the kinematic approach of the last section
to the space of pseudo-Riemannian metrics of f\/ixed signature. We obtain
a cotangent bundle and identify a Poisson algebra of functions on
this space. These functions are Hamiltonians for Hamiltonian vector
f\/ields. In the particular case of Lorentzian metrics on a f\/ixed 4-dimensional
manifold, this gives some ingredients for a possible geometric quantization
of gravity in a vacuum.

\subsection*{Convention}

From now on $M$ will no longer necessarily be symplectic. Instead
we will simply assume that~$M$ is a f\/inite dimensional manifold.

\subsection*{Notation}

Let $\text{Met}^{q}\left(M\right)$ denote the set of pseudo Riemannian
metrics on $M$ of signature $q$ and $\Gamma_{c}\left(S^{2}T^{*}M\right)$
the set of compactly supported symmetric two tensors on $M$.

Then $\text{Met}^{q}\left(M\right)$ is a smooth manifold modeled
on the space $\Gamma_{c}\left(S^{2}T^{*}M\right)$ (see Theorem~45.13
in \cite{KrieglMichor}). By identifying such metrics with their graphs
in $\Gamma\left(L\left(TM,T^{*}M\right)\right)$, we can view $\text{Met}^{q}\left(M\right)$
as a collection of submanifolds. If as in Section~\ref{sec:Momentum-Weighted-Lagrangians}
we view $\text{Met}^{q}\left(M\right)$ kinematically as the set of
possible conf\/igurations of a submanifold moving in $\Gamma\left(L\left(TM,T^{*}M\right)\right)$,
then by analogy we should be able to describe the ``conjugate momenta''
in terms of some kind of ``weightings'' on the submanifolds.

To determine what these weightings should be, we f\/irst note that since
the submanifolds in $\Gamma\left(L\left(TM,T^{*}M\right)\right)$
corresponding to pseudo-Riemannian metrics in $\text{Met}^{q}\left(M\right)$
are graphs, any sections of bundles over such manifolds can be obtained
as pullups of sections over corresponding bundles on the base $M$.
Second, the weightings should give linear functionals on each tangent
space to $\text{Met}^{q}\left(M\right)$ and be in 1-1 correspondence
with $\Gamma_{c}\left(S^{2}T^{*}M\right)$. This motivates the following
def\/inition:

\begin{definition}
A \emph{momentum weighting} of a pseudo Riemannian metric $g\in\text{Met}^{q}\left(M\right)$
of signature $q$ is a symmetric two tensor $h\in\Gamma_{c}\left(S^{2}T^{*}M\right)$.
Pairs $\left(g,h\right)$ will be called \emph{momentum weighted metrics}.
\end{definition}

\subsection*{Notation}

Let $\text{Met}_{mw}^{q}\left(M\right)$ denote the set of momentum
weighted metrics of signature $q$. That is,\[
\text{Met}_{mw}^{q}\left(M\right):=\text{Met}^{q}\left(M\right)\times\Gamma_{c}\left(S^{2}T^{*}M\right).\]

We will show that $\text{Met}_{mw}^{q}\left(M\right)$ is the appropriate
phase space for $\text{Met}^{q}\left(M\right)$. Def\/ine a 1-form~$\Theta$
on $\text{Met}_{mw}^{q}\left(M\right)$ via \[
\Theta_{\left(g,h\right)}\left(k,l\right):=\int_{M}\text{Tr}\left(g^{-1}kg^{-1}h\right)\mu\left(g\right),\]
 where $\text{Tr}\left(g^{-1}kg^{-1}h\right)=g^{jr}k_{rq}g^{qp}h_{pj}$,
$\mu\left(g\right)=\sqrt{\left|\det g\right|}\left|dx_{1}\wedge\cdots\wedge dx_{n}\right|$.

We will now explicitly calculate the exterior derivative of $\Theta$.
Given a pair of tangent vectors $\left(k_{1},l_{1}\right),\left(k_{2},l_{2}\right)\!\in\! T_{\left(g,h\right)}\text{Met}_{mw}^{q}\!(M)$,
extend each $\left(k_{p},l_{p}\right)$ to a vector f\/ield $\xi_{p}$
on $\text{Met}_{mw}^{q}\!\left(M\right)$ satisfying $\xi_{p}\left(\left(g',h'\right)\right)=\left(k_{p},l_{p}\right)$
for all $\left(g',h'\right)$ near $\left(g,h\right)$. Since as pointwise
matrices,\begin{gather*}
\left.\frac{d}{dt}\right|_{t=0}\left(g+tk\right)^{-1}   =   -g^{-1}kg^{-1},\\
\left.\frac{d}{dt}\right|_{t=0}\sqrt{\left|\det \left(g+tk\right)\right|}   =   \sqrt{\left|\det g\right|}\text{Tr}\left(g^{-1}k\right),
\end{gather*}
 it follows that\begin{gather*}
d\Theta_{\left(g,h\right)}\left(\left(k_{1},l_{1}\right),\left(k_{2},l_{2}\right)\right)   =   \left[d\Theta\left(\xi_{1},\xi_{2}\right)\right]_{\left(g,h\right)} \\
\phantom{d\Theta_{\left(g,h\right)}\left(\left(k_{1},l_{1}\right),\left(k_{2},l_{2}\right)\right)}{}
   =   \int_{M}\left[-\text{Tr}\left(g^{-1}k_{1}g^{-1}k_{2}g^{-1}h\right)
   -\text{Tr}\left(g^{-1}k_{2}g^{-1}k_{1}g^{-1}h\right)\right.\\
\left. \phantom{d\Theta_{\left(g,h\right)}\left(\left(k_{1},l_{1}\right),\left(k_{2},l_{2}\right)\right)=}{}
 +\text{Tr}\left(g^{-1}k_{2}g^{-1}h\right)\text{Tr}\left(g^{-1}k_{1}\right)\!
 +\text{Tr}\left(g^{-1}k_{2}g^{-1}l_{1}\right)\right]\mu\left(g\right)
    -\left(1\!\leftrightarrow \!2\right)\!\!\\
\phantom{d\Theta_{\left(g,h\right)}\left(\left(k_{1},l_{1}\right),\left(k_{2},l_{2}\right)\right)}{}
   =   \int_{M}\left[\text{Tr}\left(g^{-1}k_{2}g^{-1}h\right)\text{Tr}\left(g^{-1}k_{1}\right) -\text{Tr}\left(g^{-1}k_{1}g^{-1}h\right)\text{Tr}\left(g^{-1}k_{2}\right)\right.\\
\left.\phantom{d\Theta_{\left(g,h\right)}\left(\left(k_{1},l_{1}\right),\left(k_{2},l_{2}\right)\right)=}{}
+\text{Tr}\left(g^{-1}k_{2}g^{-1}l_{1}\right)-\text{Tr}\left(g^{-1}k_{1}g^{-1}l_{2}\right)\right]
 \mu\left(g\right).
 \end{gather*}

\begin{theorem}
The $2$-form $\Omega:=-d\Theta$ defined by
\begin{gather*}
\Omega_{\left(g,h\right)}\left(\xi_{1},\xi_{2}\right)   =   \int_{M}\left[\emph{Tr}\left(g^{-1}k_{1}g^{-1}h\right)\emph{Tr}\left(g^{-1}k_{2}\right)
-\emph{Tr}\left(g^{-1}k_{2}g^{-1}h\right)\emph{Tr}\left(g^{-1}k_{1}\right)\right.\\
\left.\phantom{\Omega_{\left(g,h\right)}\left(\xi_{1},\xi_{2}\right)   =}{}  +\emph{Tr}\left(g^{-1}k_{1}g^{-1}l_{2}\right)
-\emph{Tr}\left(g^{-1}k_{2}g^{-1}l_{1}\right)\right]\mu\left(g\right)
\end{gather*}
 defines an exact symplectic structure on $\emph{Met}_{mw}^{q}\left(M\right)$.
\end{theorem}

\begin{proof} If $k_{1}$ is nonzero, then by taking $k_{2}=0$ the
expression for $\Omega_{\left(g,h\right)}\left(\xi_{1},\xi_{2}\right)$
reduces to $\int_{M}\text{Tr}\left(g^{-1}k_{1}g^{-1}l_{2}\right)\mu\left(g\right)$.
Since\[
G_{g}\left(k_{1},k_{2}\right)=\int_{M}\text{Tr}\left(g^{-1}k_{1}g^{-1}k_{2}\right)\mu\left(g\right)\]
 is weakly nondegenerate (see Lemma 45.3 in \cite{KrieglMichor}),
$l_{2}$ can be chosen so that the integral does not vanish.

If $k_{1}=0$, then we are left with $-\int_{M}\text{Tr}\left(g^{-1}k_{2}g^{-1}l_{1}\right)\mu\left(g\right)$
which means we can choose $k_{2}$ so that the integral does not vanish
for the same reason.
\end{proof}

\subsection*{Notation}

For every section $r\in\Gamma_{c}\left(S^{2}T^{*}M\right)$, def\/ine
the function
\begin{gather*}
F_{r}   :  \  \text{Met}_{mw}^{q}(M)\to\mathbb{R}, \qquad
F_{r}\left(g,h\right):   =   \int_{M}\text{Tr}\left(g^{-1}rg^{-1}h\right)\mu\left(g\right).
\end{gather*}

\begin{proposition}
For every $r\in\Gamma_{c}\left(S^{2}T^{*}M\right)$, the
vector field $\xi_{F_{r}}$ defined by\[
\xi_{F_{r}}\left(g,h\right):=\left(r,rg^{-1}h+hg^{-1}r-\emph{Tr}\left(g^{-1}r\right)\cdot h\right)\]
 is a Hamiltonian vector field on $\emph{Met}_{mw}^{q}\left(M\right)$
corresponding to the Hamiltonian function $F_{r}$, i.e.
$\imath\left(\xi_{F_{r}}\right)\Omega=dF_{r}$.
\end{proposition}

\begin{proof} Given
$(k,l)\!\in\! T_{(g,h)}\text{Met}_{mw}^{q}\!(M)$,
extend $(k,l)$ constantly to a vector f\/ield $\xi$ on
$\text{Met}_{mw}^{q}\!(M)$ satisfying $\xi\left(g',h'\right)=\left(k,l\right)$.
Then\begin{gather*}
\left(dF_{r}\right)_{\left(g,h\right)}\left(\xi\right)   =   \int_{M}\Big[-\text{Tr}\left(g^{-1}kg^{-1}rg^{-1}h\right)-\text{Tr}\left(g^{-1}rg^{-1}kg^{-1}h\right) \\
 \phantom{\left(dF_{r}\right)_{\left(g,h\right)}\left(\xi\right)   =}{}  +\text{Tr}\left(g^{-1}rg^{-1}h\right)\text{Tr}\left(g^{-1}k\right)
+\text{Tr}\left(g^{-1}rg^{-1}l\right)\Big]\mu\left(g\right)
   =   \Omega_{\left(g,h\right)}\left(\xi_{F_{r}},\xi\right). \tag*{\qed}
   \end{gather*}
 \renewcommand{\qed}{}
\end{proof}

\subsection*{Notation}

The manifold topology on $\text{Met}_{mw}^{q}\left(M\right)$ is f\/iner
than the trace of the Whitney $C^{\infty}$-topology on $\Gamma\left(L\left(TM,T^{*}M\right)\right)$
(see Section~45.1 in \cite{KrieglMichor}). Let $C_{\int}^{\infty}\left(\text{Met}_{mw}^{q}\left(M\right)\right)$
denote the set of functions $F:\text{Met}_{mw}^{q}\left(M\right)\to\mathbb{R}$
such that for every $g_{0}\in\text{Met}_{mw}^{q}\left(M\right)$ there
exists a neighbourhood $U$ of $g_{0}$, a map $A:\mathbb{R}^{n}\to\mathbb{R}$,
and sections $r_{1},\ldots,r_{n}\in\Gamma_{c}\left(S^{2}T^{*}M\right)$
so that\[
F\left(g\right)=A\left(\int_{M}\text{Tr}\left(g^{-1}r_{1}g^{-1}h\right)\mu\left(g\right),\ldots,\int_{M}\text{Tr}\left(g^{-1}r_{n}g^{-1}h\right)\mu\left(g\right)\right)\]
 for all $g\in U$.

\begin{remark}
The algebra $C_{\int}^{\infty}\left(\text{Met}_{mw}^{q}\left(M\right)\right)$
contains the constant functions.
\end{remark}

By an argument similar to that in Proposition~\ref{sub:Propv_FVectorfield},
we have the following proposition:

\begin{proposition}
For every $F\in C_{\int}^{\infty}\left(\emph{Met}_{mw}^{q}\left(M\right)\right)$,
the local assignments given by $v_{F}\left(g\right):=\sum_{j=1}^{n}\frac{\partial A}{\partial y_{j}}\cdot X_{F_{r_{j}}}\left(g\right)$
on each neighbourhood $U$ define a unique vector field $v_{F}$
on $\emph{Met}_{mw}^{q}\left(M\right)$ satisfying $dF=\imath\left(v_{F}\right)\Omega$.
\end{proposition}

It follows that $\left\{ F,G\right\} :=-\Omega\left(v_{F},v_{G}\right)$
def\/ines a Poisson bracket on $\text{Met}_{mw}^{q}\left(M\right)$.

\subsection*{Acknowledgements}

I thank Eckhard Meinrenken for suggesting this project and Yael Karshon
who joint supervised this work as part of the author's PhD thesis.
I also thank Boris Khesin for his many helpful suggestions towards
improving both the content and exposition of this paper.

\pdfbookmark[1]{References}{ref}
\LastPageEnding

\end{document}